\documentclass{amsart}
\usepackage[english]{babel} 
\usepackage{amsmath,amssymb, amsthm}
\usepackage{graphicx}
\usepackage[all]{xy}
\usepackage{hyperref}

\usepackage[OT2,T1]{fontenc}
\DeclareSymbolFont{cyrletters}{OT2}{wncyr}{m}{n}
\DeclareMathSymbol{\sha}{\mathalpha}{cyrletters}{"58}

\DeclareMathOperator{\tor}{_{tors}}
\DeclareMathOperator{\free}{_{free}}

\newcommand{\mub}{\boldsymbol\mu}
\newcommand{\unr}{\textnormal{nr}}
\def\mod{\textnormal{ mod }}
\def\sel{\textnormal{Sel}}
\def\emb{\textnormal{em}}
\def\pr{\textnormal{pr}}

\def\EP{\hspace*{\fill}$\square$

\vspace{1.0ex}

}

\def\spec{\textnormal{Spec }}
\def\proof{\noindent {\it Proof. }} 
\DeclareMathOperator{\coker}{\textnormal{coker}}
\def\div{\textnormal{div}}
\def\res{\textnormal{Res}}
\def\Res{\textnormal{Res}}

\def\Gal{\textnormal{Gal}}
\def\phi{\varphi}
\def\l{\ell}

\def\m{\mathfrak{m}}
\newcommand{\p}{\mathfrak{p}}

\def\A{{\mathcal A}}
\def\B{{\mathcal B}}

\def\E{{\mathcal E}}

\def\L{{\mathcal L}}

\def\O{{\mathcal O}}

\def\CC{{\mathbb C}}

\def\FF{{\mathbb F}}
\def\GG{{\mathbb G}}

\def\QQ{{\mathbb Q}}
\def\RR{{\mathbb R}}
\def\ZZ{{\mathbb Z}}

\newtheorem{thm}{Theorem}[section]
\newtheorem{con}[thm]{Conjecture}
\newtheorem{que}[thm]{Question}
\newtheorem{cor}[thm]{Corollary}
\newtheorem{pro}[thm]{Proposition}
\newtheorem{lem}[thm]{Lemma}

\theoremstyle{definition}

\newtheorem{Set}[thm]{Setting}
\newtheorem{Rem}[thm]{Remark}
\newtheorem{Exa}[thm]{Example}
\newtheorem{Not}[thm]{Notation}

\begin{document}

\title[Examples of abelian surfaces with non-square $\sha$]{Examples of non-simple abelian surfaces over the rationals with non-square order Tate-Shafarevich group}
\author{Stefan Keil}
\thanks{The author is supported by a scholarship from the Berlin Mathematical School (BMS)}
\address{Humboldt-Universit\"at zu Berlin, Institut f\"ur Mathematik \\ Unter den Linden 6, 10099 Berlin, Germany}
\email{keil@math.hu-berlin.de}

\begin{abstract}
Let $A$ be an abelian surface over a fixed number field. If $A$ is principally
polarised, then it is known that the order of the Tate-Shafarevich group of $A$
must, if finite, be a square or twice a square. The situation for $A$ not
principally polarised remains unclear. 
For each $k \in \{1,2,3,5,6,7,10,13\}$ we construct a non-simple non-principally polarised abelian surface $B/\QQ$ whose Tate-Shafarevich group has order $k$ times a square. 
To obtain this result, we explore the invariance under isogeny of the Birch
and Swinnerton-Dyer conjecture.
\end{abstract}


\maketitle

\sloppy


\setcounter{tocdepth}{2}
\tableofcontents

\section{Introduction}
Let $A/K$ be an abelian variety over a number field $K$. Consider its Tate-Shafarevich group $\sha(A/K)$. If $A$ is an elliptic curve $E$, then the order of $\sha(E/K)$ is a perfect square, if it is finite. But in higher dimensions, even for principally polarised abelian varieties, this is no longer true in general. Denote by $A^\vee$ the dual abelian variety. The Cassels-Tate pairing \cite {Cassels-IV}, \cite{Tate_duality_thms} 
$$\langle \cdot , \cdot \rangle : \sha(A/K) \times \sha(A^\vee/K) \rightarrow \QQ/\ZZ,$$
which is non-degenerate in case $\sha(A/K)$ is finite, combined with a result of Flach \cite{flach}, gives a strong restriction on the non-square part of the order of the Tate-Shafarevich group \cite[Theorem 1.2]{Stein}.

\begin{thm}[Tate, Flach] \label{thm:nec-cond-on-non-square-part}
Assume $\sha(A/K)$ is finite. If an odd prime $p$ divides the non-square part of $\# \sha(A/K)$, then $p$ divides the degree of every polarisation of $A/K$.
\end{thm}

\begin{cor}[Poonen, Stoll]
If $A/K$ is a principally polarised abelian variety, then
$$\# \sha(A/K)=  \square \textnormal{ or } 2 \cdot \square.$$
\end{cor}

More precisely, assuming the finiteness of $\sha(A/K)$, Poonen and Stoll \cite{PS} associated to each principal polarisation $\lambda$ of $A/K$ a canonical element $c \in \sha(A/K)[2]$, and showed that the order of $\sha(A/K)$ is a square if and only if $\langle c,\lambda(c) \rangle =0$. This is clearly the case if $c=0$. They showed that $c=0$ is equivalent to the induced pairing on $\sha(A/K)$ being alternating and also equivalent to the polarisation $\lambda$ arising from a $K$-rational divisor.
It was already known by Tate \cite{Tate_duality_thms} that the order of $\sha$ is a square, if such a $K$-rational divisor exists. In case $c\neq 0$ the induced pairing is anti-symmetric, due to Flach \cite{flach}.

In 1996, Stoll constructed the first example of an abelian variety having $\# \sha=  2 \cdot \square$; see \cite{treurfeest} for some historical remarks. 
His example was the Jacobian of a genus $2$ curve over $\QQ$, a principally polarised abelian surface over $\QQ$. Thereafter, for every odd prime $p<25000$, $p\neq 37$, William Stein \cite{Stein} constructed an abelian variety $A_p/\QQ$ of dimension $p-1$, such that $\#\sha(A_p/\QQ)=p \cdot \square$. This result led Stein to make the following conjecture.

\begin{con} [William Stein]
As one ranges over all abelian varieties $A/\QQ$, every square-free natural number appears as the non-square part of the order of some $\sha(A/\QQ)$.
\end{con}

The following question is then natural.

\begin{que}
What are the possible non-square parts of the orders of finite Tate-Shafarevich groups for abelian varieties of fixed dimension over a fixed number field? Is this a finite list?
\end{que}

So far, in the case of abelian surfaces $B/\QQ$, the only known square-free positive integers $k$ which equal the non-square part of the order of some $\sha(B/\QQ)$ are $1$, $2$, and $3$. The purpose of this paper is to extend this list by $5$, $6$, $7$, $10$, and $13$. The construction we use is an isogeny applied to a product of two elliptic curves, and hence is different from the construction used by Poonen and Stoll, and by Stein. To understand the image of this isogeny we will explore an equation of Cassels and Tate, which is a consequence of the isogeny invariance of the Birch and Swinnerton-Dyer conjecture. The non-square part of the left hand side of this equation will be equal to the non-square part of the order of the Tate-Shafarevich group in question. We will explain how to calculate the right hand side and then we will give explicit examples to prove the following

\begin{thm}\label{thm:main-theorem}
For each $k \in \{1,2,3,5,6,7,10,13\}$ there exists a non-simple non-principally polarised abelian surface $B/\QQ$ such that $\# \sha(B/\QQ)=k \cdot \square$.
\end{thm}

The outline of this paper is the following. In the rest of this section we fix notation. In Section \ref{ch:controlling} we present the aforementioned equation of Cassels and Tate. This equation will break into two parts -- a local part and a global part. The remaining part of Section \ref{ch:controlling} is devoted to explaining the local part and to introduce non-simple abelian surfaces. In Section \ref{ch:examples} we will work with elliptic curves possessing a $\QQ$-rational $N$-torsion point. Such curves lead to two parameter families of abelian surfaces and we prove how to calculate the local and global part of the Cassels-Tate equation for these families. Finally, we present explicit calculations and give examples to prove the above theorem. 

\begin{Not}
Let $A/K$ be an abelian variety $A$ over a field $K$, i.e. a proper group scheme of positive dimension which is geometrically integral and of finite type over $\spec K$. Usually, $K$ is a number field, or a ($p$-adic) local field, or a finite field. Since all fields considered are perfect we do not pay attention to separability, and with $\overline K$ we denote a once and for all fixed algebraic closure of $K$. For a field $L$ containing $K$, the group of $L$-rational points is denoted by~$A(L)$, with $\O \in A(L)$ being the identity element of the group law. 
The dual abelian variety of $A/K$ is denoted by $A^\vee:= \textnormal{Pic}^0_{A/K}$ and a polarisation of $A/K$ is a symmetric isogeny $\lambda: A \rightarrow A^\vee$, such that over $\overline K$ we have $\lambda=\lambda_\L$, for an ample line bundle $\L$ on $A/\overline K$. If $\phi:A \rightarrow B$ is an isogeny between abelian varieties over a field $K$, then for a field extension $L/K$ we say that $\phi$ has a {\em $L$-kernel}, if all points in $A(\overline K)[\phi]$ are already defined over $L$, i.e. $A(\overline K)[\phi]=A(L)[\phi]$. If we do not specify the field of definition of an isogeny $\phi$ between two abelian varieties which are defined over a field~$K$, then we want $\phi$ to be also defined over $K$. 

If $K$ is a number field, then with $v$ we denote a place of $K$, and with $M_K$ the set of all places of $K$. We have the subset $M_K^0$ of all finite places (or primes) of $K$ and the subset $M_K^ \infty$ of all infinite places of $K$. With $K_v$ we denote the completion of $K$ at $v$, and with $k_v$ its residue field, i.e. the quotient of the valuation ring $\O_v$ by its maximal ideal $\m_v=\pi_v \O_v$, for a uniformiser $\pi_v$. We normalise the absolute value $|\cdot |_v$ on $K_v$ so that $|\pi_v|_v=(\#k_v)^{-1}$. If $v \in M_K^0$ is a place lying over $p \in M_\QQ^0$, we denote this by $v|p$ and call $K_v$ a $p$-adic field. 
Denote by $K_v^\unr$ the maximal unramified extension of $K_v$. It is obtained by adjoining to $K_v$ all $n$-th roots of unity, for $n$ coprime to the characteristic of $k_v$. 

The absolute Galois group of a field $K$ will be denoted by $\Gal_K$. For Galois cohomology we use the usual abbreviation $H^i(K,M):=H^i(\Gal_K,M)$, for a $K$-Galois module $M$. The Tate-Shafarevich group of $A/K$ is defined as 
$$\sha(A/K):= \ker \left( H^1(K, A(\overline K)) \rightarrow \prod_{v \in M_K} H^1( K_v, A(\overline K_v)) \right).$$

With $\l$ we denote a prime number and by $\ZZ/\l\ZZ$ we either mean a cyclic group of order $\l$ or a Galois module of order $\l$ with trivial Galois action. By $\mub_\l$ we denote the $\l$-th roots of unity as a Galois module of order $\l$, and we write $\xi=\xi_\l$ for a primitive $\l$-th root of unity. The trivial group is denoted by $0$. By $\square \in \{1,4,9,16,\ldots\}$, we denote a square natural number. 
We sometimes refer to computations carried out with the software package Sage \cite{sage}.
\end{Not}

\section{Controlling the order of Tate-Shafarevich groups modulo squares}
\label{ch:controlling}

We want to construct abelian surfaces $B/\QQ$ such that the order of their Tate-Shafarevich groups is not a square. To archieve this objective we start with an abelian surface $A/\QQ$ being the product of two elliptic curves $E_1$ and $E_2$ over $\QQ$. Hence it is known that the order of the Tate-Shafarevich group $\sha(A/\QQ)$ is a square, provided it is finite. Then we consider cyclic isogenies $\phi: A \rightarrow B$ and the goal is to understand $\sha(B/\QQ)$ in terms of $\sha(A/\QQ)$ and $\phi$. The precise situation we consider is summarised in Setting \ref{set:setting}. The isogeny $\phi$ naturally induces a group homomorphism between $\sha(A/\QQ)$ and $\sha(B/\QQ)$ which is an isomorphism between $\l$-primary parts for primes $\l$ not dividing the degree of the isogeny. 
In particular this means, that if $\sha(A/K)$ is of square order, then a necessary condition for a prime $\l$ to divide the non-square part of the order of $\sha(B/K)$ is that $\l$ divides the degree of $\phi$. 

In the next subsection we will present the Cassels-Tate equation, which expresses the relative change of the orders of the Tate-Shafarevich groups of $A$ and $B$ under $\phi$. It consists of a local and a global part, which will be called the {\em local quotient} and the {\em global quotient}. We will spend the following two subsections computing the local quotient, one for general abelian varieties and one for elliptic curves. The main result is Theorem \ref{thm:iota-p-1-of-E}. 
The isogenies constructed in Setting \ref{set:setting} are then subsumed into the so-called {\em isogenies with diagonal kernel} in the last subsection of this chapter. We also introduce the general concept of non-simple abelian surfaces and we present the Key Lemma \ref{lem:hitting-both-H-1-general-version} for the computation of the local quotient. In Section \ref{ch:examples} we will use the results obtained in this section to actually compute explicit examples of abelian surfaces $B/\QQ$ having Tate-Shafarevich group of non-square order.

\subsection{An equation of Cassels and Tate}\label{sec:c-t-equation}

In the mid 1960s, Cassels \cite{Cassels} (the elliptic curve case) and Tate \cite{Tate} (the general case) proved the following theorem to show the invariance of the Birch and Swinnerton-Dyer conjecture under isogeny.
Denote by $R_A$ the regulator and by $P_A$ the period of an abelian variety $A/K$ over a number field $K$, see pages $37$ and $52$ in \cite{jorza} for the definitions. By $c_{A,v}$ we denote the local Tamagawa number of $A$ at a finite place $v \in M_K^0$. 

\begin{thm} \label{thm:cassels-tate-equation}
Let $\phi:A\rightarrow B$ be an isogeny between two abelian varieties $A$ and $B$ over a number field $K$. Assume that at least one of $\sha(A/K)$ or $\sha(B/K)$ is finite. Then $\sha(A/K)$ and $\sha(B/K)$ are both finite, and
$$ \frac{\# \sha(A/K)}{\# \sha(B/K)} = \frac{R_B}{R_A} \cdot \frac{\# A(K)\tor \# A^\vee (K)\tor}{\# B(K)\tor \# B^\vee (K)\tor} \cdot \frac{P_B}{P_A} \cdot \prod_{v \in M_K^0}\frac{c_{B,v}}{c_{A,v}}.$$
\end{thm}

The product over the Tamagawa numbers is actually finite, since $c_{A,v}=1$ when $v$ is a place of good reduction of $A$. We define $A(K)\free$ to be the quotient group $A(K)/A(K)\tor$. Consider the following induced group homomorphisms.
$$\phi_K:A(K) \rightarrow B(K), \ \phi^\vee_K: B^\vee(K) \rightarrow A^\vee(K), \ \phi_v:A(K_v) \rightarrow B(K_v),$$
$$\phi_{K,\textnormal{tors}} :A(K)\tor \rightarrow B(K)\tor, \ \ \phi_{K,\textnormal{tors}}^\vee :B^\vee(K)\tor \rightarrow A^\vee(K)\tor,$$
$$\phi_{K,\textnormal{free}} :A(K)\free \rightarrow B(K)\free, \ \ \phi_{K,\textnormal{free}}^\vee :B^\vee(K)\free \rightarrow A^\vee(K)\free.$$
We may now reformulate the above quotients in terms of these induced group homomorphisms. This reformulation, which is part of the proof of the above theorem, turns out to be easier to handle for computational purposes, and we are going to use the Cassels-Tate equation only in this description. There are two trivial equalities, namely
$$\frac{\# A(K)\tor}{\# B(K)\tor} = \frac{\# \ker \phi_K}{\# \coker \phi_{K,\textnormal{tors}}} \ \mbox{   and   } \ \frac{\# A^\vee (K)\tor}{\# B^\vee (K)\tor} = \frac{\# \coker \phi^\vee_{K,\textnormal{tors}}}{\# \ker \phi^\vee_K},$$
and two more interesting ones, see the proof of \cite[Theorem I.7.3]{milne-arith-dual-thm}; 
$$\frac{R_B}{R_A}= \frac{\# \coker \phi^\vee_{K,\textnormal{free}}}{\# \coker \phi_{K,\textnormal{free}}} \ \mbox{   and   } \
\frac{P_B}{P_A} \cdot \prod_{v \in M_K^0}\frac{c_{B,v}}{c_{A,v}} = \prod_{v \in M_K} \frac{\# \coker \phi_v}{\# \ker \phi_v}.$$
Hence the Cassels-Tate equation becomes
\begin{equation}\label{equ:Cassels-Tate-equation}
\frac{\# \sha(A/K)}{\# \sha(B/K)} = \frac{\# \ker \phi_K}{\# \coker \phi_K}\frac{\# \coker \phi^\vee_K}{\# \ker \phi^\vee_K} \prod_{v \in M_K} \frac{\# \coker \phi_v}{\# \ker \phi_v}.
\end{equation}
In particular we have
$$\frac{R_B}{R_A} \cdot \frac{\# A(K)\tor \# A^\vee (K)\tor}{\# B(K)\tor \# B^\vee (K)\tor} = \frac{\# \ker \phi_K}{\# \coker \phi_K}\frac{\# \coker \phi^\vee_K}{\# \ker \phi^\vee_K},$$ 
and we call the right-hand side of this equation the {\em global quotient}. The global quotient clearly breaks into the {\em regulator quotient} and the {\em torsion quotient}. The product 
$$\prod_{v \in M_K} \frac{\# \coker \phi_v }{ \# \ker \phi_v}$$
runs over all places $v$ of $K$ and is called the {\em local quotient}. It is in fact a finite product, since $\# \coker \phi_v = \# \ker \phi_v$ for all but finitely many places $v$, as will be recalled in Corollary \ref{cor:product-over-quotient-of-phi_v-is-finite}. In the next two subsections we will study the local quotient $\# \coker \phi_v / \# \ker \phi_v$ for a finite place $v \in M_K^0$.

\subsection{Isogenies between abelian varieties over local fields} 
\label{sec:isogenies-over-local-fields}

In this section we will use the following notation. Let $\phi:A\rightarrow B$ be an isogeny between two abelian varieties $A$ and $B$ over a number field $K$, and let $v \in M_K^0$ be a finite place of $K$ lying over a fixed prime $p$. Consider the induced group homomorphism on $K_v$-rational points
$$\phi_v:A(K_v) \rightarrow B(K_v).$$
Our aim is to compute the quotient $\# \coker \phi_v / \# \ker \phi_v$, which mainly consists in determining the cardinality of $\coker \phi_v$, as the size of the kernel is usually obvious by the definition of $\phi$. 
On a few occasions we will focus on isogenies having a $K_v$-kernel, i.e. $ A(\overline K_v)[\phi]=A(K_v)[\phi]$, and thus $\# \ker \phi_v=\deg \phi$ and $\Gal_{K_v}$ acts trivially on $A(\overline K_v)[\phi]$. 

In general, the cokernel of $\phi_v$ can naturally be identified with a subgroup of $H^1(K_v, A(\overline K_v)[\phi])$, since the short exact sequence of Galois modules
$$\xymatrix{0 \ar[r] & A(\overline K_v)[\phi] \ar[r] & A(\overline K_v) \ar[r]^{\phi} & B(\overline K_v) \ar[r] & 0}$$
gives the long exact Galois cohomology sequence
$$\xymatrix{0 \ar[r] & \coker \phi_v \ar[r] & H^1(K_v, A(\overline K_v)[\phi]) \ar[r] & \cdots}$$

The next lemma determines the size of $H^1(K_v, A(\overline K_v)[\phi])$ and in particular shows that it is finite. Hence $\coker \phi_v$ is also finite.

\begin{lem}\label{lem:order-of-H-1-by-serre}
Let $K_v$ be a $p$-adic field and let $M$ be a finite $K_v$-Galois module of order $\# M$ and with dual $M^\vee := \textnormal{Hom}(M,\GG_m(\overline K_v))$. The size of the first cohomology group of $M$ can be computed as follows.
$$\# H^1(K_v,M) = \# H^0(K_v,M) \cdot \# H^0(K_v,M^\vee) \cdot p^{v_p(\# M) \cdot [K_v:\QQ_p]}.$$
\end{lem}
\proof
This follows from Theorems 2 and 5 in Chapter II.5 in \cite{serre-gal-coh}. Define the Euler-Poincar\'e characteristic by $\chi(K_v,M):= \# H^0(K_v,M) \cdot \# H^2(K_v,M) / \# H^1(K_v,M)$. By the duality Theorem 2 from \cite{serre-gal-coh}, we get $\# H^2(K_v,M) = \# H^0(K_v,M^\vee)$, and by Theorem 5, we get $\chi(K_v,M) = (\O_v : \# M \O_v)^{-1}$, where $\O_v$ is the valuation ring of $K_v$. Hence, $\chi(K_v,M) = p^{-v_p(\# M)\cdot [K_v:\QQ_p]}$ and we are done. \EP

\begin{cor} \label{cor:structure-of-H-1} Let $\phi$ be of prime degree $\l$. If $\phi$ or $\phi^\vee$ has a $K_v$-kernel, then
$$H^1(K_v,A(\overline K_v)[\phi])\cong \begin{cases}
    \ZZ/\l\ZZ, & v\nmid \l, \mub_\l \nsubseteq K_v\\
(\ZZ/\l\ZZ)^2, & v\nmid \l, \mub_\l \subseteq K_v\\
(\ZZ/\l\ZZ)^{[K_v:\QQ_p]+1}, & v|\l, \mub_\l \nsubseteq K_v\\
(\ZZ/\l\ZZ)^{[K_v:\QQ_p]+2}, & v|\l, \mub_\l \subseteq K_v,
                            \end{cases}$$
and if neither $\phi$ nor $\phi^\vee$ has a $K_v$-kernel, then
$$H^1(K_v, A(\overline K_v)[\phi]) \cong \begin{cases}
                                0, & v \nmid \l\\
				(\ZZ/\l\ZZ)^{[K_v:\QQ_p]}, & v |\l.
                                \end{cases}$$
\end{cor}
\proof By definition $H^1(K_v,A(\overline K_v)[\phi])$ is abelian and has exponent $\l$. By the previous lemma, for $M:=A(\overline K_v)[\phi]$, we have
$$\# H^1(K_v,M) = \begin{cases} 
                  \# H^0(K_v,M) \cdot \# H^0(K_v,M^\vee), & v\nmid \l \\ 
                  \# H^0(K_v,M) \cdot \# H^0(K_v,M^\vee) \cdot \l^{[K_v:\QQ_p]}, & v \mid \l.
                  \end{cases}$$
If $\phi$, respectively $\phi^\vee$, has a $K_v$-kernel, then $A(\overline K_v)[\phi] \cong \ZZ/l\ZZ$, respectively $\mub_\l$, as Galois modules. Since
$$H^0(K_v,\ZZ/\l\ZZ)\cong \ZZ/\l\ZZ, \text{ and } H^0(K_v,\mub_\l) \cong \begin{cases}
                            0, & \mub_\l \nsubseteq K_v\\
								    \ZZ/\l\ZZ, & \mub_\l \subseteq K_v,                                        \end{cases}$$
and $\ZZ/\l\ZZ$ and $\mub_\l$ are dual to each other, we get the first statement.

If neither $\phi$ nor $\phi^\vee$ has a $K_v$-kernel, then neither $A(\overline K_v)[\phi]$ nor its dual is isomorphic to $\ZZ/\l\ZZ$. Therefore
$$H^0(K_v,A(\overline K_v)[\phi])= H^0(K_v,A(\overline K_v)[\phi]^\vee) =0,$$
which completes the proof. \EP

\begin{cor}\label{cor:H1}
$$H^1(\QQ_p,\ZZ/\l\ZZ)\cong H^1(\QQ_p,\mub_\l)\cong\begin{cases}
    \ZZ/\l\ZZ, & p\neq \l \neq 2, p\not \equiv 1 \mod \l\\
(\ZZ/\l\ZZ)^3, & p=\l=2\\
(\ZZ/\l\ZZ)^2, & otherwise.\\
                            \end{cases}$$
\end{cor}
\proof This is immediate from Corollary \ref{cor:structure-of-H-1} upon observing that $\mub_2 \subseteq \QQ_p$ for all $p$, and $\mub_\l \nsubseteq \QQ_p$ if and only if $p \not \equiv 1 \mod \l$ and $\l\neq 2$. \EP

For a finite $K_v$-module $M$ we will now introduce the unramified Galois cohomolgy group which is an important subgroup of $H^1(K_v, M)$. Denote by $K_v^\unr$ the maximal unramified extension of $K_v$. We have that the inertia group $I_v:=\Gal_{K_v^\unr}$ is a normal subgroup of $\Gal_{K_v}$; thus the usual restriction homomorphism
$$\res_\unr : H^1(K_v, M) \rightarrow H^1(K_v^\unr, M)$$
is defined and by the Inflation-Restriction sequence its kernel is isomorphic to $H^1(\Gal(K_v^\unr / K_v), M^{I_v})$. We denote the kernel of $\res_\unr$ by $H^1_\unr(K_v, M)$ and call it the {\em unramified subgroup} of $H^1(K_v, M)$.
Consider again the following Galois cohomology sequence with respect to an isogeny $\phi:A \rightarrow B$.
$$\xymatrix{0 \ar[r] & \coker \phi_v \ar[r]^/-1em/{\delta_v} & H^1(K_v, A(\overline K_v)[\phi]) \ar[r] & \cdots}.$$
We say that $\coker \phi_v$ is {\em maximal}, respectively {\em maximally unramified}, respectively {\em trivial}, if $\delta_v$ is an isomorphism, respectively if $\delta_v$ induces an isomorphism between $\coker \phi_v$ and the unramified subgroup $H^1_\unr(K_v, A(\overline K_v)[\phi])$, respectively if $\coker \phi_v$ is the trivial group. 

\begin{Rem}\label{rem:unramified-is-good-as-local-quotient-equals-1}
If $K=\QQ$ and $(p,\l) \neq (2,2)$, the last two corollaries show that if the isogeny $\phi:A \rightarrow B$ is of prime degree $\l$ and has a $\QQ_p$-kernel, then $H^1(\QQ_p, A(\overline \QQ_p)[\phi])$ is either isomorphic to $\ZZ/\l\ZZ$ or $(\ZZ/\l\ZZ)^2$. In the former case, either $\coker \phi_p$ is the trivial group, or is isomorphic to $H^1(\QQ_p, A(\overline \QQ_p)[\phi])$. In the latter case there is a third possibility, namely that $\coker \phi_p$ has $\l$ elements and thus is one of the $\l+1$ subgroups of $H^1(\QQ_p, A(\overline \QQ_p)[\phi])$ of order $\l$. 
By the next lemma, the unramified subgroup is one of these $\l+1$ subgroups of order $\l$.

Besides merely determing the size of $\coker \phi_v$ our goal is further to specify it as a subgroup of $H^1(K_v, A(\overline K_v)[\phi])$, and hence the main purpose of this subsection is to give criteria to check whether $\coker \phi_v$ is maximally unramified.
\end{Rem}

\begin{lem} \label{lem:H1-unr}
Let $K_v$ be a $p$-adic field and let $M$ be a finite $K_v$-module. Then the order of $H^1_\unr(K_v, M)$ equals the order of $H^0(K_v, M)$. 
\end{lem}
\proof For a prime $\l$ denote by $M[\l^\infty]$ the $\l$-primary part of $M$, thus $M=\oplus_\l M[\l^\infty]$. As Galois acts on the individual $M[\l^\infty]$, we get $H^0(K_v, \oplus_\l M[\l^\infty])=\oplus_\l H^0(K_v, M[\l^\infty])$ and $H^1_\unr(K_v, \oplus_\l M[\l^\infty])=\oplus_\l H^1_\unr(K_v, M[\l^\infty])$. Now apply Lemma 4.2 in \cite{schaefer-stoll_long-version} to get that the order of $H^1_\unr(K_v, M[\l^\infty])$ equals the order of $H^0(K_v, M[\l^\infty])$. \EP

We introduce some more notation. By $\tilde \A$ we denote the reduction of $A$ modulo $v$, i.e. the special fiber at $v$ of the N\'{e}ron model $\A/\O_K$ of $A$, and by $\tilde \A_0(k_v)$ we denote the smooth part of the $k_v$-rational points of the reduction at $v$, i.e. the $k_v$-rational points of the connected component of $\tilde \A$ intersecting the zero-section. Denote by $A_0(K_v)$ the preimage of $\tilde \A_0(k_v)$ under the reduction-mod-$v$ map, and by $A_1(K_v)$ the kernel of $A_0(K_v) \rightarrow \tilde \A_0(k_v)$. We have the following two commutative diagrams with exact rows and induced group homomorphisms as vertical arrows.

\begin{equation}\label{equ:reduction_maps}
\begin{split}
  \xymatrix{
      0 \ar[r] & A_1(K_v) \ar[r]  \ar[d]_{\phi_v^1} &   A_0(K_v) \ar[r]  \ar[d]_{\phi_v^0} & \tilde \A_0(k_v) \ar[r] \ar[d]_{\tilde \phi_v^0} & 0\\
      0 \ar[r] & B_1(K_v) \ar[r]  & B_0(K_v) \ar[r]  & \tilde \B_0(k_v) \ar[r]  & 0\\
 }
\end{split}
\end{equation}
\begin{equation}\label{equ:tamagawa-number}
\begin{split}
  \xymatrix{
      0 \ar[r] & A_0(K_v) \ar[r]  \ar[d]_{\phi_v^0} &   A(K_v) \ar[r]  \ar[d]_{\phi_v} & A(K_v)/ A_0(K_v)  \ar[r] \ar[d]_{\overline \phi_v} & 0\\
      0 \ar[r] & B_0(K_v) \ar[r]  & B(K_v) \ar[r] & B(K_v)/ B_0(K_v)  \ar[r] & 0\\
 }
\end{split}
\end{equation}

The vertical maps on the right, i.e. $\tilde \phi_v^0$ and $\overline \phi_v$, are group homorphisms between finite groups, which follows from the theory of N\'eron models. Therefore, the kernels and cokernels of $\tilde \phi_v^0$ and $\overline \phi_v$ are finite groups. The kernels of $\phi_v^0$ and $\phi_v^1$ are finite as they are subgroups of $\ker \phi_v$, which is finite by definition. The cokernels of $\phi_v^0$ and $\phi_v^1$ are finite by the snake lemma, since $\coker \phi_v$ is, as seen in Corollary \ref{cor:structure-of-H-1}. Hence, all kernels and cokernels of the vertical maps in the above two diagrams are finite groups. 
In the unramified case we get the following commutative diagram with exact rows.
\begin{equation}\label{equ:unramified}
\begin{split}
  \xymatrix{
      0 \ar[r] & A_1(K_v^\unr) \ar[r]  \ar[d]_{\phi_{v,\unr}^1} &   A_0(K_v^\unr) \ar[r]  \ar[d]_{\phi_{v,\unr}^0} & \tilde \A_0(\overline k_v) \ar[r] \ar[d]_{\tilde \phi^0_{\overline k_v}} & 0\\
      0 \ar[r] & B_1(K_v^\unr) \ar[r]  & B_0(K_v^\unr) \ar[r]  & \tilde \B_0(\overline k_v) \ar[r]  & 0\\
 }
\end{split}
\end{equation}

We recall a basic fact, which follows from Lang's Theorem \cite[Theorem 1]{Lang}.

\begin{lem}\label{lem:number-of-reduced-points-is-the-same}
With notation as above, $\tilde \A_0(k_v)$ and $\tilde B_0(k_v)$ are finite groups of same cardinality.
\end{lem}
\proof The proof is given on page 561 in \cite{Lang}.
From the theory of N\'eron models it follows that $\tilde \A_0$ and $\tilde B_0$ are isogenous connected algebraic groups over the finite field $k_v$. 
Let $G/k$ be a connected algebraic group over the finite field $k$ of size $q$. Denote the group law by multiplication and define the Lang isogeny $f_G(g):=g^{-1}g^{(q)}$, for $g\in G(\overline k)$, where $g^{(q)}$ is the image of $g$ under the Frobenius morphism. Lang's Theorem \cite[Corollary of Theorem 1]{Lang} says that $f_G:G(\overline k) \rightarrow G(\overline k)$ is indeed an isogeny with kernel equal to the $k$-rational points of $G$. Now let $\phi:G \rightarrow H$ be an isogeny between connected algebraic groups $G$ and $H$ over $k$. Then $f_H \circ \phi = \phi \circ f_G$, and hence the kernels of $f_G$ and $f_H$ have the same cardinality, which proves the lemma. \EP

Now we apply the snake lemma on Diagrams (\ref{equ:reduction_maps}) and (\ref{equ:tamagawa-number}) to get a basic lemma. Recall, that the quantity $c_{A,v}$ is defined as the order of the quotient group $A(K_v)/ A_0(K_v)$ and is called the {\em local Tamagawa number} of $A$ at $v$.

\begin{lem}\label{lem:isogenies-on-k_v-points}
With notation as above, we have the equality
$$\frac{\# \coker \phi_v}{\# \ker \phi_v} = \frac{\# \coker \phi_v^1}{\# \ker \phi_v^1} \cdot \frac{c_{B,v}}{c_{A,v}}.$$
\end{lem}
\proof We have already seen the finiteness of all appearing kernels and cokernels. 
Applying the snake lemma on the kernels and cokernels in Diagram (\ref{equ:reduction_maps}) we get
$$\frac{\# \ker \phi_v^1 }{\# \coker \phi_v^1} \cdot \frac{\# \ker \tilde \phi_v^0}{\# \coker \tilde \phi_v^0} = \frac {\# \ker \phi_v^0}{\# \coker \phi_v^0}.$$
Since $\#\tilde \A_0(k_v) = \#\tilde B_0(k_v)$ by Lemma \ref{lem:number-of-reduced-points-is-the-same},
we get $\# \ker \tilde \phi_v^0 = \# \coker \tilde \phi_v^0$. Hence 
$$\frac{\# \ker \phi_v^1 }{\# \coker \phi_v^1} = \frac {\# \ker \phi_v^0}{\# \coker \phi_v^0}.$$
Applying the snake lemma on Diagram (\ref{equ:tamagawa-number}) gives
$$\frac{\# \coker \phi_v}{\# \ker \phi_v} = \frac{\# \coker \phi_v^0}{\# \ker \phi_v^0} \cdot \frac{\# \coker \overline \phi_v}{\# \ker \overline \phi_v}.$$
By definition we have
$$\frac{\# \coker \overline \phi_v}{\# \ker \overline \phi_v}=\frac{c_{B,v}}{c_{A,v}},$$
which completes the proof. \EP

We continue by examining the quotient $\# \coker \phi_v^1 / \# \ker \phi_v^1$. We start by recalling two basic lemmas, and then we deduce the well known fact that this quotient is almost always trivial, since $\phi_v^1$ will be an isomorphism for all but finitely many places $v$.  

\begin{lem}\label{lem:ker-of-red-is-a-pro-p-group}
The kernel of reduction $A_1(K_v)$ is a pro-$p$ group.
\end{lem}
\proof The multiplication-by-$\l$ endomorphism $[\l]$ on $A_1(K_v)$ is an isomorphism, for all primes $\l$ different to the characteristic $p$ of the residue field $k_v$, as $A_1(K_v)$ is isomorphic to the group $\hat A(\m_v)$ associated to the formal group $\hat A$ of $A$ defined over the valuation ring $\O_v$ of $K_v$ with maximal ideal~$\m_v$.
Hence for any subgroup~$G$ of~$A_1(K_v)$, $[\l]$ is a surjective endomorphism on $A_1(K_v)/G$, for $\l\neq p$. If in addition~$A_1(K_v)/G$ is finite, then $[\l]$ is an automorphism on $A_1(K_v)/G$, for all $\l \neq p$, and thus $A_1(K_v)/G$ is a $p$-group. Hence, $A_1(K_v)$ is a pro-$p$ group. \EP



\begin{lem}\label{lem:sufficient-condition-that-phi_v-1-is-an-iso} 
If $v \nmid \deg \phi$, then $\phi_v^1$ and $\phi_{v,\unr}^1$ are isomorphisms.
\end{lem}
\proof Denote the degree of $\phi$ by $n$. There exist isogenies $\psi:B\rightarrow A$ and $\rotatebox[origin=c]{180}{$\psi$}: A \rightarrow B$, such that $\psi \circ \phi:A\rightarrow A$ and $\rotatebox[origin=c]{180}{$\psi$} \circ \psi: B \rightarrow B$ are the multiplication-by-$n$ maps $[n]$. Hence we get the following induced group homomorphisms on the kernels of reduction.
$$\begin{xy}
  \xymatrix{
      A_1(K_v) \ar[r]^{\phi_v^1}  \ar@/^ 0.7cm/[rr]^{[n]_v^1} &  B_1(K_v) \ar[r]^{\psi_v^1}  \ar@/_ 0.7cm/[rr]^{[n]_v^1} & A_1(K_v) \ar[r]^{\rotatebox[origin=c]{180}{$\psi$}_v^1}  & B_1(K_v) \\
 }
\end{xy}$$
Since $v \nmid \deg \phi$, we have by the previous lemma that both maps $[n]_v^1$ are isomorphisms. Hence it follows that all three homomorphisms $\psi_v^1$, $\rotatebox[origin=c]{180}{$\psi$}_v^1$ and $\phi_v^1$ are isomorphisms. Now for any finite unramified extension $L_w/K_v$, we get by the same argument that $\phi_w^1$ is an isomorphism, and so also is $\phi_{v,\unr}^1$.\EP

\begin{cor}\label{cor:kernels-and-cokernels-and-Tamagawa-quotient-are-powers-of-l}
If a prime $\l$ divides the cardinality of a kernel or cokernel of one of the induced group homomorphisms $\phi_v$, $\phi_v^0$, $\phi_v^1$, $\overline \phi_v$ or $\tilde \phi_v^0$ appearing in Diagrams (\ref{equ:reduction_maps}) and (\ref{equ:tamagawa-number}), or $\l$ divides the Tamagawa quotient $c_{B,v}/c_{A,v}$, then $\l \mid \deg \phi$. Further, if $\gcd(\deg \phi, c_{A,v} \cdot c_{B,v})=1$, then $\overline \phi_v$ is an isomorphism.

In particular, if $\phi$ is of prime degree $\l$, then the cardinalities of all kernels and cokernels of $\phi_v$, $\phi_v^0$, $\phi_v^1$, $\overline \phi_v$ and $\tilde \phi_v^0$, as well as the Tamagawa quotient $c_{B,v}/c_{A,v}$, are powers of $\l$. 
\end{cor}
\proof By construction this is clear for all the kernels. If $\phi$ is the multiplication-by-$n$ endomorphism of $A$ for a positive integer $n$, then this is also clear for the cokernels. For a general isogeny $\phi$ of degree $n$, as mentioned in the proof of the above lemma, there is an isogeny $\psi:B \rightarrow A$, such that $[n]=\psi \circ \phi$. From the exact sequence
$$0 \rightarrow \ker \psi / \phi(\ker [n]) \rightarrow \coker \phi \rightarrow \coker [n] \rightarrow \coker \psi \rightarrow 0$$
we derive the statement about the cokernels of the homomorphisms induced by $\phi$. Use Lemma \ref{lem:isogenies-on-k_v-points} to get the statement about the Tamagawa quotient.


Now assume that $\gcd(\deg \phi, c_{A,v} \cdot c_{B,v})=1$. If a prime $\l$ divides the Tamagawa quotient $c_{B,c}/c_{A,v}$ then $\l \mid \deg \phi$ by the above part of this lemma, hence $\l$ does not divide the product $c_{A,v} \cdot c_{B,v}$. Therefore there are no primes $\l$ dividing $c_{B,c}/c_{A,v}$ and thus $c_{A,v}= c_{B,v}$. This implies $\# \ker \overline \phi_v = \# \coker \overline \phi_v$. If a prime $\l$ divides $\# \ker \overline \phi_v$, then $\l$ divides $\deg \phi$ and $c_{A,v}$, hence there are no such primes $\l$ and $\overline \phi_v$ is an isomorphism.  \EP

We conclude that the product over all quotients $\# \coker \phi_v / \# \ker \phi_v$ is actually a finite product. Denote by $M_K$ the set of places of $K$ and let $S$ be a finite set of $M_K$ containing the infinite places, the places of bad reduction and the places dividing the degree of the isogeny $\phi$.

\begin{cor}\label{cor:product-over-quotient-of-phi_v-is-finite}
If $v\nmid \deg \phi$ and $v$ is a place of good reduction, then
$$\frac{\# \coker \phi_v}{\# \ker \phi_v}=1;$$
thus
$$\prod_{v \in M_K} \frac{\# \coker \phi_v}{\# \ker \phi_v} =\prod_{v \in S} \frac{\# \coker \phi_v}{\# \ker \phi_v}.$$
\end{cor}
\proof Combine Lemmas \ref{lem:isogenies-on-k_v-points} and \ref{lem:sufficient-condition-that-phi_v-1-is-an-iso} with the fact that the Tamagawa quotient equals $1$ in case of good reduction. \EP

In view of the corollary, the goal of this subsection is to provide methods to compute the quotient $\# \coker \phi_v/ \# \ker \phi_v$, in case $v$ is a place of bad reduction or $v \mid \deg \phi$. If we stick to good reduction, but do not care whether $v$ divides the degree of $\phi$, then the next lemma gives a very nice criterion to check whether $\coker \phi_v$ is maximally unramified.
The notation used in part (i) of the lemma comes from the following diagram. 
$$\xymatrix{
      0 \ar[r] & A_1(\overline K_v) \ar[r]  \ar[d]_{\phi_{\overline K_v}^1} &   A_0(\overline K_v) \ar[r]  \ar[d]_{\phi_{\overline K_v}^0} & \tilde \A_0(\overline k_v) \ar[r] \ar[d]_{\tilde \phi_v^0} & 0\\
      0 \ar[r] & B_1(\overline K_v) \ar[r]  & B_0(\overline K_v) \ar[r]  & \tilde \B_0(\overline k_v) \ar[r]  & 0\\
 }$$

\begin{lem}[Criterion for maximal unramifiedness of $\coker \phi_v$ in case $v$ is a place of good reduction]\label{lem:good-reduction-is-unramified} 
Assume $v$ is a place of good reduction.

(i) If $\ker \phi^1_{\overline K_v}$ is trivial, then $\coker \phi_v$ is maximally unramified.

(ii) If $\phi$ has a $K_v$-kernel and $\phi_v^1$ is injective, then $\coker \phi_v$ is maximally unramified.
\end{lem}
\proof Part (ii) for $K_v=\QQ_p$, $\deg \phi$ being an odd prime, and $A$ and $B$ are elliptic curves is Lemma A.3 in the Appendix of \cite{DokchitserFisher} by Tom Fisher. In general, (ii) follows directly from (i), as the assumptions imply that $\ker \phi^1_{\overline K_v} = \ker \phi^1_v =0$.

For (i) note, that if $[\xi] \in H^1(K_v,A[\phi])$ is an element of $\coker \phi_v$, then $[\xi]$ lies in the kernel of $H^1(K_v,A[\phi]) \rightarrow H^1(K_v,A)$. This means that there is a point $P \in A(\overline K_v)$, such that $\xi(\sigma)=P^\sigma -P$, for all $\sigma \in \Gal_{K_v}$. As $v$ is a place of good reduction we get that $P \in A_0(\overline K_v)$. Consider the reduction-mod-$v$ map $A_0(\overline K_v) \rightarrow \tilde \A_0(\overline k_v)$, which is a group homomorphism. Hence, $\overline{P^\tau -P} = {\overline P}^\tau - \overline P = \O$, for all $\tau \in I_v$, as $I_v$ acts trivially on $\A_0(\overline k_v)$. Therefore for all $\tau \in I_v$, $P^\tau -P$ lies in the kernel of reduction $\phi^1_{\overline K_v}$. As $\phi^1_{\overline K_v}$ is assumed to be trivial we immediately deduce that $P^\tau -P=\O$, for all $\tau \in I_v$, which is equivalent to $P \in A_0(K_v^\unr)$. By definition, $[\xi]$ lies in $H^1_\unr (K_v,A[\phi])$ if it is in the kernel of $\Res_\unr$. This is clearly the case if $P \in A(K_v^\unr)$, because this makes the restriction of $\xi$ to $I_v$ to be the zero map, and thus $\coker \phi_v$ injects into $H^1_\unr (K_v,A[\phi])$. 

By Lemmas \ref{lem:H1-unr} and \ref{lem:isogenies-on-k_v-points}, $\coker \phi_v$ also surjects onto $H^1_\unr (K_v,A[\phi])$, as its order is at least the order of $H^1_\unr (K_v,A[\phi])$. \EP

We continue with presenting a reinterpretation of the quotient $\# \coker \phi_v^1 / \# \ker \phi_v^1$ given by Schaefer in \cite{schaefer}. Using these results it is quite easy to compute $\# \coker \phi_v^1 / \# \ker \phi_v^1$ for elliptic curves. First we need some notation. Assume that the abelian varieties $A$ and $B$ are of dimension $d$ and let $v \in M_K^0$ be a finite place. We can write the isogeny $\phi:A \rightarrow B$ as a $d$-tuple of power series in $d$-variables in a neighbourhood of the identity element $\O$. Let $|\phi'(0)|_v$ be the normalised $v$-adic absolute value of the determinant of the Jacobian matrix of partial derivatives of such a power series representation of $\phi$ evaluated at $0$. Note that $|\phi'(0)|_v$ is well defined.

\begin{pro} \label{pro:coker-phi-prime-relation}
With notation as above,
$$|\phi'(0)|_v^{-1}= \frac{\# \coker \phi_v^1}{\# \ker \phi_v^1};$$
hence
$$|\phi'(0)|_v= 1, \text{ if } v \nmid \deg \phi.$$ 
\end{pro}
\proof Combine \cite[Lemma 3.8]{schaefer} with Lemmas \ref{lem:isogenies-on-k_v-points} and \ref{lem:sufficient-condition-that-phi_v-1-is-an-iso}. \EP

\begin{cor}\label{cor:isogenies-on-k_v-points}
With notation as above,
$$\frac{\# \coker \phi_v}{\# \ker \phi_v} = |\phi'(0)|_v^{-1} \cdot \frac{c_{B,v}}{c_{A,v}}.$$
\end{cor}
\proof This follows from the last proposition together will Lemma \ref{lem:isogenies-on-k_v-points}. \EP

\begin{Rem} In the case of elliptic curves, $\phi'(0)$ is just the leading coefficent of the power series representation of $\phi$. One can easily compute this value: Use V\'elu's algorithm \cite{velu} to describe $\phi$ as coordinate functions $\phi(x,y)=(\tilde x(x,y), \tilde y(x,y))$ and then write $- \tilde x/ \tilde y$ as a power series in $z:=-x/y$, see \cite[IV]{silvermanI}. We will do this explicitly in Propositions \ref{pro:p-value-of-phi-prime-zero-is-one} and \ref{pro:p-value-of-phi-prime-zero-is-one_l=7}. 
\end{Rem}

Before we give our main criterion for checking that $\coker \phi_v$ is maximally unramified we give a basic lemma about $|\phi'(0)|_v$ and the maps $\phi_{v}^1$ and $\phi_{v,\unr}^1$. We will consider the ramification index $e_v$ of the place $v$ of $K$. Note that if $K_v=\QQ_p$ and $p\neq 2$, then the condition about the ramification index is fulfilled, i.e. we have $e_v<p-1$.

\begin{lem}\label{lem:eta_v-unr^1-is-iso}
With notation as above, the following holds.

(i) If $|\phi'(0)|_v= 1$ and $\phi_{v,\unr}^1$ is injective, then $\phi_{v}^1$ and $\phi_{v,\unr}^1$ are isomorphisms. Hence, if $\phi_{v}^1$ and $\phi_{v,\unr}^1$ are injective, then $\phi_v^1$ and $\phi_{v,\unr}^1$ are isomorphisms if and only if $|\phi'(0)|_v= 1$.

(ii) If the ramification index $e_v<p-1$, then $\phi_{v}^1$ and $\phi_{v,\unr}^1$ are injective.

(iii) If $K=\QQ$, then $\phi_p^1$ and $\phi_{p,\unr}^1$ are injective, unless $p=2$ and $2 \mid \deg \phi$.

(iv) If $K=\QQ$ and $|\phi'(0)|_v= 1$, then $\phi_{v}^1$ and $\phi_{v,\unr}^1$ are isomorphisms, unless $p=2$ and $2 \mid \deg \phi$.
\end{lem}
\proof 
Assume $|\phi'(0)|_v= 1$. Then we also have that $|\phi'(0)|_w= 1$ for all unramified finite field extensions $L_w/K_v$. Since $\phi_{v,\unr}^1$ is injective the maps $\phi_w^1:A_1(L_w)\rightarrow B_1(L_w)$ are also injective. 
By Proposition \ref{pro:coker-phi-prime-relation}, the size of the kernels and cokernels of $\phi_w^1$ agree and therefore they are isomorphisms. Hence $\phi_{v,\unr}^1$ is an isomorphism, which proves (i).

For (ii) use the isomorphism $A_1(K_v) \cong \hat A(\m_v)$. Then use \cite[IV. Theorem 6.1]{silvermanI} or the next lemma to conclude that $\phi_w^1$ is injective for any finite unramified field extension $L_w/K_v$. Hence $\phi_{v,\unr}^1$ is injective.

For (iii) apply (ii) in case $p\neq 2$. In case $2 \nmid \deg \phi$ this is due to Lemma \ref{lem:sufficient-condition-that-phi_v-1-is-an-iso}. Combing (i) and (iii) gives (iv). \EP

\begin{lem}\label{katz-lemma}
With notation as above, if the ramification index $e_v<p-1$ then the reduction-mod-$v$ map $A_0(K_v) \rightarrow \tilde \A_0(k_v)$ has torsion-free kernel, i.e. $A_1(K_v)$ is torsion-free. In particular, this gives an injection $A(K)\tor \hookrightarrow \tilde \A_0(k_v)$, thus if in addition $v$ is a place of good reduction there is an injection $A(K)\tor \hookrightarrow \tilde \A(k_v)$.
\end{lem}
\proof This is essentially the content of the Appendix of \cite{Katz}. \EP 

The next lemma and theorem are a slight generalisation of Lemmas 4.3 and 4.5 of \cite{schaefer-stoll_long-version}. Theorem \ref{thm:better-schaefer-stoll-lemma} provides our main criterion to check whether $\coker \phi_v$ is maximally unramified. To state the lemma we introduce the map 
$$\delta_v^0: \coker \phi_v^0 \rightarrow H^1(K_v,A(\overline K_v)[\phi]).$$
It is obtained by composing the natural map $\coker \phi_v^0 \rightarrow \coker \phi_v$ from Diagram (\ref{equ:tamagawa-number}) with the connecting homomorphism $\delta_v:\coker \phi_v \rightarrow H^1(K_v,A(\overline K_v)[\phi])$. Note that since $\coker \phi_v^0 \rightarrow \coker \phi_v$ need not be injective, $\delta_v^0$ may also not be injective. Similarly one defines the map 
$$\delta_{v,\unr}^0: \coker \phi_{v,\unr}^0 \rightarrow H^1(K_v^\unr, A(\overline K_v)[\phi]).$$

\begin{lem}\label{lem:better-schaefer-stoll-lemma} 
If $\phi_{v,\unr}^1$ is surjective, then the image of $\coker \phi_v^0$ under $\delta_v^0$ lies in $H^1_\unr(K_v,A(\overline K_v)[\phi])$. 
\end{lem}

\proof In the above Diagram (\ref{equ:unramified}), the first vertical map $\phi_{v,\unr}^1$ is surjective by assumption. The third vertical map $\tilde \phi^0_{\overline k_v}$ is surjective, since $\overline k_v$ is algebraically closed, therefore the middle vertical map $\phi_{v,\unr}^0$ is also surjective, i.e. $\coker \phi_{v,\unr}^0$ is trivial. The following diagram commutes. 
\begin{equation*}
  \xymatrix{
     \coker \phi_v^0 \ar[r]_/-1em/{\delta_v^0}  \ar[d] &   H^1(K_v, A(\overline K_v)[\phi])   \ar[d]_{\res_\unr}\\
      \coker \phi_{v,\unr}^0 \ar[r]_/-1em/{\delta_{v,\unr}^0}  &  H^1(K_v^\unr, A(\overline K_v)[\phi])
 }
\end{equation*}
As the lower left group is trivial, the image of the upper left group in the lower right group must be trivial, i.e. the image of $\delta_v^0$ lies in $H^1_\unr(K_v,A(\overline K_v)[\phi])$. \EP 

\begin{thm}[Main criterion for maximal unramifiedness of $\coker \phi_v$]\label{thm:better-schaefer-stoll-lemma}
If $\phi_{v,\unr}^1$ is surjective and $\phi_v^1$ and $\overline \phi_v$ are isomorphisms, then $\coker \phi_v$ is maximally unramified. 
\end{thm}
\proof As $\overline \phi_v$ is an isomorphism, $\coker \phi_v^0 \rightarrow \coker \phi_v$ is also an isomorphism, and so by the above lemma we have that $\coker \phi_v$ maps injectively onto a subgroup of $H^1_\unr(K_v,A(\overline K_v)[\phi])$. But these two groups have same cardinality, since $\# H^1_\unr(K_v,A(\overline K_v)[\phi]) = \# \ker \phi_v$ by Lemma \ref{lem:H1-unr}, and $\# \ker \phi_v= \# \coker \phi_v$ by Lemma \ref{lem:isogenies-on-k_v-points}, as $\overline \phi_v$ and $\phi_v^1$ are isomorphisms. \EP

Our assumptions on $\phi_v^1$ and $\phi_{v,\unr}^1$ in Lemma \ref{lem:better-schaefer-stoll-lemma} and Theorem \ref{thm:better-schaefer-stoll-lemma} replaced the assumption $v \nmid \deg \phi$ in Lemmas 4.3 and 4.5 of \cite{schaefer-stoll_long-version}. We have seen in Lemma \ref{lem:sufficient-condition-that-phi_v-1-is-an-iso} that $v \nmid \deg \phi$ is a stronger assumption, hence we can easily deduce the original result of Schaefer and Stoll.

\begin{cor}[Criterion for maximal unramifiedness of $\coker \phi_v$ in case $v \nmid \deg \phi$]\label{cor:better-schaefer-stoll-lemma}
If $v \nmid \deg \phi$ and $\gcd(\deg \phi, c_{A,v} \cdot c_{B,v})=1$ then $\coker \phi_v$ is maximally unramified. 
\end{cor}
\proof This is Lemma 4.5 of \cite{schaefer-stoll_long-version}. The corollary follows directly from Theorem \ref{thm:better-schaefer-stoll-lemma} together with Lemma \ref{lem:sufficient-condition-that-phi_v-1-is-an-iso} and  Corollary \ref{cor:kernels-and-cokernels-and-Tamagawa-quotient-are-powers-of-l}. \EP

We also want to apply Theorem \ref{thm:better-schaefer-stoll-lemma} in case $v\mid \deg \phi$. As already seen in Lemma \ref{lem:eta_v-unr^1-is-iso} we can replace $v\nmid \deg \phi$ with the condition that the ramification index $e_v<p-1$ and that $|\phi'(0)|_v= 1$. 

\begin{cor}[Criteria for maximal unramifiedness of $\coker \phi_v$ in case $v \mid \deg \phi$]\label{cor:condiction-for-unramified-cokernel}

Assume that the ramification index $e_v<p-1$.

(i) If $|\phi'(0)|_v= 1$ and $\gcd(\deg \phi, c_{A,v} \cdot c_{B,v})=1$, then $\coker \phi_v$ is maximally unramified.

(ii) If $v$ is a place of good reduction, then $\coker \phi_v$ is maximally unramified if and only if $|\phi '(0)|_v= 1$.

(iii) If $v$ is a place of good reduction and $\phi$ has a $K_v$-kernel, then $|\phi '(0)|_v= 1$ and $\coker \phi_v$ is maximally unramified.
\end{cor}
\proof 
For (i) combine Lemma \ref{lem:eta_v-unr^1-is-iso} with Theorem \ref{thm:better-schaefer-stoll-lemma} and Corollary \ref{cor:kernels-and-cokernels-and-Tamagawa-quotient-are-powers-of-l}. 

For (ii) note that $v$ being a place of good reduction implies that $c_{A,v} = c_{B,v}=1$. If $|\phi'(0)|_v= 1$, then by (i) we get that $\coker \phi_v$ is maximally unramified. Now assume that $\coker \phi_v$ is maximally unramified, hence its cardinality equals the cardinality of $\ker \phi_v$. By Corollary \ref{cor:isogenies-on-k_v-points} we get that $|\phi'(0)|_v= c_{B,v}/c_{A,v}=1$, which completes (ii). For (iii), combine (ii) with Lemmas \ref{lem:good-reduction-is-unramified} and \ref{lem:eta_v-unr^1-is-iso}. \EP

We summarise all the criteria for maximal unramifiedness for the case that $K=\QQ$. The first one is especially useful when $A$ and $B$ are elliptic curves as this makes the values $|\phi'(0)|_p$, $c_{A,p}$, and $c_{B,p}$ easily compultable.

\begin{cor}[Criteria for maximal unramifiedness of $\coker \phi_p$ in case $K=\QQ$]\label{cor:criteria-for-unramified-cokernel-if-K=QQ} Let $\phi:A \rightarrow B$ be an isogeny between two abelian varieties $A$ and $B$ over $\QQ$ and let $p$ be a prime such that $p\neq 2$ if $2 \mid \deg \phi$.

(i) If $|\phi'(0)|_p=1$ and $\gcd(\deg \phi, c_{A,p} \cdot c_{B,p})=1$, then $\coker \phi_p$ is maximally unramified.

(ii) If $p$ is a place of good reduction and $\phi$ has a $\QQ_p$-kernel, then $|\phi '(0)|_p= 1$ and $\coker \phi_p$ is maximally unramified.
\end{cor} 
\proof Follows directly from Lemma \ref{lem:eta_v-unr^1-is-iso}, Theorem \ref{thm:better-schaefer-stoll-lemma}, and Corollary \ref{cor:condiction-for-unramified-cokernel}. \EP

We end this subsection with a basic lemma about the infinite places.

\begin{lem}\label{lem:size-of-cokernel-at-infinity}
Let $L$ be either $\RR$ or $\CC$. Let $A$ and $B$ be abelian varieties over $L$ and let $\phi:A \rightarrow B$ be an isogeny and denote by $\phi_\infty :A(L) \rightarrow B(L)$ the induced group homomorphism on $L$-rational points.

(i) If $L=\CC$, then $\# \coker \phi_\infty / \# \ker \phi_\infty = 1 / \deg \phi$.

(ii) If $L=\RR$, then $\# \coker \phi_\infty=1$, if $2 \nmid \deg \phi$.
\end{lem}
\proof
The first part is obvious, as $\CC$ is algebraically closed and of characteristic $0$, hence $\phi_\infty$ is surjective and the size of the kernel equals the degree. For (ii) note that $\coker \phi_\infty$ embeds into $H^1(\RR, A[\phi])$, which is trivial if the order of $\Gal_\RR$ is coprime to $A[\phi]$. \EP

In the next subsection we will consider the special case of $A$ and $B$ being elliptic curves $E$ and $E'$ and the isogeny being of prime degree $\l$ and having a $K_v$-kernel. Isogenies between elliptic curves will usually be denoted by $\eta$ and as before we are interested in whether $\coker \eta_v$ is maximal, maximally unramified, or trivial. The goal is to classify $\coker \eta_v$ with respect to the reduction type at $v$. To show maximal unramifiedness we will use the criteria established above, and we will use the equation 
$$\# \coker \eta_v = \l \cdot |\eta'(0)|_v^{-1} \cdot \frac{c_{E',v}}{c_{E,v}}.$$
from Lemma \ref{cor:isogenies-on-k_v-points} to compute the order of $\coker \eta_v$ directly to see whether it is trivial or maximal. In any case we want to have a way to compute $|\eta'(0)|_v$ and the Tamagawa numbers $c_{E',v}$ and $c_{E,v}$ with respect to the reduction type at $v$.

\subsection{Isogenies of prime degree between elliptic curves over local fields}

The notation for this subsection is the following. Fix two prime numbers $p$ and $\l$. Note that $p=\l$ is allowed. Let $E$ and $E'$ be elliptic curves over a $p$-adic field $K_v$ and let $\eta:E\rightarrow E'$ be an isogeny of prime degree $\l$. We will use the notations from Diagrams (\ref{equ:reduction_maps}), (\ref{equ:tamagawa-number}), and (\ref{equ:unramified}) with $A=E$ and $B=E'$. Assuming that $\eta$ has a $K_v$-kernel, the goal of this subsection is to determine under which further assumptions the reduction type of $E$ at $v$ determines whether $\coker \eta_v$ is maximal, maximally unramified, or trivial. In the case when $K_v=\QQ_p$ and $\l \geq 5$ we can give a complete classification, which will be stated in our main theorem \ref{thm:iota-p-1-of-E} at the end of this subsection. We will see that if the reduction type at $p$ is not additive, then we can weaken the condition on $\l$.

We start with computing the quotient $c_{E'}/c_E$ of Tamagawa numbers with respect to the reduction type at $v$. In most cases the Tamagawa quotient of isogenous elliptic curves can easily be computed with Tate's algorithm and the theory of Tate curves. See for example the Appendix of \cite{DokchitserFisher} by Tom Fisher or \cite[\textsection 6 and \textsection 9]{DokchitserLocal} by Tim and Vladimir Dokchitser. 

\begin{lem}\label{lem:Tamagawa-quotient-is-1}
Suppose that $E/K_v$ has
\begin{enumerate}
 \item good reduction, or
\item non-split multiplicative reduction and $\l \neq 2$, or
\item additive reduction and $\l \geq 5$.
\end{enumerate}
Then the group homomorphism $\overline \eta_v$ is an isomorphism, and hence $c_{E'}/c_{E}=1$. 
\end{lem}
\proof 
In case of good reduction this is clear, since $c_{E,v}=c_{E',v}=1$. From Tate's algorithm \cite{Tate_algorithm} it follows that $c_{E,v}$ and $c_{E',v}$ are at most $4$ in the additive case, and at most $2$ in the non-split multiplicative case. Hence $\gcd(\deg \eta,c_{E,v} \cdot c_{E',v})=1$, and the result follows directly from Corollary \ref{cor:kernels-and-cokernels-and-Tamagawa-quotient-are-powers-of-l}. \EP

To calculate the Tamagawa quotient in the case of split multiplicative reduction we use the theory of Tate curves.
 
 \begin{thm}\label{thm:tate-curve} \textnormal{(Tate)}
 Assume that $E/K_v$ has split multiplicative reduction. Then there is a unique $q \in K_v^*$, s.t. $v(q)>0$, and we have the following Galois equivariant $p$-adic analytic isomorphism
 $$E(L) \cong L^*/q^\ZZ,$$
 for all algebraic field extension $L/K_v$. Moreover, $c_{E,v} = v(q)$.
 \end{thm}
 
 \proof See \cite[V Theorem 5.3]{silvermanII}. The last statement follows from the proof of the surjectivity of the Tate map \cite[V.4]{silvermanII}.\EP

\begin{Rem}
If $E/K_v$ is an elliptic curve having split multiplicative reduction we have $E(\overline K_v) \cong \overline K_v^*/q^\ZZ$, for $q \in K_v^*$ and $v(q)>0$. We want to classify which Galois invariant subgroups of prime order $\l$ exist and whether they are contained in the connected component of the identity $E_0(\overline K_v)\cong \overline K_v^*$. Since they are all subgroups of $E(\overline K_v)[\l]\cong \ZZ/\l\ZZ \times \ZZ/\l\ZZ$, there are at most $\l+1$ of such groups. The $\l$-th roots of unity $\{ \xi_\l, \xi_\l^2, \ldots, \xi_\l^{\l-1},1 \}$ in $\overline K_v^*$ form a Galois invariant subgroup of $\overline K_v^*/q^\ZZ$, which is contained in the connected component of the identity, and a generator is defined over $K_v$ if and only if $\mub_\l \subseteq K_v$. Hence this subgroup of $E(\overline K_v)[\l]$ is isomorphic to $\mub_\l$ as a Galois module. The remaining $\l$ subgroups are defined over $K_v(\sqrt[\l]{q}, \mub_\l)$. 
None of these $\l$ subgroups are contained in the connected component of the identity. They are generated by $\xi_\l^i \sqrt[\l]{q}$, for $i=0, \ldots,\l-1$.
The  elements of such a subgroup have different minimal polynomials, hence if the subgroup is Galois invariant then it is ismorphic to $\ZZ/\l\ZZ$ as a Galois module.
\end{Rem}

\begin{lem}\label{lem:split-reduction}
With notation as above, if $E/K_v$ has split multiplicative reduction, then
$$\frac{c_{E'}}{c_{E}} = \begin{cases}
		    1/\l, & \ker \eta_v \nsubseteq E_0(\overline K_v)\\
                    \l, & \ker \eta_v \subseteq E_0(\overline K_v).\\
                   \end{cases}$$
\end{lem}
\proof This is Lemma A.2 of the appendix of \cite{DokchitserFisher} by Tom Fisher. 
To identify our lemma with Lemma A.2, use the last remark to see that $\ker \eta_v \nsubseteq E_0(\overline K_v)$ implies that $\ker \eta_v \cong \ZZ/\l\ZZ$, and $\ker \eta_v \subseteq E_0(\overline K_v)$ implies that $\ker \eta_v \cong \mub_\l$ as Galois modules. 

We give now a slitely longer version of Tom Fisher's proof. By theorem \ref{thm:tate-curve} we have $E(\overline K_v) \cong \overline K_v^*/q_1^\ZZ$ and $E'(\overline K_v) \cong \overline K_v^*/q_2^\ZZ$. If $\ker \eta_v \nsubseteq E_0(\overline K_v)$ then $\ker \eta_v = \langle [\xi_\l^i \sqrt[\l]{q_1}] \rangle$, for an $i \in \{0, \ldots, \l-1\}$, and $\eta_v:\overline K_v^*/q_1^\ZZ\rightarrow \overline K_v^*/q_2^\ZZ$ is given by $[x]\mapsto [x]$ and $q_2=\xi_\l^i \sqrt[\l]{q_1}$. Therefore 
 $$\frac{c_{E'}}{c_{E}} = \frac{v(q_2)}{v(q_1)}= \frac{v(\xi_\l^i \sqrt[\l]{q_1})}{v(q_1)}=1/\l.$$
 If $\ker \eta_v \subseteq E_0(\overline K_v)$ then $\ker \eta_v = \langle [\xi_\l] \rangle$ and $\eta_v:\overline K_v^*/q_1^\ZZ\rightarrow \overline K_v^*/q_2^\ZZ$ is given by $[x]\mapsto [x^\l]$ and $q_2=q_1^\l$. Therefore 
 $$\frac{c_{E'}}{c_{E}} = \frac{v(q_2)}{v(q_1)}= \frac{v(q_1^\l)}{v(q_1)}=\l,$$
which completes the proof. \EP

Now we study the implications of $\ker \eta_v$ being or not being part of the connected component of the identity $E_0(\overline K_v)$. The result is essentially a corollary of Tate's algorithm \cite{Tate_algorithm} and explores the fact that $\eta$ is of prime degree $\l$.

\begin{lem} \label{lem:reduction-type-vs-roots-of-unity}
With notation as above, we have:

(i) If $\ker \eta_v \nsubseteq E_0(\overline K_v)$, then $\eta$ has a $K_v$-kernel, $\eta_v^1$ is an isomorphism, $|\eta'(0)|_v=1$, $\l \mid c_E$, and exactly one of the following three cases holds. 

\begin{itemize}
\item $E$ has split multiplicative reduction,
\item $E$ has non-split multiplicative reduction and $\l=2$,
\item $E$ has additive reduction and either $\l=2$ or $\l=3$.
\end{itemize}

(ii) If $\ker \eta_v \subseteq E_0(\overline K_v)$, assume additionally that $\eta$ has a $K_v$-kernel and that $\eta_v^1$ is injective. We have the following implications.
  \begin{itemize}
  \item $E$ has multiplicative reduction $\Rightarrow \ v \nmid \l$ and $|\eta'(0)|_v=1$,
\item $E$ has split multiplicative reduction $\Rightarrow \ \mub_\l \subseteq K_v$,
\item $E$ has non-split multiplicative reduction and $\l\neq 2$ $\Rightarrow \ \mub_\l \nsubseteq K_v$,
\item $E$ has additive reduction $\Rightarrow \ v|\l$.
  \end{itemize}
\end{lem}
\proof 
If $\ker \eta_v$ is trivial, then it is clearly contained in $E_0(\overline K_v)$. Hence $\ker \eta_v \nsubseteq E_0(\overline K_v)$ implies that $\ker \eta_v$ is non-trivial, and therefore $\eta$ has a $K_v$-kernel as its degree is prime. It also implies that $\eta_{\overline K_v}^0$, $\eta_v^0$, and thus $\eta_v^1$ are injective. From the triviality of $\eta_{\overline K_v}^0$ is follows that $H^1(K_v, E_0(\overline K_v)[\eta])$ is trivial and hence $\coker \eta_v^0$ is also trivial.
We deduce that $\eta_v^0$ is an isomorphism, and therefore $\tilde \eta_v^0$ is surjective. By Lemma \ref{lem:number-of-reduced-points-is-the-same}, $\tilde \eta_v^0$ is an isomorphism, as its kernel and cokernel have equal cardinalities. This immediately implies that $\eta_v^1$ is an isomorphism, which gives $|\eta'(0)|_v=1$ by Proposition \ref{pro:coker-phi-prime-relation}. Again by the fact that $\eta_v^0$ is an isomorphism, it follows that $\#\ker \overline \eta_v=\l$, which gives that $\l$ divides the Tamagawa number $c_E$. In particular, the reduction type is bad. By \cite{Tate_algorithm}, $c_E$ is $\leq 2$ in the non-split multiplicative case and $\leq 4$ in the additive case, giving (i). 

For (ii) let $P \in E(K_v)$ be a generator of $\ker \eta_v$. If $\ker \eta_v \subseteq E_0(\overline K_v)$ and $\eta$ has a $K_v$-kernel, then $P$ generates $\ker \eta_v^0$. Since we assumed $\eta_v^1$ to be injective, the order of $\overline P$ is $\l$. Set $|k_v|=:p^f$. The order of $\overline P$ divides the cardinality of $\tilde \E_0(k_v)$, which is either $p^f-1$, $p^f+1$, or $p^f$, depending on whether the reduction type is split multiplicative, non-split multiplicative, or additive, respectively \cite[\textsection 7]{Tate_algorithm}. Therefore we get the following implications
\begin{itemize}
\item multiplicative $\Rightarrow \ p^f\not \equiv 0 \mod \l \ \ \Rightarrow \ p \neq \l \ \ \Rightarrow \ v \nmid \l,$
\item split $\Rightarrow \ p^f\equiv 1 \mod \l \ \ \Rightarrow \ \mub_\l \subseteq k_v \ \ \Rightarrow \ \mub_\l \subseteq K_v,$
\item non-split and $\l \neq 2$ $\Rightarrow \ p^f \not \equiv 0,1 \mod \l  \Rightarrow \ \mub_\l \nsubseteq k_v  \Rightarrow \ \mub_\l \nsubseteq K_v,$
\item additive $\Rightarrow \ p^f\equiv 0 \mod \l \ \ \Rightarrow \ p = \l,  \ \ \Rightarrow \ v \mid \l.$
\end{itemize}
By Proposition \ref{pro:coker-phi-prime-relation}, $v \nmid \l$ implies $|\eta'(0)|_v=1$, which completes (ii). \EP

We summarise the results and state under which further assumptions $\coker \eta_v$ is trivial, maximally unramified, or maximal in the case of multiplicative reduction.

\begin{cor}[Criteria to classify $\coker \eta_v$ is case of multiplicative reduction]\label{cor:multipl-red-gives-trivial-or-maximal-coker}
(i) If the reduction type of $E/K_v$ is split multiplicative and $\ker \eta_v \nsubseteq E_0(\overline K_v)$, then $|\eta'(0)|_v=1$ and $\coker \eta_v$ is trivial.

(ii) If the reduction type of $E/K_v$ is split multiplicative, $\ker \eta_v \subseteq E_0(\overline K_v)$, $\eta$ has a $K_v$-kernel, and $\eta_v^1$ is injective, then $v \nmid \l, \mub_\l \subseteq K_v$,  $|\eta'(0)|_v=1$, and $\coker \eta_v$ is maximal.

(iii) If the reduction type of $E/K_v$ is non-split multiplicative, $\l \neq 2$, $\eta$ has a $K_v$-kernel, and $\eta_v^1$ is injective, then $v \nmid \l, \mub_\l \nsubseteq K_v$, $|\eta'(0)|_v=1$ and $\coker \eta_v$ is maximally unramified.

(iv) If the reduction type of $E/K_v$ is non-split multiplicative, $\l = 2$, $v \nmid \l$, and $\eta$ has a $K_v$-kernel, then $\mub_\l \subseteq K_v$ and $|\eta'(0)|_v=1$. Further $\coker \eta_v$ is trivial if $c_{E'} / c_{E}=1/2$, $\coker \eta_v$ is maximal if $c_{E'} / c_{E}=2$, and $\coker \eta_v$ is maximally unramified if $c_{E} = c_{E'}=1$.
\end{cor}
\proof Lemma \ref{lem:reduction-type-vs-roots-of-unity} already contains everything of (i)-(iii) but the statement whether $\coker \eta_v$ is trivial, maximally unramified, or maximal. In (iv) we get $|\eta'(0)|_v=1$ and $\mub_\l \subseteq K_v$, as $v \nmid \l$ and $\l=2$. It remains to classify $\coker \eta_v$. 

By Corollary \ref{cor:isogenies-on-k_v-points} we obtain the equation $\# \coker \eta_v = \l \cdot  c_{E'} / c_{E}$, and the size of $H^1(K_v,E(\overline K_v)[\eta])$ is given by Corollary \ref{cor:H1}. The Tamagawa quotient in (i)-(iii) can be computed with Lemmas \ref{lem:Tamagawa-quotient-is-1} and \ref{lem:split-reduction} This shows triviality of $\coker \eta_v$ in (i) and the first case in (iv), and maximality in all other cases but the third case of (iv). Note that in (iii), $H^1(K_v,E(\overline K_v)[\eta])$ equals the unramified subgroup, as its cardinality is $\l$. To get the maximal unramifiedness in the third part of (iv) use Corollary \ref{cor:better-schaefer-stoll-lemma}. \EP

We finish with the main theorem of this subsection. Recall that we call $\coker \eta_p$ {\em maximal} if it equals the full $H^1(\QQ_p,E(\overline \QQ_p)[\eta])$, and {\em maximally unramified} if it equals the unramified subgroup $H^1_\unr(\QQ_p,E(\overline \QQ_p)[\eta])$; see the discussion before Remark \ref{rem:unramified-is-good-as-local-quotient-equals-1}. The definition of $|\eta'(0)|_p$ is given before Proposition \ref{pro:coker-phi-prime-relation} and having a $\QQ_p$-kernel means that $E(\overline \QQ_p)[\eta]=E(\QQ_p)[\eta]$.

\begin{thm}[Criteria to classify $\coker \eta_p$ is case $\eta$ has a $\QQ_p$-kernel and is of prime degree] \label{thm:iota-p-1-of-E}
Let $E$ and $E'$ be elliptic curves over $\QQ_p$ and let $\eta:E\rightarrow E'$ be an isogeny of prime degree $\l$, and assume that $\eta$ has a $\QQ_p$-kernel. Then the left column of the table below implies the two columns to the right and in all but the last row we also get that $|\eta'(0)|_p=1$.
\begin{center}
    \begin{tabular}{|c||c|c|}
    \hline
    reduction type of $E/\QQ_p$,& $p=$ or $\neq \l$, & $\coker \eta_p$\\
    plus further assumptions & $\mub_\l \subseteq$ or $\nsubseteq \QQ_p$ & is\\
    \hline
\text{split multiplicative}, $\ker \eta_p \nsubseteq E_0(\overline \QQ_p)$ & no implications   &  \text{trivial}\\
    \text{split multipl.}, $\ker \eta_p\subseteq E_0(\overline \QQ_p)$, $p\neq 2$ \text{ or } $\l \neq 2$ & $p\neq \l, \mub_\l \subseteq \QQ_p$ &  \text{maximal}\\ 
       \hline
\text{non-split multiplicative}, $\l \neq 2$ &  $p\neq \l, \mub_\l \nsubseteq \QQ_p$   & \text{max. unramified}\\ 
\text{non-split multiplicative}, $\l = 2\neq p$, $c_{E'} / c_{E}=1/2$&  $p\neq \l, \mub_\l \subseteq \QQ_p$   & \text{trivial}\\ 
\text{non-split multiplicative}, $\l = 2 \neq p$, $c_{E'} / c_{E}=2$&  $p\neq \l, \mub_\l \subseteq \QQ_p$   & \text{maximal}\\ 
\text{non-split multiplicative}, $\l = 2 \neq p$, $c_{E}= c_{E'}=1$&  $p\neq \l, \mub_\l \subseteq \QQ_p$   & \text{max. unramified}\\ 
   \hline
\text{good}, $p\neq 2$ \text{ or } $\l \neq 2$ &  no implications  &   \text{max. unramified}\\
    \hline
\text{additive}, $\l \geq 5$, $|\eta'(0)|_p=1$& $p = \l, \mub_\l \nsubseteq \QQ_p$    &  \text{max. unramified}\\ 
\text{additive}, $\l \geq 5$, $|\eta'(0)|_p\neq 1$& $p = \l, \mub_\l \nsubseteq \QQ_p$    & \text{maximal}\\
    \hline
    \end{tabular}
\end{center}
\end{thm}
\proof For all but the first row of the table we use Lemma \ref{lem:eta_v-unr^1-is-iso} to deduce that $\eta_p^1$ is injective. Then the six cases of multiplicative reduction are contained in the last corollary and the case of good reduction is covered by Lemma \ref{lem:good-reduction-is-unramified}. 

In the additive case, due to $\l \geq 5$, we get that $p=\l$ by Lemma \ref{lem:reduction-type-vs-roots-of-unity} and hence $\mub_\l \nsubseteq \QQ_p$. This implies that $\# H^1(\QQ_p, E(\overline \QQ_p)[\eta])=\l^2$ by Corollary \ref{cor:H1}. Further $\overline \eta_p$ is an isomorphism by Lemma \ref{lem:Tamagawa-quotient-is-1} and thus by Corollary \ref{cor:isogenies-on-k_v-points} we have $\# \coker=\l \cdot |\eta'(0)|_p^{-1}$. We know that $|\eta'(0)|_p^{-1} \geq 1$ as $\eta_p^1$ is injective. Hence, there are two possibilities. Firstly, $\# \coker \eta_p=\l$, which is equivalent to $|\eta'(0)|_
p=1$, and secondly, $\# \coker \eta_p=\l^2$, which is 
equivalent to $|\eta'(0)|_p\neq 1$, and which implies that $\coker \eta_p$ is maximal. It remains to show that $\coker \eta_p$ is maximally unramified in case the reduction type is additive and $|\eta'(0)|_p=1$. At this point we apply Theorem \ref{thm:better-schaefer-stoll-lemma}. All conditions are fulfilled due to Lemma \ref{lem:eta_v-unr^1-is-iso}. This finishes the proof. \EP 

In the next subsection we present the general concept of non-simple abelian varieties and of isogenies with diagonal kernel.

\subsection{Non-simple abelian varieties and isogenies with diagonal kernel}

In this subsection $K$ will always denote a field of characteristic $0$. An abelian variety $B/K$ is called \emph{non-simple} if it is isogenous to a product of two abelian varieties $A_1/K$ and $A_2/K$:
$$\phi: A_1 \times A_2 \rightarrow B.$$
Recall, that if we do not specify the field of definition of an isogeny $\phi$ between two abelian varieties which are defined over a field $K$, then want $\phi$ to be defined over $K$, too.
Let $A_1,A_2$ and $B$ be abelian varieties over a field $K$ and let $\phi:A_1 \times A_2 \rightarrow B$ be an isogeny. We say that $\phi$ has \emph{diagonal kernel}, or simply say that $\phi$ is {\em diagonal}, if there is a finite group scheme $G$ over $K$ contained in both $A_i$, together with fixed embeddings $\iota_i: G \hookrightarrow A_i$, such that the kernel of $\phi$ is the natural embedding of $G$ into the product $A_1 \times A_2$ via $\iota_1 \times \iota_2$. We denote the image of $G$ in $A_i$ by $G_i:=\iota_i(G)$. Clearly, $G$, $G_1$, and $G_2$ are pairwise isomorphic as finite group schemes and the $\overline K$-rational points of $G_1$ and $G_2$ form isomorphic Galois modules; hence there is a Galois equivariant isomorphism $\alpha: G_1 \rightarrow G_2$ such that $\iota_2=\alpha \circ \iota_1$ and that $\ker \phi$ equals the graph of $\alpha$. 
Further, both $A_i$ possess an isogeny $\eta_i:A_i \rightarrow A_i'$ which is defined through its kernel by setting $\ker \eta_i := G_i$ and $A'_i:=A_i/G_i$. Clearly $\ker \eta_1 \cong \ker \eta_2 \cong \ker \phi$.

We start with presenting our Key Lemma to controll the local quotient for isogenies with diagonal kernel. First we state a basic lemma about Galois cohomology.

\begin{lem}\label{lem:alpha-maps-unramified-to-unramified}
Let $K$ be a field and let $G_1$ and $G_2$ be two finite $K$-Galois modules. Assume $\alpha:G_1 \rightarrow G_2$ is a Galois equivariant homomorphism. Then the map 
$$\alpha^*:H^1(K, G_1) \rightarrow H^1(K,G_2), \ \ [\xi] \mapsto [\alpha \circ \xi],$$
is a well-defined group homomorphism. 
If in addition $\alpha$ is an isomorphism, then $\alpha^*$ is an isomorphism, too. Further, the isomorphism $\alpha^*$ respects the Inflation-Restriction sequence for normal subgroups of $\Gal_K$, i.e. for any Galois extension $L/K$, $\alpha^*$ induces an isomorphism $H^1(\Gal(L/K),G_1^{\Gal_L}) \rightarrow H^1(\Gal(L/K),G_2^{\Gal_L})$ and an isomorphism $H^1(L,G_1) \rightarrow H^1(L,G_2)$ and these isomorphisms commute with the Inflation-Restriction sequence. 

In particular, if $K=K_v$ is a local field, for every Galois equivariant isomorphism $\alpha:G_1 \rightarrow G_2$, the isomorphism $\alpha^*:H^1(K_v, G_1) \rightarrow H^1(K_v,G_2)$ induces an isomorphism between the unramified subgroups $H^1_\unr(K_v,G_1)$ and $H^1_\unr(K_v,G_2)$.
\end{lem}
\proof Follows directly from the functoriality of Galois cohomology. \EP

\begin{lem}[Key Lemma to compute the local quotient for isogenies with diagonal kernel]\label{lem:hitting-both-H-1-general-version}
Let $A_1$ and $A_2$ be two abelian varieties over a number field $K$ and let $\phi:A_1 \times A_2 \rightarrow B$ be an isogeny with diagonal kernel. Denote by $\eta_i:A_i \rightarrow A'_i$ the isogenies for which there is a Galois equivariant isomorphism $\alpha:\ker \eta_1 \rightarrow \ker \eta_2$ whose graph equals $\ker \phi$. Let $v \in M_K^0$ be a finite place of $K$. Then we have the following.

(i) $\coker \phi_v$ is maximal if and only if $\coker \eta_{1,v}$ and $\coker \eta_{2,v}$ are both maximal.

(ii) $\coker \phi_v$ is trivial if either $\coker \eta_{1,v}$ or $\coker \eta_{2,v}$ is trivial.

(iii) $\coker \phi_v$ is maximally unramified if either $\coker \eta_{1,v}$ or $\coker \eta_{2,v}$ is maximally unramified and the other one is maximally unramified or maximal.
\end{lem}
\proof Define the two Galois equivariant isomorphisms 
$$\gamma_1 := (id, \alpha) : \ker \eta_1 \rightarrow \ker \phi \ \mbox{  and  } \ \gamma_2 :=(\alpha^{-1},id): \ker \eta_2 \rightarrow \ker \phi.$$
By the above lemma we get two group isomorphisms 
$\gamma_i^* : H^1(K_v, A_i[\eta_i]) \rightarrow H^1(K_v,(A_1 \times A_2)[\phi]).$
Thus, for any $[\xi] \in H^1(K_v,(A_1 \times A_2)[\phi])$ there is a unique $[\xi_1] \in H^1(K_v, A_1[\eta_1])$ and a unique $[\xi_2] \in H^1(K_v, A_2[\eta_2])$ such that $\gamma_1^*([\xi_1])=\gamma_2^*([\xi_2])=[\xi]$. It follows that $\xi(\sigma) = (\xi_1(\sigma),\alpha(\xi_1(\sigma)))$ and $\xi(\sigma) = (\alpha^{-1}(\xi_2(\sigma)),\xi_2(\sigma))$, and hence
$\xi(\sigma) = (\xi_1(\sigma),\xi_2(\sigma)), \mbox{ for all } \sigma \in \Gal_{K_v}.$
Thus
\begin{equation}\label{equ:key-lemma}
[\xi] \in \coker \phi_v \ \Leftrightarrow \ [\xi_1] \in \coker \eta_{1,v} \mbox{ and } [\xi_2] \in \coker \eta_{2,v},
\end{equation}
since both assertions are equivalent to
the existence of $P_1 \in A_1(\overline K_v)$ and  $P_2 \in A_2(\overline K_v)$, such that for all $\sigma \in \Gal_{K_v}$ we have $\xi_1(\sigma)=P_1^\sigma-P_1$ and $\xi_2(\sigma)=P_2^\sigma-P_2$.
For (ii), recall that $[\xi]$ is the trivial class if and only if $[\xi_1]$ and $[\xi_2]$ are both the trivial class. For (iii) use the above lemma again to get that $[\xi] \in H^1_\unr(K_v,(A_1 \times A_2)[\phi])$ if and only if $[\xi_1] \in H^1_\unr(K_v,A_1[\eta_1])$ and $[\xi_2]\in H^1_\unr(K_v,A_2[\eta_2])$. Now everything follows directly from (\ref{equ:key-lemma}). \EP

\begin{Rem}
The Key Lemma shows that if one knows whether $\coker \eta_{1,v}$ and $\coker \eta_{2,v}$ are maximal, maximally unramified, or trivial, then one knows whether $\coker \phi_v$ is maximal, maximally unramified, or trivial. This is enough to determine the quotient $\# \coker \phi_v / \# \ker \phi_v$, as seen in Section \ref{sec:isogenies-over-local-fields}. For all examples of cyclic isogenies $\phi:E_1/\QQ \times E_2/\QQ \rightarrow B/\QQ$ we are going to consider we will compute $\coker \phi_p$ by first computing $\coker \eta_{i,p}$ and then applying the Key Lemma. This is the reason why Theorem \ref{thm:iota-p-1-of-E} is so important, as it classifies the $\coker \eta_{i,p}$. 
\end{Rem}

For fixed abelian varieties $A_1/K$ and $A_2/K$ and fixed isomorphic finite subgroup schemes $G_1/K \subset A_1$ and $G_2/K \subset A_2$, we can define an isogeny $\phi:A_1 \times A_2 \rightarrow B$ with diagonal kernel by setting the kernel of $\phi$ to be equal to the graph of $\alpha$. Note that $\phi$ and $B$ depend on the choice of $\alpha$, which we may denote by $\phi_\alpha$ and $B_\alpha$ to emphasise it. We will now show that the order of $\sha(B_\alpha/K)$ is independent of $\alpha$ if $\phi$ is a cyclic isogeny. Recall, that an isogeny is called {\em cyclic}, if the $\overline K$-rational points of its kernel form a cyclic group. 

\begin{pro}\label{pro:order-of-sha-independent-of-alpha}
Let $A_1$ and $A_2$ be two abelian varieties over a number field $K$, such that there are isomorphic finite cyclic $K$-subgroup schemes $G_1 \subseteq A_1$ and $G_2 \subseteq A_2$. Choose a Galois equivariant isomorphism $\alpha:G_1 \rightarrow G_2$ and let $\phi_\alpha:A_1 \times A_2 \rightarrow B_\alpha$ be the cyclic isogeny with diagonal kernel such that $\ker \phi_\alpha$ equals the graph of $\alpha$. Then $\# \sha(B_\alpha/K)$ is independent of the choice of $\alpha$.
\end{pro}
\proof What we will do is to take the Cassels-Tate equation 
$$ \frac{\# \sha(A_1 \times A_2/K)}{\# \sha(B_\alpha/K)} =\frac{\# \ker \phi_{\alpha,K}}{\# \coker \phi_{\alpha,K}}\frac{\# \coker \phi^\vee_{\alpha,K}}{\# \ker \phi^\vee_{\alpha,K}} \prod_{v \in M_K} \frac{\# \coker \phi_{\alpha,v} }{ \# \ker \phi_{\alpha,v}}$$
and then we show that the cardinality of all occurring kernels and cokernels on the right hand side are independent of $\alpha$. The set of $\overline K$-rational points of the kernels of the isogenies $\phi_\alpha: A_1 \times A_2 \rightarrow B_\alpha$ and $\phi^\vee_\alpha: B^\vee_\alpha \rightarrow A_1^\vee \times A_2^\vee$ depend on $\alpha$. But the isomorphism class of $\ker \phi_\alpha$ and of $\ker \phi^\vee_\alpha$ as a Galois module is fixed, hence it is clear that the size of all occurring kernels in the Cassels-Tate equation are unaffected by $\alpha$. It remains to consider the cokernels.

Fix two Galois equivariant isomorphisms $\alpha:G_1 \rightarrow G_2$ and $\alpha':G_1 \rightarrow G_2$. Then there is a Galois equivariant automorphism $\beta_2$ of $G_2$, such that $\alpha' = \beta_2 \circ \alpha$. The Galois equivariant automorphism $\gamma_2:=id \times \beta_2$ of $G_1 \times G_2$ induces a Galois equivariant isomorphism between $\ker \phi_\alpha$ and $\ker \phi_{\alpha'}$. As $G_2$ is cyclic, the automorphism $\beta_2$ is multipication by some factor, hence there is an endomorphism $B_2$ of $A_2$ such that the restriction of $B_2$ to $G_2$ equals $\beta_2$. The following diagram has exact rows and commutes, with the vertical maps $id \times B_2$ and $[id \times B_2]$ being isogenies.

$$  \xymatrix{
      0 \ar[r] & \ker \phi_\alpha \ar[r]  \ar[d]_{\gamma_2} &   A_1 \times A_2 \ar[r]^{\phi_\alpha}  \ar[d]_{id \times B_2} & B_\alpha \ar[r] \ar[d]_{[id \times B_2]} & 0\\
      0 \ar[r] & \ker \phi_{\alpha'} \ar[r]  & A_1 \times A_2 \ar[r]^{\phi_{\alpha'}}  & B_{\alpha'} \ar[r]  & 0\\
 }$$

Applying Galois cohomology yields the following commutative diagram with exact rows, where $L$ is either the number field $K$ or one of its completions $K_v$.

$$  \xymatrix{
      0 \ar[r] & \coker \phi_{\alpha,L} \ar[r]  \ar[d] &   H^1(L, \ker \phi_\alpha) \ar[r]^{\iota^1_\alpha}  \ar[d]_{\gamma_2^*} & H^1(L,A_1 \times A_2) \ar[d]_{(id \times B_2)^*} \ar[r] & \ldots\\
      0 \ar[r] & \coker \phi_{\alpha',L} \ar[r]  & H^1(L, \ker \phi_{\alpha'}) \ar[r]^{\iota^1_{\alpha'}}  & H^1(L,A_1 \times A_2) \ar[r]  & \ldots \\
 }$$

The homomorphism $\gamma_2^*$ is an isomorphism by Lemma \ref{lem:alpha-maps-unramified-to-unramified}. As the diagram commutes we get that $\gamma_2^*$ induces an injection $\ker \iota^1_\alpha \hookrightarrow \ker \iota^1_{\alpha'}$. Switching the roles of $\alpha$ and $\alpha'$, for which we need the liftability of $\beta_2^{-1}$, gives an injection $\ker \iota^1_{\alpha'} \hookrightarrow \ker \iota^1_\alpha$. Thus $\coker \phi_{\alpha,L}$ and $\coker \phi_{\alpha',L}$ have same cardinality. 
Now consider the dual picture, where $\gamma_2^\vee$ is an isomorphism making the diagram commutative.

$$  \xymatrix{
      0 \ar[r] & \ker \phi^\vee_{\alpha'} \ar[r]  \ar[d]_{\gamma_2^\vee} & B^\vee_{\alpha'}   \ar[r]^{\phi^\vee_{\alpha'}}  \ar[d]_{[id \times B_2]^\vee} & A_1^\vee \times A_2^\vee \ar[r] \ar[d]_{id \times B_2^\vee} & 0\\
      0 \ar[r] & \ker \phi^\vee_{\alpha} \ar[r]  & B^\vee_{\alpha} \ar[r]^{\phi^\vee_{\alpha}}  & A_1^\vee \times A_2^\vee \ar[r]  & 0\\
 }$$

With the same argument as before, one gets a bijection between $\coker \phi^\vee_{\alpha,K}$ and $\coker \phi^\vee_{\alpha',K}$ and thus they have the same number of elements. This finishes the proof of the proposition. \EP

Now we have a look at the special case of $A_1$ and $A_2$ being elliptic curves $E_1$ and $E_2$ over $\QQ$, i.e. we focus on non-simple abelian surfaces $B/\QQ$. The following setting will be used throughout the next section.

\begin{Set} \label{set:setting} Let $N$ be a positive integer and let $E_1$ and $E_2$ be two elliptic curves over $\QQ$, each having a $\QQ$-rational point $P_i$ of exact order $N$. The point $P_i$ generates a finite subgroup scheme $G_i:= \langle P_i \rangle$ in $E_i$. Denote by $E_i':=E_i/G_i$ the quotient and by $\eta_i:E_i \rightarrow E_i'$ the corresponding quotient isogeny.
Define in $E_1 \times E_2$ the finite subgroup scheme $\tilde G := \langle (P_1,nP_2)\rangle$, for some $n \in (\ZZ/N\ZZ)^*$. Let $(E_1 \times E_2)/\tilde G$ be the quotient and denote the corresponding isogeny by $\phi:E_2 \times E_2\rightarrow B$. Hence, $\phi$ is a cyclic isogeny with diagonal kernel of degree $N$. Further, $\phi$ has a $\QQ$-kernel and thus $\phi_p$ has a $\QQ_p$-kernel for every place $p$ of $\QQ$. Denote by $\eta_1 \times \eta_2 : E_1 \times E_2 \rightarrow E_1' \times E_2'$ the isogeny having as kernel $G_1 \times G_2$. We let $\psi:B\rightarrow E_1' \times E_2'$ be the isogeny satisfying $\eta_1 \times \eta_2 = \psi \circ \phi$. Note that, as elliptic curves are principally polarised, we have $E_1 \times E_2\cong (E_1 \times E_2)^\vee$ and $E_1' \times E_2' \cong (E_1' \times E_2')^\vee$. To summarise, we have a commutative diagram:
\[ \xymatrix{
&B\ar[rd] ^\psi\\
A=E_{1}\times E_{2} \ar[ur]^\varphi  \ar@/^1pc/[rr]^{\eta_1\times \eta_2}&   &A'=E_{1}'\times E_{2}'\ar[ld]^{\psi^\vee}  \ar@/^1pc/[ll]^{\eta_1^\vee \times \eta_2^\vee}\\
& B^\vee \ar[lu]^{\varphi^\vee}
}\]

By construction $\ker \eta_1 \cong \ker \eta_2 \cong \ker \phi \cong \ZZ/ N\ZZ$, therefore $\ker (\eta_1 \times \eta_2) \cong \ZZ/ N\ZZ \times \ZZ/ N\ZZ$ and $\ker \psi \cong \ZZ/ N\ZZ$. Since the kernels of the dual isogenies are the Cartier duals, we have $\ker \eta_1^\vee \cong \ker \eta_2^\vee \cong \ker \phi^\vee \cong \ker \psi^\vee \cong \mub_N$ and $\ker (\eta_1^\vee \times \eta_2 ^\vee) \cong \mub_N \times \mub_N$.  
\end{Set}

\begin{Rem}
(i) Let $G$ be a finite group scheme being isomorphic to the isomorphic group schemes $G_1$ and $G_2$, i.e. $G \cong \ZZ/N\ZZ$. Fix a generating point $P$, i.e. $G=\langle P \rangle$. Then there are natural embeddings $\iota_i$ of $G$ into $E_i$ with image $G_i$ given by $\iota_1(P):=P_1$ and $\iota_2(P):=nP_2$, such that $\tilde G$ is the embedding of $G$ into $E_1 \times E_2$ with respect to $\iota_1 \times \iota_2$. The Galois equivariant isomorphism $\alpha:G_1 \rightarrow G_2$ fulfilling the condition $\iota_2=\alpha \circ \iota_1$ is defined by $P_1 \mapsto nP_2$. In other words the choice of $n$ is equivalent to the choice of $\alpha$. As we have seen in Proposition \ref{pro:order-of-sha-independent-of-alpha}, the order of $\sha(B/\QQ)$ is indendent of that choice.

(ii) Due to Mazur's classification of possible torsion points of elliptic curves over $\QQ$, Theorem 7.5 in \cite{silvermanI} or \cite{mazur77} and \cite{mazur78}, the only possible values for $N$ in Setting \ref{set:setting} are $N=1,2,3,4,5,6,7,8,9,10,12$. 

(iii) If $\# \sha(B/\QQ)=k \cdot \square$, with $k$ square-free, then $k$ has to divide $N$. Thus, the only possible values for $k$ that one can obtain with Setting \ref{set:setting} are $k=1,2,3,5,6,7,10$. In the next section we will see that indeed all these values for $k$ are possible.
\end{Rem}

The next lemma tells us that the abelian surface $B/\QQ$ from Setting \ref{set:setting} has the interesting property that every polarisation it possesses has degree divisible by $\l$, in case $\deg \phi=N=\l$ is a prime and $E_1$ and $E_2$ are not isogenous. The proof we present follows a sketch of Brian Conrad.

\begin{lem} \label{lem:brian-conrads-construction}
Let $K$ be a field and let $E_1$ and $E_2$ be two non-isogenous elliptic curves over $K$. Let $G$ be a finite cyclic group scheme of prime order $\l$ over $K$ together with fixed embeddings $\iota_1:G \hookrightarrow E_1$ and $\iota_2:G \hookrightarrow E_2$. Thus the map $\iota_1 \times \iota_2$ is a diagonal embedding of $G$ into the product $E_1 \times E_2$. Denote its image in $E_1 \times E_2$ by $\tilde G$. Then any polarisation of the quotient $B:=(E_1 \times E_2)/\tilde G$ has degree divisible by $\l$.
\end{lem}
\proof Set $A:=E_1 \times E_2$ and let $\lambda:B \rightarrow B^\vee$ be any polarisation and consider the quotient map $\phi:A\rightarrow B$ and its dual $\phi^\vee:B^\vee \rightarrow A^\vee =A$. The composition 
$$\Psi : A \overset{\phi}{\rightarrow} B \overset{\lambda}{\rightarrow} B^\vee \overset{\phi^\vee}{\rightarrow} A$$
is a polarisation of $A$. Denote by $\emb_i:E_i \hookrightarrow A$ the natural embedding of $E_i$ into the product, and by $\pr_i:A \rightarrow E_i$ the natural projection. Define homomorphisms
$$\Psi_1 : E_1 \overset{\emb_1}{\rightarrow} A \overset{\Psi}{\rightarrow} A \overset{\pr_1}{\rightarrow} E_1 \ \ \mbox{ and } \ \ \Psi_2 : E_2 \overset{\emb_2}{\rightarrow} A \overset{\Psi}{\rightarrow} A \overset{\pr_2}{\rightarrow} E_2.$$
We claim that $\Psi =\Psi_1 \times \Psi_2$. The claim is equivalent to $\pr_2 \circ \Psi \circ \emb_1:E_1 \rightarrow E_2$ and $\pr_1 \circ \Psi \circ \emb_2:E_2 \rightarrow E_1$ being the zero map. By assumption $E_1$ and $E_2$ are non-isogenous, hence the only homomorphism between them is the zero map, which gives the claim.

Now we proceed as follows: for $i=1$ and $i=2$ we get that $\Psi_i$ is a polarisation of $E_i$ having $\iota_i(G)$ in its kernel. As the degree of a polarisation is always a square and $\l$ is a prime we get that $\l^2$ divides the degree of $\Psi_1$ and of $\Psi_2$. Therefore, $\l^4$ divides the degree of $\Psi$. We conclude that $\l^2$ divides the degree of the polarisation $\lambda$, as $\deg \Psi = \deg \phi \cdot \deg \lambda  \cdot \deg \phi^\vee  = \l^2 \cdot \deg \lambda$, which completes the proof. \EP


Now we give a remark which says that it is enough to be able to compute the Cassels-Tate equation for isogenies of prime power degree. This enables us to deal with Setting \ref{set:setting} for the composite cases $N=6$ and $N=10$. 

\begin{Rem}\label{rem:can-decompose-an-isogeny-into-prime-power-degree-isogenies}
Let $A$ and $B$ be abelian varieties over a field $K$ and let $\phi:A \rightarrow B$ be an isogeny. Denote by $\prod_i \l_i^{e_i}$ the prime factorisation of $\deg \phi$, with the $\l_i$ being pairwise different primes. The $\l_i$-primary part of the $\overline K$-rational points of $\ker \phi$ forms a Galois invariant subgroup. Hence for each $\l_i$, $\phi$ factors through an isogeny $\phi_{\l_i}: A \rightarrow B_{\l_i}$ of degree $\l_i^{e_i}$ by defining $\ker \phi_{\l_i}$ to be the subgroup scheme of $\ker \phi$ of order $\l_i^{e_i}$. Therefore, there is an isogeny $\psi_{\l_i}: B_{\l_i} \rightarrow B$ of degree coprime to $\l_i$, such that $\phi=\psi_{\l_i} \circ \phi_{\l_i}$. Thus, the $\l_i$-primary parts of $\sha(B_{\l_i}/K)$ and $\sha(B/K)$ are isomorphic. For the dual isogeny we get an analogous decomposition $\phi^\vee=\psi^\vee_{\l_i} \circ \phi^\vee_{\l_i}$. Note that $\phi^\vee_{\l_i}:= (\phi^\vee)_{\l_i} \neq (\phi_{\l_i})^\vee$. Now let $K$ be a number field. Hence, in order to compute the 
Cassels-Tate equation (\ref{equ:Cassels-Tate-equation}) for $\phi$ it suffices to compute all the Cassels-Tate equations for the $\phi_{\l_i}$. As the degrees of all $\phi_{\l_i}$ are pairwise coprime we get that
$$\coker \phi_K = \prod_i \coker \phi_{\l_i,K}, \ \ \coker \phi^\vee_K = \prod_i \coker \phi^\vee_{\l_i,K}, \ \ \coker \phi_v = \prod_i \coker \phi_{\l_i,v}.$$
The same is true for the kernels and hence we can compute
$$\frac{\#\ker \phi_K}{\#\coker \phi_K} = \prod_i \frac{\# \ker \phi_{\l_i,K}}{\# \coker \phi_{\l_i,K}}, \ \ \frac{\#\coker \phi^\vee_K}{\#\ker \phi^\vee_K} = \prod_i \frac{\# \coker \phi^\vee_{\l_i,K}}{\# \ker \phi^\vee_{\l_i,K}},$$
$$\frac{\#\coker \phi_v}{\#\ker \phi_v} = \prod_i \frac{\# \coker \phi_{\l_i,v}}{\# \ker \phi_{\l_i,v}}.$$
In case $\phi:A_1 \times A_2 \rightarrow B$ is an isogeny with diagonal kernel then all the $\phi_{\l_i}$ also have diagonal kernel.
\end{Rem}

\section[Constructing non-simple abelian surfaces over $\QQ$ with non-square order Tate-Shafarevich groups]{Constructing non-simple abelian surfaces over $\QQ$ with non-square order Tate-Shafarevich groups using elliptic curves with a rational $N$-torsion point}
\label{ch:examples}

In this section we will construct non-simple non-principally polarised abelian surfaces $B/\QQ$, such that $\#\sha(B/\QQ)=k \cdot \square$, for $k=1,2,3,5,6,7,10$. All these examples are obtained via an isogeny $\phi:E_1 \times E_2 \rightarrow B$ as constructed in Setting \ref{set:setting} with respect to $\deg \phi=N=5,6,7,10$. The elliptic curves $E_1/\QQ$ and $E_2/\QQ$ have a $\QQ$-rational $N$-torsion point, thus they correspond to points on the modular curve $X_1(N)$. The genus of $X_1(N)$ equals 0 if and only if $N=1, \ldots, 10,12$. In this case the set of $\QQ$-rational points of $X_1(N)$ is non-empty, hence there are infinitely many elliptic curves over $\QQ$ possessing a $\QQ$-rational point of order $N$ and these curves can be parametrised by a rational number $d \in \QQ$. 
The parametisations we use can be found in Proposition 1.1.2 of \cite{KloostermanMasterThesis} and Section 6 of \cite{KloSche}, also see Exercise 8.13 and Remark 7.8 of Chapter VIII of \cite{silvermanI}. The goal is to express the local and the global quotient of the Cassels-Tate equation (\ref{equ:Cassels-Tate-equation}) with respect to such a parametrisation, i.e. with respect to two rational numbers $d_1$ and $d_2$, which represent the two elliptic curves $E_1$ and $E_2$. Therefore, for fixed $N$ we will look at a two parameter family of abelian surfaces $B/\QQ$.

In the first two subsections we will compute the local and the global quotient of the Cassels-Tate equation (\ref{equ:Cassels-Tate-equation}) with respect to Setting \ref{set:setting} with a focus on $N$ being a prime number $\l$. For the local quotient we provide a formula which computes it with respect to the reduction type of $E_1$ and $E_2$ at the primes $p$. For the global quotient we explain how to obtain two functions with which one can compute the global quotient as long as a Mordell-Weil basis for $E_1$ and $E_2$ is known.

In the two prime cases $N=5$ and $N=7$, then the results of Chapter 2 enable us to give a formula computing the local and the torsion quotient for any given pair of rational numbers $(d_1,d_2)$ that correspond to the two elliptic curves via the chosen parametrisation. Further we compute the two functions to determine the global quotient once a Mordell-Weil basis of $E_1$ and $E_2$ is known. This will be discussed in the third subsection and provides examples of non-simple abelian surfaces $B$ over $\QQ$, such that $\#\sha(B/\QQ)=k \cdot \square$, for $k=5,7$. Since for any given pair $(d_1,d_2)$ we can compute whether $\#\sha(B/\QQ)$ is five or seven times a square provided we have the corresponing Mordell-Weil bases, we are able to obtain comprehensive numerical results about the occurrence of non-square order Tate-Shafarevich groups in these two families of abelian surfaces. We did so for $N=5$ and the results are presented in \cite{ANTS}.

The fourth subsection treats with the composite cases $N=6$ and $N=10$ and we will give examples of non-simple abelian surfaces $B$ over $\QQ$, such that $\#\sha(B/\QQ)=k \cdot \square$, for $k=1,2,3,6,10$.

In an appendix we will have a brief look at cyclic isogenies $\phi:E_1 \times E_2 \rightarrow B$ with diagonal kernel of degree $13$ to show that the case $\#\sha(B/\QQ)=13 \cdot \square$ is also possible.

\subsection{The local quotient}

We want to compute the quotients $\# \coker \phi_p / \# \ker \phi_p$ with respect to Setting \ref{set:setting}. If $p$ is the place at infinity this is often very easy. The induced map on $\RR$-rational points is denoted by $\phi_\infty$.

\begin{lem} \label{lem:cokernel-at-infinity-for-non-simple-abelian-surfaces}
Let $E_1$ and $E_2$ be elliptic curves over $\RR$ and $\phi:E_1 \times E_2 \rightarrow B$ a diagonal cyclic isogeny of degree $N$ having a $\RR$-kernel, i.e. $\# \ker \phi_\infty=N$. 

(i) If $2 \nmid N$, then $\coker \phi_\infty$ is trivial, and thus $\# \coker \phi_p / \# \ker \phi_p = 1/N$. 

(ii) If $2 \mid N$ assume further that both elliptic curves have negative discriminant. Then $\# \coker \phi_\infty= 2$, and thus $\# \coker \phi_p / \# \ker \phi_p = 2/N$. 
\end{lem}
\proof 
The first part follows directly from Lemma \ref{lem:size-of-cokernel-at-infinity}. For (ii) note that by the long exact sequence of Galois cohomology, $\coker \phi_\infty$ injects into $H^1(\RR, (E_1 \times E_2)[\phi])$. By assumption, the Galois action on the kernel of $\phi$ is trivial hence $H^1(\RR, (E_1 \times E_2)[\phi])$ is just the group of homomorphisms from $\ZZ/2\ZZ$ to $\ZZ/N\ZZ$, which  has $2$ elements, if $2 \mid N$. In case both discrimanants of the two elliptic curves are negative, we have that $H^1(\RR,(E_1 \times E_2)(\CC))$ is trivial, by Theorem V.2.4 in \cite{silvermanI}, which implies that $\coker \phi_\infty$ surjects onto $H^1(\RR, (E_1 \times E_2)[\phi])$, implying that it consists of two elements. \EP

Now we state the main theorem about the local quotient with respect to Setting \ref{set:setting} for $\deg \phi=N=\l$ being prime. It expresses $\# \coker \phi_p/ \# \ker \phi_p$ in terms of the type of reduction of both $E_i$ at $p$. In case the reduction type is split multiplicative we additionally have to consider whether $\ker \eta_{i,p} \subseteq (E_i)_0(\QQ_p)$, and in case the reduction type is non-split multiplicative we also have to consider the value of the Tamagawa quotient $c(E_i')_p / c(E_i)_p$.
In case the reduction type is additive, the local quotient also depends on the values of $|\eta_i'(0)|_p$. Further, we have to do some restrictions on $p$ or $\l$. If $\l \geq 5$, then the theorem determines the size of $\# \coker \phi_p / \# \ker \phi_p$ for any $p$ and any combination of reduction types of the two elliptic curves.

\begin{thm}\label{thm:cardinality-of-cokernel}
Assume Setting \ref{set:setting} for $\deg \phi=N=\l$ being prime and let $p \in M_\QQ^0$ be a finite place. Then the local quotient at $p$ can be computed as follows in case $\l \geq 5$. 
$$\frac{\# \coker \phi_p}{\# \ker \phi_p} = \begin{cases}
			    1/\l, & \text{at least one elliptic curve }  E_i \text{ has split multiplicative} \\
& \ \ \ \ \ \text{reduction at } p \text{ with }  \ker \eta_{i,p} \nsubseteq (E_i)_0(\QQ_p)\\
			    \l, & \text{both elliptic curves have split multiplicative reduction at } p \\
& \ \ \ \ \ \text{and both } \ker \eta_{i,p} \subseteq (E_i)_0(\QQ_p)\\
			    \l, & \text{both elliptic curves have additive reduction at } p \\
& \ \ \ \ \ \text{and both satisfy } |\eta_i'(0)|_p \neq 1\\
			    1, & \text{otherwise}.\\
                            \end{cases}$$
In case $\l=3$ we get the following equality.
$$\frac{\# \coker \phi_p}{\# \ker \phi_p} = \begin{cases}
			    1/3, & \text{at least one elliptic curve }  E_i \text{ has split multiplicative} \\
& \ \ \ \ \ \text{reduction at } p \text{ with }  \ker \eta_{i,p} \nsubseteq (E_i)_0(\QQ_p)\\
			    3, & \text{both elliptic curves have split multiplicative reduction at } p \\
& \ \ \ \ \ \text{and both } \ker \eta_{i,p} \subseteq (E_i)_0(\QQ_p)\\
			    1, & \text{all other cases, such that neither elliptic curve }\\
& \ \ \ \ \ \text{has additive reduction at } p.\\
                            \end{cases}$$
And in case $\l=2 \neq p$ the situation is the following.
$$\frac{\# \coker \phi_p}{\# \ker \phi_p} = \begin{cases}
			    1/2, & \text{at least one elliptic curve }  E_i \text{ has split multiplicative} \\
& \ \ \ \ \ \text{reduction at } p \text{ with }  \ker \eta_{i,p} \nsubseteq (E_i)_0(\QQ_p)\\
			    1/2, & \text{at least one elliptic curve }  E_i \text{ has non-split multiplicative} \\
& \ \ \ \ \ \text{reduction at } p \text{ with }  c(E_i')_p / c(E_i)_p=1/2,\\
			    2, & \text{both elliptic curves have either split multiplicative reduction} \\
& \ \ \ \ \ \text{at } p \text{ with } \ker \eta_{i,p} \subseteq (E_i)_0(\QQ_p) \text{ or non-split multiplicative}\\
& \ \ \ \ \ \text{reduction at } p \text{ with } c(E_i')_p / c(E_i)_p=2,\\
			    1, & \text{all other cases, such that both elliptic curves do not }\\
& \ \ \ \ \ \text{have additive reduction at } p, \text{ and } (c(E_i')_p, c(E_i)_p) \neq (2,2)\\
& \ \ \ \ \ \text{in case } E_i \text{ has non-split multiplicative reduction}.
                            \end{cases}$$
In case $\l=2 = p$ we get that $\# \coker \phi_p/\# \ker \phi_p=1/2$, if at least one elliptic curve $E_i$ has split multiplicative reduction at $p$ with $\ker \eta_{i,p} \nsubseteq (E_i)_0(\QQ_p)$.
\end{thm}
\proof Use Theorem \ref{thm:iota-p-1-of-E} and the Key Lemma \ref{lem:hitting-both-H-1-general-version} to deduce from the reduction type of both $E_i$ at $p$ plus the stated further conditions whether $\coker \phi_p$ is maximal, maximally unramified or trivial, i.e. by Corollary \ref{cor:H1} has order $\l^2$, $\l$, or $1$ respectively. As $\# \ker \phi_p=\l$ we are done. \EP 

Note, that it is not possible that one of the elliptic curves has split multiplicative reduction with $\ker \eta_{i,p} \subseteq (E_i)_0(\QQ_p)$ and the other curve has additive reduction with $|\eta_i'(0)|_p \neq 1$, as the former case implies $p \neq \l$ and the latter case implies $p=\l$.

\subsection{The global quotient}\label{sec:global-quotient}

Now we investigate the global quotient
$$\frac{\# \ker \phi_\QQ}{\# \coker \phi_\QQ}\frac{\# \coker \phi^\vee_\QQ}{\# \ker \phi^\vee_\QQ}$$
with respect to Setting \ref{set:setting} for $\deg \phi=N=\l$ being prime. As $\phi$ has a $\QQ$-kernel, the order of the kernels are clear by construction, and we need a strategy to compute the co\-kernels. We will not come up with a formula as for the local quotient, but instead we will describe a method for computing the global quotient in case one knows generators of the cokernels of $\eta_{i,\QQ}$ and $\eta_{i,\QQ}^\vee$. Clearly, one knows such generators in case one has a Mordell-Weil basis for $E_i(\QQ)$ and $E'_i(\QQ)$. 
%
%
%
By construction the maps $\ker \eta_{1,\QQ}^\vee \times \ker \eta_{2,\QQ}^\vee \rightarrow \ker \phi^\vee_\QQ$ and $\ker \eta_{1,\QQ} \times \ker \eta_{2,\QQ} \rightarrow \ker \psi_\QQ$ are surjective, therefore we have two short exact sequences of the cokernels.
$$0 \rightarrow \coker \psi^\vee_\QQ \rightarrow \coker \eta_{1,\QQ}^\vee \times \coker \eta_{2,\QQ}^\vee \rightarrow \coker \phi^\vee_\QQ \rightarrow 0$$
$$0 \rightarrow \coker \phi_\QQ \rightarrow \coker \eta_{1,\QQ} \times \coker \eta_{2,\QQ} \rightarrow \coker \psi_\QQ \rightarrow 0$$
We first compute $\coker \phi^\vee_\QQ$, which is simpler than the computation of $\coker \phi_\QQ$. We have the following long exact sequences of Galois cohomology.
$$\xymatrix{ 0 \ar[r] & \coker \eta_{1,\QQ}^\vee \times \coker \eta_{2,\QQ}^\vee \ar[r] & H^1(\QQ, (E'_1 \times E'_2)(\overline \QQ)[\eta_1^\vee \times \eta_2^\vee]) \ar[r] & \ldots }$$
$$\xymatrix{ 0 \ar[r] & \coker \phi^\vee_\QQ \ar[r] & H^1(\QQ, B^\vee(\overline \QQ)[\phi^\vee]) \ar[r] & \ldots }$$
The Kummer sequence for a number field $K$ and Hilbert's Theorem 90 yield
$$\delta_K:  H^1(K, \mub_\l) \cong K^*/K^{*\l}.$$
Since $E'_i(\overline \QQ)[\eta_i^\vee]$ and $B^\vee(\overline \QQ)[\phi^\vee]$ are isomorphic to $\mub_\l$ as Galois modules for $\Gal_\QQ$, we obtain isomorphisms from $H^1(\QQ, E'_i(\overline \QQ)[\eta_i^\vee])$ and $H^1(\QQ, B^\vee(\overline \QQ)[\phi^\vee])$ to $H^1(\QQ, \mub_\l)$. 
Composing with $\delta_\QQ$ we get natural injective group homomorphisms
$$\coker \eta_{i,\QQ}^\vee \hookrightarrow  \QQ^*/\QQ^{*\l}, \ \coker \phi^\vee_\QQ \hookrightarrow \QQ^*/\QQ^{*\l}.$$
Note that the images of these embeddings are independent of all choices made. 
We get the following commutative diagram.

\begin{equation}\label{equ:coker-phi-dual}
\begin{split}
  \xymatrix{
\coker \eta_{1,\QQ}^\vee \times \coker \eta_{2,\QQ}^\vee \ar@^{(->}[r] \ar@{->>}[d] & \QQ^*/\QQ^{*\l} \times \QQ^*/\QQ^{*\l} \ar[d] \\
\coker \phi^\vee_\QQ \ar@^{(->}[r] & \QQ^*/\QQ^{*\l}
  }
\end{split}
\end{equation}

In this diagram the natural surjection $\coker \eta_{1,\QQ}^\vee \times \coker \eta_{2,\QQ}^\vee \twoheadrightarrow \coker \phi^\vee_\QQ$ becomes $(x,y) \mapsto x^{m}/y$ as a map from $\QQ^*/\QQ^{*\l} \times \QQ^*/\QQ^{*\l}$ to $\QQ^*/\QQ^{*\l}$, for a suitable $m\in \{1,\ldots,\l-1\}$. Note that $m$ depends on $n$, but it is clear that the image of $\coker \eta_{1,\QQ}^\vee \times \coker \eta_{2,\QQ}^\vee$ in the lower right group $\QQ^*/\QQ^{*\l}$ is independent of $m$ and $n$, and for determining the image we can simply set $m=1$. The next proposition explains how to calculate the images of $\coker \eta_{i,\QQ}^\vee$ in $\QQ^*/\QQ^{*\l}$, i.e. how to calculate the upper horizontal map. Combining afterwards with $(x,y)\mapsto x/y$ gives $\coker \phi^\vee_\QQ$ as a subset of $\QQ^*/\QQ^{*\l}$. 

\begin{pro}\label{pro:dual_coker_in_K*_mod_l-th_powers}
Let $E$ and $E'$ be elliptic curves over a number field $K$ and $\eta:E \rightarrow E'$ an isogeny of prime degree $\l$. Assume that $\eta$ has a $K$-kernel, i.e. $E(\overline K)[\eta]=E(K)[\eta]$, and let $P \in E(K)$ be a generator of the kernel. Let $f_P \in K(E)$ be a $K$-rational function on $E$ such that $\div(f_P)=\l(P)-\l(\O)$. Then the following holds. 

(i) There exists a unique constant $c=c(f_P) \in K^*/K^{*\l}$ such that 
$$\coker \eta_K^\vee \rightarrow K^*/K^{*\l}$$
$$[Q] \mapsto c \cdot f_P(Q) \mod K^{*\l}, \text{ for } Q \in E(K) \text{ with } Q \neq \O, P,$$
is a well-defined and injective group homomorphism.

(ii) The image of the map $c\cdot f_P$ is independent of the choice of the point $P$ and function $f_P$ and agrees with the image of the natural injection $\coker \eta_K^\vee \hookrightarrow K^*/K^{*\l}$ described above. 

(iii) The image of the map $c\cdot f_P$ lies in the finite set
$$K(S,\l):=\{x \in K^*/K^{*\l} \mid \ v_\p(x)\equiv 0 \bmod \l, \ \forall \p \notin S\},$$
where $S$ is the set of all primes $\p \subset \O_K$, such that $\p$ divides the degree of $\eta$ or $\p$ is a prime of bad reduction of $E$.
\end{pro}
\proof This is Exercise 10.1 in \cite{silvermanI}. \EP 

\begin{Rem}
By the Riemann-Roch Theorem, the vector space of functions $f_P \in K(E)$ with $\div(f_P)=\l(P)-\l(\O)$ is $1$-dimensional, hence such a function always exists. Given such a $f_P$ it is easy to determine $c \in K^*/K^{*\l}$ and to find the value for the image of $P$ in $K^*/K^{*\l}$, by using the fact that the map $c \cdot f_P \mod K^{*\l}$ is a group homomorphism. We will do this explicitly in Propositions \ref{pro:dual_coker_for_l=5} and \ref{pro:dual_coker_for_l=7} and Lemmas \ref{lem:torsion-for-k-equal-to-6} and \ref{lem:torsion-for-k-equal-to-10}.
\end{Rem}

Now we consider the remaining case, i.e. determining $\coker \phi_\QQ$. There is no natural injection of $\coker \eta_{i,\QQ}$ into $\QQ^*/\QQ^{*\l}$ as before, since $E_i(\overline \QQ)[\eta_i]$ is not isomorphic to $\mub_\l$ as a Galois module for $\Gal_\QQ$. But $E_i(\overline \QQ)[\eta_i]$ is isomorphic to $\mub_\l$ as a Galois module for $\Gal_L$, for $L:=\QQ(\mub_\l)$. Note that the natural restriction map
$$H^1(\QQ, E_i(\overline \QQ)[\eta_i]) \rightarrow H^1(L, E_i(\overline \QQ)[\eta_i])$$
is injective, as the kernel, which equals $H^1(\Gal(L/\QQ), E_i(\overline \QQ)[\eta_i])$, is trivial, since $[L:\QQ]=\l-1$ is coprime to $\# E_i(\overline \QQ)[\eta_i]=\l$. Using the isomorphism $\delta_L$ we get natural injections $\coker \eta_{i,\QQ} \hookrightarrow H^1(\QQ, E_i(\overline \QQ)[\eta_i]) \hookrightarrow H^1(L, E_i(\overline \QQ)[\eta_i]) \cong L^*/L^{*\l}$ and hence we obtain the following commutative diagram.
\begin{equation}\label{equ:coker-phi}
\begin{split}
  \xymatrix{
\coker \phi_\QQ \ar@^{(->}[r] \ar@^{(->}[d] & L^*/L^{*\l} \ar[d]\\
\coker \eta_{1,\QQ} \times \coker \eta_{2,\QQ} \ar@^{(->}[r] \ar@{->>}[d] & L^*/L^{*\l} \times L^*/L^{*\l} \ar[d] \\
\coker \psi_\QQ \ar@^{(->}[r] & L^*/L^{*\l}
  }
\end{split}
\end{equation}
In this diagram the natural surjection $\coker \eta_{1,\QQ} \times \coker \eta_{2,\QQ} \twoheadrightarrow \coker \psi_\QQ$ is $(x,y) \mapsto x^m/y$ as a map from $L^*/L^{*\l} \times L^*/L^{*\l}$ to $L^*/L^{*\l}$, for a suitable $m\in \{1,\ldots,\l-1\}$. As before, all images are independent of $m$ and $n$, and so we can simply set $m=1$ in our computations. Hence $\coker \phi_\QQ$ is easy to determine provided we know the images of $\coker \eta_{i,\QQ}$ in $L^*/L^{*\l}$.

To obtain a map which computes the images of $\coker \eta_{i,\QQ}$ in $L^*/L^{*\l}$, we note that the dual isogeny $\eta_i^\vee:E'_i \rightarrow E_i$ has a $L$-kernel. Hence by Proposition \ref{pro:dual_coker_in_K*_mod_l-th_powers} we need a generator $\check P \in E_i'(L)$ of $E_i'[\eta_i^\vee]$ and a $L$-rational function $f_{\check P} \in L(E'_i)$, such that $\div(f_{\check P})=\l(\check P)-\l(\O)$. Again, the image of $\coker \eta_{i,\QQ}$ in $L^*/L^{*\l}$ lies in the finite set
$$L(S,\l):=\{x \in L^*/L^{*\l} \mid \ v_\p(x)\equiv 0 \bmod \l, \ \forall \p \notin S\},$$
where $S$ is the set of all primes $\p \subset \O_L$, such that $\p$ divides the degree of $\eta$ or $\p$ is a prime of bad reduction of $E_i/L$.

\subsection{$N=5$ and $N=7$ ($k=5,7$)}

Recall, that for a prime $\l\neq 2$, Mazur's theorem \cite{mazur77} tells us that the rational $\l$-primary part $E(\QQ)[\l^\infty]$ of an elliptic curve $E/\QQ$ is either trivial, or cyclic of order $\l$, where the non-trivial case can only happen if $\l=3,5$ or $7$. Further, if $E/\QQ$ has a $\QQ$-rational $5$-torsion point, then its torsion subgroup is cyclic of order $5$ or $10$ and if $E/\QQ$ has a $\QQ$-rational $7$-torsion point, then its torsion subgroup is cyclic of order $7$. We will use these standard facts without explicitly refering to them.

\subsubsection{$N=5$}

It is a well-known fact that all elliptic curves $E$ over a number field $K$ with a $K$-rational $5$-torsion point $P$ are parametrised by the Weierstra\ss \ equation
$$E: Y^2 + (d+1)XY + dY = X^3 + dX^2, \ P=(0,0),$$
for $d \in K$. Clearly the discriminant
$$\Delta = -d^5(d^2+11d-1)$$
has to be different from zero. For $K=\QQ$ this is exactly the case when $d \neq 0$ holds. Using V\'elu's algorithm \cite{velu} one can show that the curve $E$ is isogenous to the elliptic curve
$$E': Y^2 + (d+1)XY + dY = X^3 +dX^2 + (5d^3 - 10d^2 -5d)X + (d^5-10d^4 - 5d^3 - 15d^2 -d),$$
$$\Delta' = -d(d^2+11d-1)^5,$$
via the isogeny $\eta:E \rightarrow E'$ whose kernel is generated by $P$. Note that
$$\langle P \rangle =\{\O, \ P=(0,0), \ 2P=(-d, d^2), \ 3P=(-d, 0), \ 4P=(0, -d) \}.$$
If we write $d=u/v$, with $u,v \in \ZZ$ coprime, then $E$ is isomorphic to
$$E_{u,v}: Y^2 +(u+v)XY+uv^2Y=X^3+uvX^2, \ P=(0,0),$$
$$\Delta_{u,v} = -(uv)^5(u^2+11uv-v^2),$$
and $E'$ is isomorphic to
$$E_{u,v}': Y^2 +(u+v)XY+uv^2Y=$$
$$X^3+uvX^2 + (5u^3v - 10u^2v^2 -5uv^3)X + (u^5v-10u^4v^2 - 5u^3v^3 - 15u^2v^4 -uv^5),$$
$$\Delta_{u,v}' = -uv(u^2+11uv-v^2)^5,$$
$$c_{4,u,v}'=u^4 - 228 u^3 v + 494 u^2 v^2 + 228 u v^3 + v^4,$$
where $c_{4,u,v}'$ is the usual coefficients of a short Weierstra\ss \ equation of $E'_{u,v}$ as given for example in \cite[III.1]{silvermanI}.

We want to use Theorem \ref{thm:cardinality-of-cokernel} to determine the local quotient, thus for each prime $p$ we have to know the reduction type of $E$ at $p$ and the value $|\eta'(0)|_p$.

\begin{lem}\label{lem:red-type-for-5-torsion}
Let $E$ be an elliptic curve as above parametrised by $d=u/v \in \QQ \setminus \{0\}$, with $u,v \in \ZZ$ coprime, and let $p$ be a prime number.

(i) If $p|uv$ then $E$ has split multiplicative reduction at $p$ with $\ker \eta_p \nsubseteq E_0(\QQ_p)$.

(ii) If $p|u^2+11uv-v^2$ then $\ker \eta_p \subseteq E_0(\QQ_p)$. Further, $E$ has split multiplicative reduction at $p$ if and only if $p \equiv 1 \mod 5$, additive reduction if and only if $p=5$, and otherwise non-split multiplicative reduction with $p \equiv -1 \mod 5$.

(iii) a) $v_5(u^2+11uv-v^2) \in \{0,2,3\},$

b) $v_5(u^2+11uv-v^2)=0 \ \Leftrightarrow \ u \not\equiv 2v \mod 5,$

c) $v_5(u^2+11uv-v^2)=3 \ \Leftrightarrow \ u \equiv 7v \mod 25,$

d) $u \equiv 2v \mod 5\ \Rightarrow 5^4 \mid  c_{4,u,v}'$
\end{lem}
\proof Most of part (i) and (ii) follow from Lemma 1.4 and the comment thereafter of \cite{fisher__}. We give a more detailed proof. Consider the reduction-mod-$p$ map $E(\QQ_p)\rightarrow \tilde \E(\FF_p)$ and the point $P=(0,0)$, which generates $\ker \eta_p$. If $p|uv$ then $\tilde \E:\overline Y^2 + \alpha \overline{XY} = \overline X^3$, for a non-zero $\alpha \in \ZZ/p\ZZ$. In particular $\overline P$ is a node of $\tilde \E$ and the tangent cone is generated by $\overline X = -\alpha \overline Y$ and by $\overline Y=0$. Thus the reduction type is split multiplicative and $P \notin E_0(\QQ_p)$, which proves (i).
 
If $p|u^2+11uv-v^2$ then $\overline P$ is non-singular, hence $\ker \eta_p \subseteq E_0(\QQ_p)$. Also $\overline P$ is non-trivial, therefore it has order $5$. Since the order of $\overline P$ divides $\# \tilde \E_0(\FF_p)$, which equals $p-1$ if the reduction is split multiplicative, $p+1$ if the reduction is non-split multiplicative, and $p$ if the reduction is additive, we get (ii).
 
Part (iii) is an easy calculation. Note that any pair of integers $u$ and $v$ making the expression $u^2+11uv-v^2$ divisible by $5^4$ are not coprime, since $u$ and $v$ will both be divisible by $5$. \EP

\begin{pro}\label{pro:p-value-of-phi-prime-zero-is-one}
Let $\eta:E\rightarrow E'$ be the isogeny described above, for the parameter $d=u/v \in \QQ \setminus \{0\}$, with $u,v \in \ZZ$ coprime. Then $$|\eta'(0)|_p= \begin{cases}
                                                                                 1/5, & p=5 \text{ and } u \equiv 7v \mod 25\\
										 1, & \text{otherwise}.                                                                                \end{cases}$$
\end{pro}
\proof It is clear that $|\eta'(0)|_p$ equals $1$ if $p\neq 5$ or if $p$ is a place of good or multiplicative reduction; see Theorem \ref{thm:iota-p-1-of-E}. So it only remains to deal with the case $p=5$ and $p$ is additive. If $p=5$ is additive, combining Lemma \ref{lem:red-type-for-5-torsion} with \cite[Exercise 7.1]{silvermanI} gives that the Weierstra\ss \ equation for $E_{u,v}$ is minimal and the one for $E_{u,v}'$ is not minimal if and only if $u \equiv 7v \mod 25$. In this case $v_5(\Delta_{u,v}')=15$ and $c_{4,u,v}'$ is divisible at least by $5^4$, so the Weierstra\ss \ equation of $E_{u,v}'$ will become minimal if we make the following change of variables: $X \mapsto X/5^2$ and $Y \mapsto Y/5^3$. Assume that the equation for $E_{u,v}'$ is minimal. We will now compute the $p$-adic valuation of the 
leading coefficent of the power series representation of $\eta$. We claim that $\eta(Z)=Z+...$ as a power series in $Z$ in a neighbourhood of $\O$.
Set $\eta(X,Y)=:(\tilde X(X,Y), \tilde Y(X,Y))$. Then, by \cite{velu}, we have 
$- \frac{\tilde X(X,Y)}{\tilde Y(X,Y)}=\frac{p(X)}{q(X,Y)}$, for
$$p(X):=X (d + X) [d^4 + (3 d^3  + d^4) X + (3 d^2 + 3 d^3) X^2$$
$$+ (d + 3 d^2 -  d^3) X^3 + 2 d X^4 + X^5],$$
$$q(X,Y):=d^6 +( 5 d^5  + 2 d^6) X +( 10 d^4 + 8 d^5 + d^6 )X^2 +( 
 10 d^3  + 13 d^4 + 4 d^5 )X^3 $$
$$+( 5 d^2  + 10 d^3  + 
 4 d^4) X^4 +( d  + 3 d^2 + d^3 - d^4) X^5 $$
$$+Y[  2 d^5  +  (7 d^4   + d^5) X  +( 9 d^3 + 3 d^4 )X^2 + (5 d^2  + 
 3 d^3  + d^4) X^3$$
$$+ (d  - d^2  - d^3 )X^4 -  3 d X^5 - X^6].$$
For $Z:=-X/Y$, we have $X(Z)=Z^{-2}+\ldots$ and $Y(Z)=-Z^{-3}+ \ldots$ as Laurent series for $X$ and $Y$, see \cite[IV.1]{silvermanI}, therefore
$\eta(Z)=\frac{Z^{-14}+\ldots}{Z^{-15}+\ldots}=Z+\ldots $
as power series in $Z$. Hence $\eta'(0)=1$, and therefore $|\eta'(0)|_p=1$.

In case the equation for $E'_{u,v}$ was not minimal, we have to replace $Z$ by $5Z$, which gives $\eta(Z)=5Z+\ldots$, and therefore $\eta'(0)=5$. Hence $|\eta'(0)|_5=1/5$. \EP

Combining both the above lemma and proposition with Lemma \ref{lem:cokernel-at-infinity-for-non-simple-abelian-surfaces} and Theorem \ref{thm:cardinality-of-cokernel} gives complete control of the local quotient.

\begin{thm}\label{thm:cardinality-of-cokernel_l=5}
Assume Setting \ref{set:setting} with $N=5$. Let $E_i$ be given by $d_i=u_i/v_i$, for $d_i \in \QQ \setminus \{0\}$, with $u_i, v_i \in \ZZ$ coprime. If $p \in M_\QQ$ is a place, then 
$$\frac{\# \coker \phi_p}{\# \ker \phi_p} = \begin{cases}
                                            1/5, & p = \infty\\
                                            1/5, & p \mid u_1v_1u_2v_2\\
					    5, & p \mid \gcd (u_1^2+11u_1v_1-v_1^2,u_2^2+11u_2v_2-v_2^2), \ p \equiv 1 (5)\\
					    5, & u_1 \equiv 7v_1 \mod 25, u_2 \equiv 7v_2 \mod 25, \ p=5\\
					    1, & \text{otherwise}.
                                            \end{cases}$$
\end{thm}

Next comes the global quotient. As the $\eta_i$ have $\QQ$-kernels with the generating point $P_i=(0,0)$, we will use Proposition \ref{pro:dual_coker_in_K*_mod_l-th_powers} to calculate $\coker \eta_{i,\QQ}^\vee$ in $\QQ^*/\QQ^{*5}$. Note that $\coker \eta_{i,\QQ,\textnormal{tors}}^\vee$ is generated by the point $P_i$, independent of the structure of $E_{i}(\QQ)_\textnormal{tors}$.

\begin{pro}\label{pro:dual_coker_for_l=5}
For $P=(0,0)$ set 
$$f_P:=-X^2+XY+Y \in K(E).$$
The image of the natural embedding $\coker \eta^\vee_{\QQ} \hookrightarrow \QQ^*/\QQ^{*5}$ equals the image of 
$$f_P(X,Y) \mod \QQ^{*5}, \text{ for } Q=(X,Y) \neq \O, P.$$
By linearity $f_P(P)=d^4$, and $f_P(\coker \eta^\vee_{ \QQ,\textnormal{tors}})=\langle d \rangle$ in  $\QQ^*/\QQ^{*5}$.
\end{pro}
\proof For functions $X,Y,X+Y+d \in K(E)$, one easily sees that $\div(X)=(P)+(4P)-2(\O)$, $\div(Y)=2(P)+(3P)-3(\O)$, and $\div(X+Y+d)=2(3P)+(4P)-3(\O)$. Multiplying $(XY^2)/(X+Y+d)$ with $(-Y-dX)/(-Y-dX)$ yields $-X^2+XY+Y$ in $K(E)$, and thus $\div(f_P)=5(P)-5(\O)$. By Proposition \ref{pro:dual_coker_in_K*_mod_l-th_powers} we obtain that $c \cdot f_P$ is the function we are looking for. Since $(f_P(2P))^2=f_P(4P)$, we deduce $c=1$ and that $f_P(P)\equiv f_P(2P)^3\equiv d^4\bmod \QQ^{*5}$. \EP

\begin{cor}\label{cor:dual_5-torsion}
With notation as above, $E'(\QQ)[5]\cong \ZZ/5\ZZ\ \Leftrightarrow \ d \in \QQ^{*5}$.
\end{cor}
\proof We have that $E'(\QQ)[5]$ is non-trivial if and only if $\coker \eta_{\QQ}^\vee$ is trivial on the torsion part, i.e., the injective map $\eta_{\QQ, \textnormal{tors}}^\vee :E'(\QQ)\tor \rightarrow E(\QQ)\tor$ is an isomorphism. The cokernel of $\eta_{\QQ,\textnormal{tors}}^\vee$ is generated by $d$ in $\QQ^*/\QQ^{*5}$. Hence $E'(\QQ)[5]$ is non-trivial if and only if $d$ is trivial in $\QQ^*/\QQ^{*5}$. \EP
                                                                                                                                                                                                                                                                                                                                                                                                                                                                                                                                                                                                                                                                                                   
Now we calculate $\coker {\eta}_\QQ$ in $L^*/L^{*5}$, for $L:=\QQ(\xi)$, with $\xi \in \mub_5$ a primitive fifth root of unity. Fix a generator $\check P$ of  $E'(\overline \QQ)[\eta^\vee]$. Since $\check P \in E'(L)$, we have that $E'(L)[\eta^\vee]\cong \ZZ/5\ZZ$ and hence $(E', \check P)$ is isomorphic over $L$ to a pair $(E_{\tilde d}, (0,0))$, where $E_{\tilde d}$ denotes the elliptic curve over $L$ with a $L$-rational $5$-torsion point $(0,0)$ with respect to the parameter $\tilde d \in L$. Such a $L$-isomorphism $\epsilon:(E', \check P) \tilde \rightarrow (E_{\tilde d}, (0,0))$ is given by four values $r,s,t \in L$ and $w \in L^*$ and has the form $X=w^2X'+r$ and $Y=w^3Y'+w^2sX'+t$ \cite[III.1]{silvermanI}. Having such an isomorphism $\epsilon$ and the formula of $f_P$ from Proposition \ref{pro:dual_coker_for_l=5}, we can determine $f_{\check P}$, since
$$f_{\check P}(X,Y)\equiv \epsilon^* f_P(X',Y') \mod L^{*5}.$$                                                                                                                                                                                                        
To obtain $\epsilon$ we use \cite[III Table 1.2]{silvermanI}.  As $a_6$ of the Weierstra\ss \ equation of $E_{\tilde d}$ vanishes, we get $(r,t)= \check P$. The kernel polynomial of the dual isogeny $\eta^\vee:E' \rightarrow E$ is
$$X^2+(d^2+d+1)X+\frac{1}{5}(d^4-3d^3-26d^2+8d+1);$$
thus, for $\vartheta:=\xi+\xi^{-1}=(\sqrt 5 -1)/2$, we may choose
$$r= \frac{1}{5}[(-\vartheta-3)d^2 + (-11 \vartheta-8)d + (\vartheta-2)]\in \QQ(\vartheta)=\QQ(\sqrt 5),$$
$$t=\frac{1}{5}[(\xi^2 + 2\xi + 2)d^3 + (\xi^3 + 10\xi^2 + 23\xi + 11)d^2$$
$$ + (11\xi^3 - 12\xi^2 + 9\xi + 2)d + (-\xi^3 + \xi^2 - \xi + 1)] \in L.$$
Since $a_4$ of $E_{\tilde d}$ also vanishes we deduce 
$$s=\frac{1}{5}[(-4\xi^3 - 3\xi^2 - 7\xi - 6)d + (3\xi^3 - 4\xi^2 - \xi - 3)],$$
and since $a_3=a_2$ we deduce
$$w=\frac{1}{5}[(-\xi^3 - 7\xi^2 - 8\xi - 4)d + (7\xi^3 - \xi^2 + 6\xi + 3)].$$
Also one can use the conditions on the $a_i$ to calculate $\tilde d=\frac{(5\vartheta-3)d + 1}{d - (5\vartheta-3)}$. All in all we have described an algorithm to compute $f_{\check P}$. If one multiplies the obtained result by $w^5$ to get rid of denominators one obtains
$$f_{\check P}(X,Y)= \frac{1}{25} [(3 + 6\xi - \xi^2 + 7\xi^3) + (80 + 235\xi - 60\xi^2 + 245\xi^3)d$$
$$+(220 + 465\xi + 185\xi^2 + 205\xi^3)d^2 + (15 + 55\xi - 55\xi^2 + 160\xi^3)d^3$$
$$ +(140 + 280\xi + 245\xi^2 + 35\xi^3)d^4 + (-4 - 8\xi - 7\xi^2 - \xi^3)d^5] $$
$$+[(-1 + \xi - \xi^2) + (3 + 9\xi + 2\xi^2 + 2\xi^3)d + (2 + 6\xi + 8\xi^2 - 3\xi^3)d^2 $$
$$+(-1 - \xi + \xi^3)d^3]X + [(-\xi + \xi^2 - 2\xi^3) + (2 + 3\xi + 2\xi^2 + \xi^3)d]X^2 $$
$$+[(-3 - 2\xi^2 - 2\xi^3) + (-1 - 3\xi^2 - 3\xi^3)d + (-1 + 2\xi^2 + 2\xi^3)d^2]Y + XY \in L(E').$$

Now we can state the torsion quotient in terms of the pair $(d_1,d_2)$. Recall that if the two elliptic curves $E_{i}$ have rank equal to $0$, then the regulator quotient equals $1$, hence the global quotient is just the torsion quotient. If elliptic curves of positive rank are involved, we need generators for the cokernels of $\eta_{i,\QQ}$ and $\eta_{i,\QQ}^\vee$ in order to use the above described procedure to calculate the global quotient, but if both curves have trivial rank we can just use the formula given below to compute it.

\begin{pro}\label{pro:torsion-quotient-for-l=5}
Assume Setting \ref{set:setting} with $N=5$. Let $E_i$ be given by $d_i \in \QQ \setminus \{0\}$. Then the following holds.
$$\frac{\# A(\QQ)\tor \# A^\vee (\QQ)\tor}{\# B(\QQ)\tor \# B^\vee (\QQ)\tor} = \begin{cases}
                                                      1 \textnormal{ or } 5, & d_1,d_2 \in \QQ^{*5}\\
						      5^2, & d_i \in \QQ^{*5}, d_j \notin \QQ^{*5}, i \neq j\\
						      5^2, & \langle 1 \rangle \neq \langle d_1 \rangle = \langle d_2 \rangle \neq \langle 1 \rangle \textnormal{ in } \QQ^*/\QQ^{*5}\\
						      5^3, & \langle 1 \rangle \neq \langle d_1 \rangle \neq \langle d_2 \rangle \neq \langle 1 \rangle \textnormal{ in } \QQ^*/\QQ^{*5}.\\
                                                     \end{cases}$$
To be more precise, in case both $d_i\in \QQ^{*5}$, set $d_i=:\delta_i^5$, for $\delta_i \in \QQ^*$, and define $\zeta_1:=-\xi^4(\xi+1)$, $\zeta_2:=-\xi(\xi+1)$, $\zeta_3:=-\xi^3(\xi+1)$ and $\zeta_4:=-(\xi+1)$, where $\xi \in \mub_5$ is a primitive fifth root of unity. Then the torsion quotient equals $1$ if and only if
$$\left\langle \prod_{j=1}^4(\delta_1 + \zeta_j)^j(\delta_1 -1/\zeta_j)^j \right\rangle = \left\langle \prod_{j=1}^4(\delta_2 + \zeta_j)^j(\delta_2 -1/\zeta_j)^j \right\rangle \textnormal{ in } L^*/L^{*5}.$$
\end{pro}
\proof Recall that the torsion quotient equals $5 \cdot \# \coker \phi_{\QQ, \textnormal{tors}}^\vee / \# \coker \phi_{\QQ, \textnormal{tors}}$, and that $\coker \phi_{\QQ, \textnormal{tors}}^\vee$ equals the image of $\coker \eta_{1,\QQ,\textnormal{tors}}^\vee \times \coker \eta_{2,\QQ,\textnormal{tors}}^\vee$ in $\QQ^*/\QQ^{*5}$. As $\coker \eta_{i,\QQ,\textnormal{tors}}^\vee$ is generated by $d_i \mod \QQ^{*5}$ and the map on $\coker \phi_{\QQ, \textnormal{tors}}^\vee$ is $(x,y)\mapsto x/y$, we get $$\# \coker \phi_{\QQ, \textnormal{tors}}^\vee = \begin{cases} 
                                                      1, & d_1,d_2 \in \QQ^{*5}\\
						      5, & d_i \in \QQ^{*5}, d_j \notin \QQ^{*5}, i \neq j\\
						      5, & \langle 1 \rangle \neq \langle d_1 \rangle = \langle d_2 \rangle \neq \langle 1 \rangle \textnormal{ in } \QQ^*/\QQ^{*5}\\
						      5^2, & \langle 1 \rangle \neq \langle d_1 \rangle \neq \langle d_2 \rangle \neq \langle 1 \rangle \textnormal{ in } \QQ^*/\QQ^{*5}.\\
                                                     \end{cases}$$

We have seen above that $E'(\QQ)[5]\cong \ZZ/5\ZZ$ if and only if $d \in \QQ^{*5}$, hence $\coker \eta_{i,\QQ,\textnormal{tors}}$ is trivial in case $d_i \notin \QQ^{*5}$, otherwise it is $1$-dimensional. Looking at the kernel of $(x,y)\mapsto x/y$ gives
$$\# \coker \phi_{\QQ, \textnormal{tors}} = \begin{cases} 
                                                      1 \textnormal{ or } 5, & d_1,d_2 \in \QQ^{*5}\\
						      1, & \textnormal{otherwise}\\
                                                     \end{cases}$$
which finishes the first part.

For the second part, note that if $d_i=\delta_i^5$, then $E'(\QQ)[5]$ is generated by the point $P'_i=(x_i,y_i)$, where
$$x_i=\delta_i +2\delta_i^2+3\delta_i^2+5\delta_i^4+2\delta_i^5+2\delta_i^6-\delta_i^7+\delta_i^8,$$
$$y_i=\delta_i^2+3\delta_i^3+5\delta_i^4 +11\delta_i^5 +13\delta_i^6 +10\delta_i^7 + \delta_i^8 -\delta_i^{10}+\delta_i^{11} +\delta_i^{12}.$$
The image of $\langle P'_i \rangle$ under $f_{\check T}$ in $L^*/L^{*5}$, i.e. the image of $\coker \eta_{i,\QQ,\textnormal{tors}}$ in $L^*/L^{*5}$, is
$$\left\langle \prod_{j=1}^4(\delta_i + \zeta_j)^j(\delta_i -1/\zeta_j)^j \right\rangle,$$
which completes the second part. \EP

Finally, we give two unconditional examples of an abelian surface $B/\QQ$ of rank $0$, respectively of rank $1$, such that $\# \sha(B/\QQ)=5$.

\begin{Exa}
If $d_1=u_1/v_1=1/11$, $d_2=u_2/v_2=2/9$, then $\# \sha(B/\QQ)=5$.
\end{Exa}
\proof We start with the local quotient. There are three different primes dividing $u_1v_1u_2v_2=2\cdot 3^2\cdot 11$. We also have the contribution of the prime at infinity, and no contribution from any other prime, as $u_i \not \equiv 7 \cdot v_i \mod 25$ for both $i$, and $\gcd(u_1^2+11u_1v_1-v_1^2,u_2^2+11u_2v_2-v_2^2)=1$. Hence the local quotient equals $1/5^4$. Both elliptic curves $E_{i}$ have analytic rank equal to $0$, hence we know that $\sha(A/\QQ)$ and $\sha(B/\QQ)$ are finite and that the global quotient equals the torsion quotient. Thus the global quotient equals $5^3$. We conclude that $\# \sha(B/\QQ)=5\cdot \# \sha(A/\QQ)$.

It remains to show that both $\sha(E_i/\QQ)$ are trivial. The predicted size by the Birch and Swinnerton-Dyer formula is $1$. Both $E_i$ are non-CM curves of conductor $\leq 1000$, hence we can apply \cite[Theorem 3.31 and Theorem 4.4]{stein_bsd}. This gives us that $\# \sha(E_i/\QQ)[p^\infty]=1$, for all primes $p\neq5$. (The primes occurring as the degrees of cyclic isogenies or dividing any Tamagawa number are only $2$ and $5$.) Now use \cite[Theorem 1 or Table 3 in the Appendix]{fisher__} to calculate $\sel^{\eta_i}(E_i/\QQ)=0$ and $\sel^{\eta^\vee_i}(E'_i/\QQ)\cong \ZZ/5\ZZ$, for both $i$. As $\coker \eta_{i,\QQ}=0$ and $\coker \eta^\vee_{i,\QQ} \cong \ZZ/5\ZZ$ we have $\sha(E_i/\QQ)[\eta_i]=\sha(E'_i/\QQ)[\eta^\vee_i]=0$ and thus $\sha(E_i/\QQ)[5]=0$. Hence $\sha(E_i/\QQ)$ is trivial. \EP

\begin{Exa}
If $d_1=u_1/v_1=1/10$, $d_2=u_2/v_2=3/1$, then $\# \sha(B/\QQ)=5$.
\end{Exa}
\proof We have $u_1v_1u_2v_2=2\cdot3\cdot5$, $u_i \not \equiv 7 \cdot v_i \mod 25$, for both $i$, and $\gcd(u_1^2+11u_1v_1-v_1^2,u_2^2+11u_2v_2-v_2^2)=1$. Hence the local quotient equals $1/5^4$. The elliptic curve $E_1$ is of analytic rank $0$ and $E_2$ of analytic rank $1$. A generator of the free part of $E_2(\QQ)$ is the point $(-6,12)$. We will now determine $\coker \eta_{i,\QQ}^\vee$ as a subset of $\QQ^*/\QQ^{*5}$. For the first curve this equals just the torsion part of the cokernel, hence $\coker \eta_{1,\QQ}^\vee$ is generated by $\{2 \cdot 5\}$. The second cokernel is generated by the image of the torsion point, which is $3$, and by the image of $(-6,12)$ under $f=-X^2+XY+Y$, which is $-3\cdot 2^5\equiv 3 \mod \QQ^{*5}$. Therefore  $\coker \eta_{2,\QQ}^\vee$ is generated only by $\{3\}$ and hence $\coker \phi_\QQ^\vee$ has dimension equal to $2$. Since neither $d_i$ are fifth powers, we get that the dimension of $\coker \eta_{1,\QQ}$ equals $0$ and the dimension of $\coker \eta_{2,\QQ}$ equals $0$ or 
$1$, 
thus the dimension of $\coker \phi_\QQ$ equals $0$. We conclude that the global quotient equals $5^3$, which gives $\# \sha(B/\QQ)=5\cdot \# \sha(A/\QQ)$. Now one can use a similar strategy as in the previous example to show that $\sha(A/\QQ)$ is trivial. \EP

\subsubsection{$N=7$}

The situation for $N=7$ is very similar to the case $N=5$, so we mostly just state the results. The elliptic curves $E$ with a rational $7$-torsion point $P$ are parametrised by the Weierstra\ss \ equation
$$E: Y^2 + (1+d-d^2)XY + (d^2-d^3)Y = X^3 + (d^2-d^3)X^2, \ P=(0,0),$$
$$\Delta = -d^7(1-d)^7(d^3-8d^2+5d+1).$$
Thus for $K=\QQ$ we have $d\neq 0, 1$. The isogenous curve is
$$E': Y^2 + (1+d-d^2)XY + (d^2-d^3)Y = $$
$$X^3 +(d^2-d^3)X^2 + (5d-35d^2+70d^3-70d^4+35d^5-5d^7)X$$
$$ + (d-19d^2+94d^3-258d^4+393d^5-343d^6+202d^7-107d^8+46d^9-8d^{10}-d^{11}),$$
$$\Delta' = -d(1-d)(d^3-8d^2+5d+1)^7,$$
and the points in the kernel of $\eta:E \rightarrow E'$ are
$$\langle P\rangle =\{\O, \ P=(0,0), \ 2P=(d^3 - d^2, d^5 - 2d^4 + d^3), \ 3P=(d^2 - d, d^3 - 2d^2 + d),$$
$$4P=(d^2 - d, d^4 - 2d^3 + d^2), \ 5P= (d^3 - d^2, 0), \ 6P= (0, d^3 - d^2) \}.$$
If we write $d=u/v$, with $u,v \in \ZZ$ coprime, we get
$$E_{u,v}: Y^2 +((v-u)(v+u)+uv)XY+(v-u)u^2v^3Y=X^3+(v-u)u^2vX^2, \ P=(0,0),$$
$$\Delta_{u,v} = -(uv)^7(v-u)^7(u^3-8u^2v+5uv^2+v^3),$$
$$E_{u,v}': Y^2 +((v-u)(v+u)+uv)XY+(v-u)u^2v^3Y=$$
$$X^3+(v-u)u^2vX^2+(-5 u^7 v + 35 u^5 v^3 - 70 u^4 v^4 + 70 u^3 v^5 - 35 u^2 v^6 + 5 u v^7)X$$
$$-u^{11} v - 8 u^{10} v^2 + 46 u^9 v^3 - 107 u^8 v^4 + 202 u^7 v^5 - 343 u^6 v^6$$
$$ + 393 u^5 v^7 - 258 u^4 v^8 + 94 u^3 v^9 - 19 u^2 v^{10} + u v^{11},$$
$$\Delta_{u,v}' = -uv(v-u)(u^3-8u^2v+5uv^2+v^3)^7.$$
$$c_{4,u,v}'=u^8 + 228 u^7 v + 42 u^6 v^2 - 1736 u^5 v^3 + 3395 u^4 v^4 $$
$$- 3360 u^3 v^5 + 1666 u^2 v^6 - 236 u v^7 + v^8.$$

As before, to determine the local quotient we have to know the reduction type of $E$ at $p$ and the value $|\eta'(0)|_p$.

\begin{lem}\label{lem:red-type-for-7-torsion}
Let $E$ be an elliptic curve as above parametrised by $d=u/v \in \QQ \setminus \{0,1\}$, with $u,v \in \ZZ$ coprime, and let $p$ be a prime number.

(i) If $p|uv(v-u)$ then $E$ has split multiplicative reduction at $p$ with $\ker \eta_p \nsubseteq E_0(\QQ_p)$.

(ii) If $p|u^3-8u^2v+5uv^2+v^3$ then $\ker \eta_p \subseteq E_0(\QQ_p)$. Further, $E$ has split multiplicative reduction at $p$ if and only if $p \equiv 1 \mod 7$, additive reduction if and only if $p=7$, and otherwise non-split multiplicative reduction with $p \equiv -1 \mod 7$.

(iii) a) $v_7(u^3-8u^2v+5uv^2+v^3) \in \{0,2\}$,

b) $v_7(u^3-8u^2v+5uv^2+v^3)=2 \ \Leftrightarrow u \equiv 5v \mod 7,$

c) $u \equiv 5v \mod 7 \ \Rightarrow 7^6 \mid c_{4,u,v}'.$
\end{lem}
\proof Analogous to the proof of Lemma \ref{lem:red-type-for-5-torsion}. \EP
%
%

\begin{pro}\label{pro:p-value-of-phi-prime-zero-is-one_l=7}
Let $\eta:E\rightarrow E'$ be the isogeny described above, for the parameter $d=u/v \in \QQ \setminus \{0,1\}$, with $u,v \in \ZZ$ coprime. Then $$|\eta'(0)|_p= \begin{cases}
                                                                                 1/7, & p=7 \text{ and } u \equiv 5v \mod 7\\
										 1, & \text{otherwise}.                                                                                \end{cases}$$
\end{pro}
\proof Analogous to the proof of Proposition \ref{pro:p-value-of-phi-prime-zero-is-one}. \EP

Hence, for the local quotient we have the following

\begin{thm}\label{thm:cardinality-of-cokernel_l=7}
Assume Setting \ref{set:setting} with $N=7$. Let $E_i$ be given by $d_i=u_i/v_i$, for $d_i \in \QQ \setminus \{0,1\}$, with $u_i, v_i \in \ZZ$ coprime. If $p \in M_\QQ$ is a place, then
$$\frac{\# \coker \phi_p}{\# \ker \phi_p} = \begin{cases}
                                            1/7, & p =\infty\\
                                            1/7, & p \mid u_1v_1u_2v_2(v_1-u_1)(v_2-u_2)\\
					    7, & p \mid \gcd (u_1^3-8u_1^2v_1+5u_1v_1^2+v_1^3, u_2^3-8u_2^2v_2+5u_2v_2^2+v_2^3), \ p \equiv 1 (7)\\
					    7, & u_1 \equiv 5v_1 \mod 7, u_2 \equiv 5v_2 \mod 7, \ p=7\\
					    1, & \text{otherwise}.
                                            \end{cases}$$
\end{thm}

Next comes the global quotient. 

\begin{pro}\label{pro:dual_coker_for_l=7}
For $P=(0,0)$ set 
$$f_P:=d^2 X^2 + X^3 + d X^3 - d^2 Y - X Y - 2 d X Y - X^2 Y \in K(E).$$
The image of the natural embedding $\coker \eta^\vee_\QQ \hookrightarrow \QQ^*/\QQ^{*7}$ equals the image of 
$$f_P(X,Y) \mod \QQ^{*7}, \text{ for } Q=(X,Y) \neq \O, P.$$
By linearity $f_P(P)=d^3(d-1)^6$, and $f_P(\coker \eta^\vee_{\QQ,\textnormal{tors}})=\langle d(d-1)^2 \rangle$ in  $\QQ^*/\QQ^{*7}$.
\end{pro}
\proof We have that $\div(X)=(P)+(6P)-2(\O)$, $\div(Y)=2(P)+(5P)-3(\O)$, $\div(Xd-1)-Y)=(P)+2(3P)-3(\O)$, and $\div(X+Y-d^3+d^2)=(3P)+(5P)+(6P)-3(\O)$, hence $\div(Y^2X^2(X(d-1)-Y)/(X+Y-d^3+d^2)^2)=7(P)-7(\O)$. Multiplying with $(-Y-(1+d-d^2)X-(d^2-d^3))/(-Y-(1+d-d^2)X-(d^2-d^3))$ gives $d^2 X^2 + X^3 + d X^3 - d^2 Y - X Y - 2 d X Y - X^2 Y$. Proceed as in Proposition \ref{pro:dual_coker_for_l=5}. \EP

\begin{cor}\label{cor:dual_7-torsion}
With notation as above, $E'(\QQ)[7]=0$.
\end{cor}
\proof As in Corollary \ref{cor:dual_5-torsion}, $E'(\QQ)[7]$ is non-trivial if and only if $d(d-1)^2$ is trivial in $\QQ^*/\QQ^{*7}$, which is equivalent to $d$ and $d-1$ being a seventh power, for $d\in \QQ \setminus \{0,1\}$. But Fermat's Last Theorem for exponent $7$ says that this can never happen. \EP

Now set $L:=\QQ(\xi)$, for $\xi \in \mub_7$ a primitive seventh root of unity. As in case $N=5$, we want to compute a function $f_{\check P}$, which calculates the image of $\coker {\eta}_\QQ$ in $L^*/L^{*7}$, and which depends on a point $\check P=(r,t)\in E'(\overline \QQ)[\eta^\vee]$. The coefficients $r,t,s,w$ for the $L$-isomorphism $\epsilon:(E', \check P) \tilde \rightarrow (E_{\tilde d}, (0,0))$ can be computed in the same manner as before. The kernel polynomial of the dual isogeny $\eta^\vee:E' \rightarrow E$ is 
$$\frac{1}{7}(d^{12} + 3d^{11} - 51d^{10} + 185d^9 - 767d^8 + 2097d^7 - 2835d^6$$
$$ + 1738d^5 - 295d^4 - 116d^3 + 55d^2 - 15d + 1)$$
$$+ (d^8 - d^7 - 14d^6 + 32d^5 - 29d^4 + 7d^3 + 11d^2 - 7d + 1)X$$
$$ +(2d^4 - 5d^3 + 6d^2 - 3d + 2)X^2 + X^3,$$
hence for $\vartheta:=\xi + \xi^{-1}$ we may chose
$$r=\frac{1}{7}[(3\vartheta^2 + 2\vartheta -9)d^4 + (-25\vartheta^2 -19\vartheta +47)d^3 $$
$$+ (23\vartheta^2 +34\vartheta -41)d^2 + (-2\vartheta^2 -13\vartheta +6)d +(-\vartheta^2 -3\vartheta-4)] \in \QQ(\vartheta), $$
$$t=\frac{1}{7}[(-3\xi^5 - 6\xi^4 - \xi^3 - \xi^2 - 5\xi - 5)d^6 + (28\xi^5 + 59\xi^4 + 7\xi^3 + 10\xi^2 + 45\xi + 33)d^5$$
$$ + (-52\xi^5 - 119\xi^4 + 6\xi^3 - 16\xi^2 - 62\xi - 51)d^4 + (56\xi^5 + 54\xi^4 - 35\xi^3 - 37\xi^2 - 9\xi + 13)d^3 $$
$$+ (-13\xi^5 + 30\xi^4 + 54\xi^3 + 75\xi^2 + 60\xi + 32)d^2 + (-10\xi^5 - 16\xi^4 - 22\xi^3 - 25\xi^2 - 22\xi - 17)d $$
$$+ (-\xi^5 - 3\xi^4 - 5\xi^3 - 6\xi^2 - 5\xi - 1)] \in L.$$
Using the conditions on the $a_i$ gives
$$s=\frac{1}{7}[(3\xi^5 + 6\xi^4 - 5\xi^3 - 2\xi^2 + \xi + 4)d^2 + (-16\xi^5 - 11\xi^4 - 6\xi^3 - \xi^2 - 17\xi - 12)d $$
$$+ (5\xi^5 + 3\xi^4 + 8\xi^3 + 6\xi^2 + 11\xi + 2)],$$
$$w=\frac{1}{7}[(-3\xi^5 - 6\xi^4 - \xi^3 - \xi^2 - 5\xi - 5)d^6 + (28\xi^5 + 59\xi^4 + 7\xi^3 + 10\xi^2 + 45\xi + 33)d^5 $$
$$+ (-52\xi^5 - 119\xi^4 + 6\xi^3 - 16\xi^2 - 62\xi - 51)d^4 + (56\xi^5 + 54\xi^4 - 35\xi^3 - 37\xi^2 - 9\xi + 13)d^3$$
$$ + (-13\xi^5 + 30\xi^4 + 54\xi^3 + 75\xi^2 + 60\xi + 32)d^2 + (-10\xi^5 - 16\xi^4 - 22\xi^3 - 25\xi^2 - 22\xi - 17)d $$
$$+ (-\xi^5 - 3\xi^4 - 5\xi^3 - 6\xi^2 - 5\xi - 1)],$$
$$\tilde d = \frac{(\vartheta^2+3\vartheta+2)d - (\vartheta^2+3\vartheta+1)}{d - (\vartheta^2+3\vartheta+2)}.$$
Now putting everything together gives 
$$f_{\check P}\equiv w^7 \cdot f_P((X-r)/w^2,(Y-t-s(X-r))/w^3)$$
$$=w^3 \tilde d^2 (X-r)^2 + w (X-r)^3 + w \tilde d (X-r)^3 - w^4 \tilde d^2 (Y-t-s(X-r)) - w^2 (X-r)(Y-t-s(X-r)),$$
which yields a page long formula for $f_{\check P}$.

For the torsion quotient we get the following

\begin{pro}\label{pro:torsion-quotient-for-l=7}
Assume Setting \ref{set:setting} with $N=7$. Let $E_i$ be given by $d_i \in \QQ \setminus \{0,1\}$. Then
$$\frac{\# A(\QQ)\tor \# A^\vee (\QQ)\tor}{\# B(\QQ)\tor \# B^\vee (\QQ)\tor} = \begin{cases}
                                                        7^2, & \langle d_1(d_1-1)^2 \rangle = \langle d_2(d_2-1)^2 \rangle \textnormal{ in }  \QQ^*/\QQ^{*7}\\
						      7^3, & \textnormal{otherwise}.\\
                                                     \end{cases}$$
\end{pro}
\proof Since $A(\QQ)[7^\infty]\cong (\ZZ/7\ZZ)^2$ and $A'(\QQ)[7^\infty]=0$ we have $B(\QQ)[7^\infty]\cong \ZZ/7\ZZ$, and hence 
$$\# \coker \phi_{\QQ, \textnormal{tors}} = 1.$$
We know that $\coker \eta_{i,\QQ,\textnormal{tors}}^\vee$ is generated by $d_i(d_i-1)^2$ in $\QQ^*/\QQ^{*7}$ and as the product of these two cokernels maps surjectively onto $\coker \phi_{\QQ, \textnormal{tors}}^\vee $ via the map $(x,y)\mapsto x/y$, we conclude that 
$$\# \coker \phi_{\QQ, \textnormal{tors}}^\vee = \begin{cases} 
                                                      7, & \langle d_1(d_1-1)^2 \rangle = \langle d_2(d_2-1)^2 \rangle \textnormal{ in }  \QQ^*/\QQ^{*7}\\
						      7^2, & \textnormal{otherwise},\\
                                                     \end{cases}$$
which completes the proof. \EP

We finish by giving an unconditional example of an abelian surface $B/\QQ$ of rank equal to $0$, such that $\# \sha(B/\QQ)=7$.

\begin{Exa}
If $d_1=u_1/v_1=1/3$, $d_2=u_2/v_2=1/4$, then $\# \sha(B/\QQ)=7$.
\end{Exa}
\proof We have $u_1v_1u_2v_2(v_1-u_1)(v_2-u_2)=2^3\cdot3^2$, $u_1 \equiv 5 \cdot v_1 \mod 7$, $u_2 \not \equiv 5 \cdot v_2 \mod 7$, and $\gcd(u_1^3-8u_1^2v_1+5u_1v_1^2+v_1^3, u_2^3-8u_2^2v_2+5u_2v_2^2+v_2^3)=1$. Hence the local quotient equals $1/7^3$. Both elliptic curves $E_i$ have analytic rank equal to $0$, hence we know that $\sha(A/\QQ)$ and $\sha(B/\QQ)$ are finite and that the global quotient equals the torsion quotient. For $a=4$ we have that $d_1^a(d_1-1)^{2a}\equiv 2\cdot 3^2 \equiv d_2(d_2-1)^2 \mod \QQ^{*7}$, thus the global quotient equals $7^2$. We conclude that $7\cdot \# \sha(A/\QQ)=\# \sha(B/\QQ)$. As in the examples of $N=5$, one can use \cite{stein_bsd} and \cite{fisher__} to show that $\sha(A/\QQ)$ is trivial. \EP

\subsection{$N=6$ and $N=10$ ($k=1,2,3,6,10$)}

We start with the case $N=6$ and will give examples for $k=1,2,3,6$. Then we have a look at the case $N=10$ to give an example for $k=10$.

\subsubsection{$N=6$}
The elliptic curves $E$ over a number field $K$ having a rational $6$-torsion point $P$ are parametrised by
$$E: Y^2 + (d+1)XY -d(d-1)Y=X^3 -d(d-1)X^2, \ P=(0,0),$$
$$\Delta= d^6(9d - 1)(d - 1)^3,$$
for $d \in K \setminus\{0,1,1/9\}$. If $d=u/v$, with $u,v \in \ZZ$ coprime, then $E$ is isomorphic to
$$E_{u,v}: Y^2 + (u+v)XY -uv(u-v)Y =X^3 -u(u-v)X^2, \ P=(0,0),$$
$$\Delta_{u,v}= u^6 v^2(9u - v)(u - v)^3.$$

Denote by $\eta:E \rightarrow E'$ the cyclic isogeny of degree $6$, whose kernel is $\langle P \rangle$. Then
$$E': Y^2 +(d+1)XY -d(d-1)Y =X^3 -d(d-1)X^2 $$
$$-5(3d^3 - 4d^2 + d + 1)dX -(19d^5 - 33d^4 + 18d^3 - 22d^2 + 14d + 1)d.$$
The points of the kernel of $\eta$ are
$$\{ \O, \ P=(0,0), \ 2P=(d(d-1),-d^2(d-1)), \ 3P=(-d,d^2), $$
$$4P=(d(d-1),0), \ 5P=(0,d(d-1)) \}.$$
Let $\check P$ denote a generator of the kernel of the dual isogeny $\eta^\vee:E' \rightarrow E$. Then the two points of order $6$ in $\ker \eta^\vee$ are
$$\pm \check P=\left(-2d^2+4d-1, \ d^3 - \frac{1}{2}d^2 - 2d + \frac{1}{2}  \pm \frac{1}{2}(d - 1)(9d - 1)\sqrt{-3}\right).$$
The point of order $2$ in $\ker \eta^\vee$ is 
$$3\check P=\left(\frac{19}{4}d^2-\frac{14}{4}d-\frac{1}{4}, \ -\frac{19}{8}d^3-\frac{1}{8}d^2+\frac{11}{8}d+\frac{1}{8}\right)$$ 
and the two points of order $3$ are 
$$\pm 2 \check P=\left(-2d^2-2d-\frac{1}{3}, \ d^3 + \frac{5}{2}d^2 + \frac{2}{3}d + \frac{1}{6} \pm \frac{1}{18}(9d - 1)^2\sqrt{-3}\right).$$

Note that $\QQ(\mub_3)=\QQ(\sqrt{-3})$ has degree $2$ and class number 
$1$ and that the only prime which ramifies is $3$. We introduce some further notation. As $\eta$ is cyclic of degree $6$, $E$ also possesses a cyclic isogeny of degree $2$ whose kernel is generated by $3P$ and which we denote by $\eta_{\l=2}:E \rightarrow E'_{\l=2}$, and $E$ possesses a cyclic isogeny of degree $3$ whose kernel is generated by $2P$ and which we denote by $\eta_{\l=3}:E \rightarrow E'_{\l=3}$. Using V\'elu's algorithm \cite{velu} one easily computes that 
$$E'_{\l=2}:Y^2 +(d+1)XY -d(d-1)Y =X^3 -d(d-1)X^2 -5d^3X + (3d^2 + d - 1)d^3 $$
$$\Delta'_{\l=2}= d^3(9d - 1)^2(d - 1)^6.$$
Before we give the explicit examples for $k=1,2,3,6$ we provide two lemmas. One which computes the torsion quotient and one which computes the reduction type of $E$ at $p$ and further data. 

\begin{lem}\label{lem:torsion-for-k-equal-to-6}
Let $E/\QQ$ be an elliptic curves with a rational $6$-torsion point $P=(0,0)$ corresponding to the parameter $d\in \QQ \setminus\{0,1,1/9\}$ as given above. Assume that $E(\QQ)\tor = \langle P \rangle \cong \ZZ/6\ZZ$ and $E'(\QQ)\tor = \langle 3\check P \rangle \cong \ZZ/2\ZZ$. Then 

(i) $\coker \eta^\vee_{\QQ, \textnormal{tors}}$ can be identified with $\langle d \rangle \times \langle d^2(d-1)\rangle $ in $\QQ^*/\QQ^{*2} \times \QQ^*/\QQ^{*3}$, and 

(ii) $\coker \eta_{\QQ, \textnormal{tors}}$ can be identified with $\langle (9d-1)(d-1)\rangle$ in $\QQ(\mub_3)^*/\QQ(\mub_3)^{*2}$.
\end{lem}
\proof Note that by Remark \ref{rem:can-decompose-an-isogeny-into-prime-power-degree-isogenies} we have that $\coker \eta^\vee_{\QQ, \textnormal{tors}} = \coker \eta^\vee_{\l=2,\QQ, \textnormal{tors}} \times \coker \eta^\vee_{\l=3,\QQ, \textnormal{tors}}$ and also $\coker \eta_{\QQ, \textnormal{tors}} = \coker \eta_{\l=2,\QQ, \textnormal{tors}} \times \coker \eta_{\l=3,\QQ, \textnormal{tors}}$. 

Let $f_2:=X+d$ and $f_3:=Y+2dX-d^2(d-1)$ be two functions in the function field of $E/\QQ$. Then $\div(f_2)=2(3P)-2(\O)$ and $\div(f_3)=3(2P)-3(\O)$. Using Proposition \ref{pro:dual_coker_in_K*_mod_l-th_powers}, 
one easily checks that the embeddings 
$$\coker \eta^\vee_{\l=2,\QQ, \textnormal{tors}} \hookrightarrow H^1(\QQ,E'[\eta_{\l=2}^\vee]) \cong \QQ^*/\QQ^{*2} $$
$$\coker \eta^\vee_{\l=3,\QQ, \textnormal{tors}} \hookrightarrow H^1(\QQ,E'[\eta_{\l=3}^\vee]) \cong \QQ^*/\QQ^{*3}$$
are given by $f_2$, respectively $f_3$. 
As $f_2(P)\equiv d \mod \QQ^{*2}$, $f_3(P)\equiv d^2(d-1) \mod \QQ^{*3}$, $3P$ generates $\coker \eta^\vee_{\l=2,\QQ, \textnormal{tors}}$, and $2P$ generates $\coker \eta^\vee_{\l=3,\QQ, \textnormal{tors}}$, we get (i).

For (ii) note that by assumption on the torsion groups of $E$ and $E'$, we get that $\coker \eta_{\l=3,\QQ, \textnormal{tors}}$ is trivial. For $\coker \eta_{\l=2,\QQ, \textnormal{tors}}$ we have to work in the field extension $\QQ(\mub_3)$ making the action of Galois on $\ker \eta_{\l=2}^\vee$ trivial. Thus, we only need to compute the embedding
$$\coker \eta_{\l=2,\QQ, \textnormal{tors}} \hookrightarrow H^1(\QQ(\mub_3),E[\eta_{\l=2}]) \cong \QQ(\mub_3)^*/\QQ(\mub_3)^{*2}.$$
The map is given by $f_2^\vee:= X-19/4d^2+14/4d+1/4$, as $\div(f_2^\vee)=2(3\check P)-2(\O)$. From $f_2^\vee(\check P)=-3/4(d-1)(9d-1)$ it follows that $f_2^\vee(3\check P)\equiv (d-1)(9d-1) \mod \QQ(\mub_3)^{*2}$. As $3\check P$ generates $\coker \eta_{\l=2,\QQ, \textnormal{tors}}$ we get (ii). \EP

\begin{lem}\label{lem:red-type-for-k-equal-to-6}
Let $E/\QQ$ be an elliptic curves with a rational $6$-torsion point $P=(0,0)$ corresponding to the parameter $d=u/v\in \QQ \setminus\{0,1,1/9\}$, with $u,v \in \ZZ$ coprime.
Let $p$ be a prime number.

(i) If $p|u$ then $\ker \eta_{\l=2,p} \nsubseteq E_0(\QQ_p)$ and $\ker \eta_{\l=3,p} \nsubseteq E_0(\QQ_p)$ and $E$ has split multiplicative reduction at $p$.

(ii) If $p|u-v$ then $\ker \eta_{\l=2,p} \subseteq E_0(\QQ_p)$ and $\ker \eta_{\l=3,p} \nsubseteq E_0(\QQ_p)$. Further, if $p\neq 2$ then $E$ has split multiplicative reduction at $p$.

(iii) If $p|v$ then $\ker \eta_{\l=2,p} \nsubseteq E_0(\QQ_p)$ and $\ker \eta_{\l=3,p} \subseteq E_0(\QQ_p)$. Further, if $p \equiv 1 \mod 3$ then $E$ has split multiplicative reduction at $p$, and if $p \equiv 2 \mod 3$ then $E$ has non-split multiplicative reduction at $p$ with
$$c(E'_{\l=2})_p/c(E)_p=\begin{cases}
                        1/2, & v_p(v) \text{ odd}\\
                        2/2, & v_p(v) \text{ even}.
                       \end{cases}$$
                       
(iii) If $p|9u-v$, $p \neq2$, and $p\neq 3$ then $\ker \eta_{\l=2,p} \subseteq E_0(\QQ_p)$ and $\ker \eta_{\l=3,p} \subseteq E_0(\QQ_p)$. Further, if $p \equiv 1 \mod 3$ then $E$ has split multiplicative reduction at $p$, and if $p \equiv 2 \mod 3$ then $E$ has non-split multiplicative reduction at $p$ with
$$c(E'_{\l=2})_p/c(E)_p=\begin{cases}
                        2/1, & v_p(v) \text{ odd}\\
                        2/2, & v_p(v) \text{ even}.
                       \end{cases}$$
\end{lem}
\proof This is an easy exercise and will be left to the reader. First check whether $3P$ and $2P$ reduce to a non-singular or singular point to deduce whether $\ker \eta_{\l=2,p}$ and $\ker \eta_{\l=3,p}$ lie on $E_0(\QQ_p)$. For most of the statements of the lemma one can proceed as in the proof of Lemma \ref{lem:red-type-for-5-torsion}. For the statement about the Tamagawa quotient one actually has to perform Tate's algorithm \cite{Tate_algorithm} on $E$ and $E'_{\l=2}$. \EP

In the following four examples we will give two parameters $d_1=u_1/v_1$, $d_2=u_2/v_2 \in \QQ \setminus\{0,1,1/9\}$ that correspond to two elliptic curves $E_1$ and $E_2$ over $\QQ$ having a rational $6$-torsion point. Hence, $E_1$ and $E_2$ fulfill Setting \ref{set:setting} with $N=6$ and we get $\phi:E_1 \times E_2 \rightarrow B$ (with respect to some $n \in (\ZZ/6\ZZ)^*$). Recall that by Proposition \ref{pro:order-of-sha-independent-of-alpha} the order of $\sha(B/\QQ)$ is independent of the choice of $n$, thus we simply set $n=1$. Further we get the corresponding isogenies $\phi_{\l=2}$ and $\phi_{\l=3}$, which are introduced in Remark \ref{rem:can-decompose-an-isogeny-into-prime-power-degree-isogenies}.

In all four examples the analytic rank of both elliptic curves $E_1$ and $E_2$ is $0$ and the discriminant of both curves is negative. Hence, all Tate-Shafarevich groups are finite, the regulator quotient is $1$, and the local quotient at infinity for $\phi_{\l=2}$ is $1$ and for $\phi_{\l=3}$ is $1/3$ by Lemma \ref{lem:cokernel-at-infinity-for-non-simple-abelian-surfaces}. 

Also, in all four examples the rational torsion of $E_1$ and $E_2$ is isomorphic to $\ZZ/6\ZZ$ and the rational torsion of $E_1'$ and $E_2'$ is isomorphic to $\ZZ/2\ZZ$. Note that by construction $\# \ker \phi_{\QQ} / \# \ker \phi_{\QQ}^\vee = 3$. Hence we can apply the last lemma to compute the torsion quotient.   

Finally, in all four examples for both elliptic curves the reduction types at all primes $p$ are 'nice', in the sense that we can apply Theorem \ref{thm:cardinality-of-cokernel}. Thus we can compute the local quotients for $\phi_{\l=2}$ and $\phi_{\l=3}$ for all finite primes $p$.

\begin{Exa} ($k=6$) 
Choose $d_1=u_1/v_1=2/7$ and $d_2=u_2/v_2=4/17$. Then $\#\sha(B/\QQ)=6\square$.

The Cremona label of $E_1$ is $770g1$ and of $E_2$ is $8398i1$. We start with the torsion quotient. Recall that $\# \ker \phi_{\QQ} / \# \ker \phi_{\QQ}^\vee = 3$. By Lemma \ref{lem:torsion-for-k-equal-to-6}, we get 
$$\coker \eta^\vee_{1,\QQ, \textnormal{tors}} = \langle 2 \cdot 7 \rangle \times \langle 2^2 \cdot 5 \rangle, \ \coker \eta^\vee_{2,\QQ, \textnormal{tors}} = \langle 17 \rangle \times \langle 2 \cdot 13 \rangle,$$
$$\coker \eta_{1,\QQ, \textnormal{tors}} = \langle -5 \cdot 11 \rangle, \ \coker \eta_{2,\QQ, \textnormal{tors}} = \langle -13 \cdot 19 \rangle.$$ 
From Diagrams (\ref{equ:coker-phi-dual}) and (\ref{equ:coker-phi}) we conclude that 
$$\# \coker \phi^\vee_{\QQ, \textnormal{tors}}= 2^2 \cdot 3^2 \ \ \mbox{ and } \ \  \# \coker \phi_{\QQ, \textnormal{tors}} = 1.$$
Therefore, the torsion quotient equals $2^2 \cdot 3^3$.

It remains to calculate the local quotient. The conductor of $E_1$ is $2 \cdot 5 \cdot 7 \cdot 11$ and of $E_2$ is $2 \cdot 13 \cdot 17 \cdot 19$. From Lemma \ref{lem:red-type-for-k-equal-to-6} we deduce the first two rows of the following table for the finite places $p$, and with Theorem \ref{thm:cardinality-of-cokernel} we get the third and fourth row. In case the reduction type at $p$ is split multiplicative, then we indicate whether $\ker \eta_{i,\l=2,p}$ and $\ker \eta_{i,\l=3,p}$ are contained in $(E_i)_0(\QQ_p)$. In case the reduction type at $p$ is non-split multiplicative then we give the Tamagawa quotient at $p$ for $\eta{i,\l}=2$. Remember that the local quotient at infinity equals $1$ if $\l=2$, and $1/3$ if $\l=3$, by Lemma \ref{lem:cokernel-at-infinity-for-non-simple-abelian-surfaces}.
\begin{center}
    \begin{tabular}{|r||r|r|r|r|r|r|r|r|}
    \hline
    $p=$               & 2    & 5  & 7  & 11 & 13 & 17 & 19 & $\infty$\\
    \hline
red. type of $E_1$ & $\nsubseteq$, $\nsubseteq$  & $\subseteq$, $\nsubseteq$  & $\nsubseteq$, $\subseteq$  & $\frac{c'}{c}=2$ & good & good & good &\\ 
\hline
red. type of $E_2$ & $\nsubseteq$, $\nsubseteq$   & good & good & good & $\subseteq$, $\nsubseteq$  & $\frac{c'}{c}=1/2$  & $\subseteq$, $\subseteq$ &\\
\hline
$\frac{\# \coker \phi_{\l=2,p}}{\# \ker \phi_{\l=2,p}}=$  & 1/2    & 1  & 1/2  & 1  & 1  & 1/2  & 1 & 1\\
\hline
$\frac{\# \coker \phi_{\l=3,p}}{\# \ker \phi_{\l=3,p}}=$  & 1/3    & 1/3  & 1  & 1  & 1/3  & 1  & 1 & 1/3\\
\hline
         \end{tabular}
\end{center}

Hence, the local quotient equals $2^{-3} \cdot 3^{-4}$, since by Remark \ref{rem:can-decompose-an-isogeny-into-prime-power-degree-isogenies} we get that 
$$\#\coker \phi_p/\#\ker \phi_p = \# \coker \phi_{\l=2,p}/\# \ker \phi_{\l=2,p} \cdot \# \coker \phi_{\l=3,p}/\# \ker \phi_{\l=3,p}.$$
In total we have $\#\sha(B/\QQ) = 6 \cdot \#\sha(E_1 \times E_2/\QQ) = 6 \square$.
\end{Exa}

\begin{Exa} ($k=3$) 
Choose $d_1=2/7$ and $d_2=2/13$. Then $\#\sha(B/\QQ)=3\square$.

The Cremona label of $E_1$ is $770g1$ and of $E_2$ is $1430g1$ and the conductor of $E_1$ is $2 \cdot 5 \cdot 7 \cdot 11$ and of $E_2$ is $2 \cdot 5 \cdot 11 \cdot 13$. By Lemma \ref{lem:torsion-for-k-equal-to-6}, we get 
$$\coker \eta^\vee_{1,\QQ, \textnormal{tors}} = \langle 2 \cdot 7 \rangle \times \langle 2^2 \cdot 5 \rangle, \ \coker \eta^\vee_{2,\QQ, \textnormal{tors}} = \langle 2 \cdot 13 \rangle \times \langle 2^2 \cdot 11\rangle,$$
$$\coker \eta_{1,\QQ, \textnormal{tors}} = \langle -5 \cdot 11 \rangle, \ \coker \eta_{2,\QQ, \textnormal{tors}} = \langle -5 \cdot 11\rangle,$$
$$\# \coker \phi^\vee_{\QQ, \textnormal{tors}} = 2^2\cdot 3^2 \ \ \mbox{ and } \ \ \# \coker \phi_{\QQ, \textnormal{tors}} =2.$$
Hence, the torsion quotient equals $2 \cdot 3^3$. The next table follows from Lemma \ref{lem:red-type-for-k-equal-to-6} and Theorem \ref{thm:cardinality-of-cokernel} and implies that the local quotient equals $2 \cdot 3^{-4}$.
\begin{center}
    \begin{tabular}{|r||r|r|r|r|r|r|}
    \hline
    $p=$               & 2    & 5  & 7  & 11 & 13  & $\infty$ \\
    \hline
red. type of $E_1$ & $\nsubseteq$, $\nsubseteq$  & $\subseteq$, $\nsubseteq$  & $\nsubseteq$, $\subseteq$  & $\frac{c'}{c}=2$ & good  &\\ 
\hline
red. type of $E_2$ & $\nsubseteq$, $\nsubseteq$   & $\frac{c'}{c}=2$ & good & $\subseteq$, $\nsubseteq$  &  $\nsubseteq$, $\subseteq$ &\\
\hline
$\frac{\# \coker \phi_{\l=2,p}}{\# \ker \phi_{\l=2,p}}=$  & 1/2    & 2  & 1/2  & 2  & 1/2 & 1\\
\hline
$\frac{\# \coker \phi_{\l=3,p}}{\# \ker \phi_{\l=3,p}}=$  & 1/3    & 1/3  & 1  & 1/3  & 1 & 1/3\\
\hline
         \end{tabular}
\end{center}
In total we have $\#\sha(B/\QQ) = 3 \cdot \#\sha(E_1 \times E_2/\QQ) = 3 \square$.
\end{Exa}

\begin{Exa} ($k=2$) 
Choose $d_1=2/7$ and $d_2=6/7$. Then $\#\sha(B/\QQ)=2\square$.

The Cremona label of $E_1$ is $770g1$ and of $E_2$ is $1974l1$ and the conductor of $E_1$ is $2 \cdot 5 \cdot 7 \cdot 11$ and of $E_2$ is $2 \cdot 3 \cdot 7 \cdot 47$.  By Lemma \ref{lem:torsion-for-k-equal-to-6}, we get 
$$\coker \eta^\vee_{1,\QQ, \textnormal{tors}} = \langle 2 \cdot 7 \rangle \times \langle 2^2 \cdot 5 \rangle, \ \coker \eta^\vee_{2,\QQ, \textnormal{tors}} = \langle 2 \cdot 3 \cdot 7 \rangle \times \langle 2 \cdot 3 \rangle ,$$ 
$$\coker \eta_{1,\QQ, \textnormal{tors}} = \langle -5 \cdot 11 \rangle, \ \coker \eta_{2,\QQ, \textnormal{tors}} = \langle -47\rangle,$$ 
$$\#\coker \phi^\vee_{\QQ, \textnormal{tors}} = 2^2\cdot 3^2 \ \ \mbox{ and } \ \ \# \coker \phi_{\QQ, \textnormal{tors}}=1.$$
Hence, the torsion quotient equals $2^2 \cdot 3^3$. The next table follows from Lemma \ref{lem:red-type-for-k-equal-to-6} and Theorem \ref{thm:cardinality-of-cokernel} and implies that the local quotient equals $2^{-3} \cdot 3^{-3}$.
\begin{center}
    \begin{tabular}{|r|r|r|r|r|r|r|r|}
    \hline
    $p=$               & 2  & 3  & 5  & 7  & 11 & 47  & $\infty$\\
    \hline
red. type of $E_1$ & $\nsubseteq$, $\nsubseteq$ & good & $\subseteq$, $\nsubseteq$  & $\nsubseteq$, $\subseteq$  & $\frac{c'}{c}=2$ & good  &\\ 
\hline
red. type of $E_2$ & $\nsubseteq$, $\nsubseteq$   & $\nsubseteq$, $\nsubseteq$ & good & $\nsubseteq$, $\subseteq$ & good  &  $\frac{c'}{c}=2$ &\\
\hline
$\frac{\# \coker \phi_{\l=2,p}}{\# \ker \phi_{\l=2,p}}=$  & 1/2    & 1/2  & 1  & 1/2  & 1 & 1 & 1\\
\hline
$\frac{\# \coker \phi_{\l=3,p}}{\# \ker \phi_{\l=3,p}}=$  & 1/3    & 1/3  & 1/3  & 3  & 1 & 1 & 1/3\\
\hline
         \end{tabular}
\end{center}
In total we have $\#\sha(B/\QQ) = 2 \cdot \#\sha(E_1 \times E_2/\QQ) = 2 \square$.
\end{Exa}

\begin{Exa} ($k=1$) 
Choose $d_1=2/7$ and $d_2=8/13$. Then $\#\sha(B/\QQ)=\square$.

The Cremona label of $E_1$ is $770g1$ and of $E_2$ is $7670i1$ and the conductor of $E_1$ is $2 \cdot 5 \cdot 7 \cdot 11$ and of $E_2$ is $2 \cdot 5 \cdot 13 \cdot 59$. By Lemma \ref{lem:torsion-for-k-equal-to-6}, we get 
$$\coker \eta^\vee_{1,\QQ, \textnormal{tors}} = \langle 2 \cdot 7 \rangle \times \langle 2^2 \cdot 5 \rangle, \ \coker \eta^\vee_{2,\QQ, \textnormal{tors}} = \langle 2 \cdot 13 \rangle \times \langle 5\rangle ,$$ 
$$\coker \eta_{1,\QQ, \textnormal{tors}} = \langle -5 \cdot 11 \rangle, \ \coker \eta_{2,\QQ, \textnormal{tors}} = \langle -5 \cdot 59 \rangle .$$
$$\#\coker \phi^\vee_{\QQ, \textnormal{tors}}=2^2\cdot 3^2 \ \ \mbox{ and } \ \ \#\coker \phi_{\QQ, \textnormal{tors}}=1.$$ 
Hence, the torsion quotient equals $2^2 \cdot 3^3$. The next table follows from Lemma \ref{lem:red-type-for-k-equal-to-6} and Theorem \ref{thm:cardinality-of-cokernel} and implies that the local quotient equals $2^{-2} \cdot 3^{-3}$.
\begin{center}
    \begin{tabular}{|r|r|r|r|r|r|r|r|}
    \hline
    $p=$               & 2    & 5  & 7  & 11 & 13 & 59  & $\infty$\\
    \hline
red. type of $E_1$ & $\nsubseteq$, $\nsubseteq$  & $\subseteq$, $\nsubseteq$  & $\nsubseteq$, $\subseteq$  & $\frac{c'}{c}=2$ & good  & good &\\ 
\hline
red. type of $E_2$ & $\nsubseteq$, $\nsubseteq$   & $\subseteq$, $\nsubseteq$ & good & good & $\nsubseteq$, $\subseteq$  &  $\frac{c'}{c}=2$ &\\
\hline
$\frac{\# \coker \phi_{\l=2,p}}{\# \ker \phi_{\l=2,p}}=$  & 1/2    & 2  & 1/2  & 1  & 1/2 & 1 & 1\\
\hline
$\frac{\# \coker \phi_{\l=3,p}}{\# \ker \phi_{\l=3,p}}=$  & 1/3    & 1/3  & 1  & 1  & 1 & 1 & 1/3\\
\hline
         \end{tabular}
\end{center}
In total we have $\#\sha(B/\QQ) = \#\sha(E_1 \times E_2/\QQ) = \square$.
\end{Exa}

\begin{Rem}
(i) The predicted size of the Tate-Shafarevich group of all occuring elliptic curves $E_1$ and $E_2$ is $1$, hence under this assumption we have provided examples of non-simple non-principally polarized abelian surfaces over $\QQ$ such that the order of their Tate-Shafarevich groups are precisely $1,2,3,6$.

(ii) It would be interesting to investigate how often the four different values $k=1,2,3,6$ do occur while ranging through all pairs of elliptic curves having a rational $6$-torsion point, sorted by increasing height. See \cite{ANTS} for numercial tests concerning the case $N=5$.

(iii) The first curve $E_1$ is the same curve in all four examples for $N=6$, hence one could ask a related density question about the distribution of the four different values for $k$ while fixing one of the two elliptic curves and varying the other one. 
\end{Rem}

\subsubsection{$N=10$}
Finally, we will give an example for $k=10$. The elliptic curves over a number field $K$ with a rational $10$-torsion point $P$ are given by
$$E: Y^2+(-d^3+d^2+d+1)XY-d^2(d-1)(d+1)^2Y=X^3-d^2(d-1)(d+1)X^2,$$
$$P=(d^3-d,(d^3-d)^2),$$
$$\Delta=d^{10} (d - 1)^5 (d + 1)^5 (d^2 - 4d - 1)(d^2 + d - 1)^2.$$
Thus if $K=\QQ$, then $d \in \QQ \setminus \{-1,0,1\}$. 
As usual we denote the isogeny having $\langle P \rangle$ as kernel by $\eta:E \rightarrow E'$. The other nine points in the kernel of $\eta$ are
$$2P=(0,0), \ 3P=(d^2(d-1)(d+1)^2,-d^4(d-1)(d+1)^2),$$
$$4P=(d^2(d-1)(d+1),d^4(d-1)(d+1)^2),\ 5P=(-d^2,d^4),$$
$$6P=(d^2(d-1)(d+1),0), \ 7P=(d^2(d-1)(d+1)^2,d^3(d-1)^2(d+1)^3),$$
$$8P=(0,d^2(d-1)(d+1)^2), \ 9P=(d(d-1)(d+1),d(d-1)^2(d+1)).$$
The coefficients $a_1',a_2',a_3'$ for the dual curve $E'$ are the same as for $E$. The other two coefficients are
$$a_4'=-5d^{11} - 30d^{10} - 15d^9 + 40d^8 + 65d^7 - 25d^6 - 65d^5 + 40d^4 + 15d^3 - 30d^2 + 5d,$$
$$a_6'=-d^{17} - 18d^{16} - 56d^{15} - 40d^{14} + 180d^{13} + 151d^{12} - 207d^{11} - 79d^{10} + 65d^9$$
$$ - 144d^8 + 127d^7 + 221d^6 - 170d^5 - 70d^4 + 61d^3 - 18d^2 + d.$$

Let $\check P$ denote a generator of the kernel of the dual isogeny $\eta^\vee:E' \rightarrow E$. The point of order $2$ in $\ker \eta^\vee$ is
$$5\check P=(-1/4 \cdot (d^6 + 14d^5 - 5d^4 - d^2 - 14d + 1),$$
$$-1/8 \cdot (d^9 + 13d^8 - 20d^7 - 10d^6 - 14d^5 - 12d^4 + 20d^3 + 18d^2 + 13d - 1)).$$
As before, $E$ possesses a cyclic isogeny of degree $2$ whose kernel is generated by $5P$ and which we denote by $\eta_{\l=2}:E \rightarrow E'_{\l=2}$, and $E$ possesses a cyclic isogeny of degree $5$ whose kernel is generated by $2P$ and which we denote by $\eta_{\l=5}:E \rightarrow E'_{\l=5}$. 
Using V\'elu's algorithm \cite{velu} one easily computes that 
$$E'_{\l=2}:Y^2+(-d^3+d^2+d+1)XY-d^2(d-1)(d+1)^2Y=X^3-d^2(d-1)(d+1)X^2$$
$$-5d^5(d^2 + d - 1)X -d^5(d^2 + d - 1)(d^6 - 2d^5 - 5d^4 + 2d + 1),$$
$$\Delta'_{\l=2}= d^5 (d - 1)^{10}(d + 1)^{10} (d^2 - 4d - 1)^2(d^2 + d - 1).$$
Before we give the explicit example for $k=10$ we provide a lemma to compute the torsion quotient and one which computes the reduction type of $E$ at $p$ and further data. Note that the Galois extension $\QQ(\mub_5)$ has degree $4$ and class number $1$ and that the only prime that ramifies is $5$. 

\begin{lem}\label{lem:torsion-for-k-equal-to-10}
Let $E/\QQ$ be an elliptic curves with a rational $10$-torsion point $P=(d^3-d,(d^3-d)^2)$ corresponding to the parameter $d\in \QQ \setminus\{-1,0,1\}$ as above. Assume that $E'(\QQ)\tor \cong \ZZ/2\ZZ$. Then 

(i) $\coker \eta^\vee_{\QQ, \textnormal{tors}}$ can be identified with $\langle d(d^2+d-1) \rangle \times \langle d^4(d-1)(d+1)^3 \rangle$ in $\QQ^*/\QQ^{*2} \times \QQ^*/\QQ^{*5}$, and 

(ii) $\coker \eta_{\QQ, \textnormal{tors}}$ can be identified with $\langle(d - 1)(d + 1)(d^2 - 4d - 1)\rangle$ in $\QQ(\mub_5)^*/\QQ(\mub_5)^{*2}$.
\end{lem}
\proof We proceed as in Lemma \ref{lem:torsion-for-k-equal-to-6}. Let $f_2:=X+d^2$ and $f_5:=XY^2/(Y + (d+1)X - (d^5+d^4-d^3-d^2))$ be two functions in the function field of $E$. Then $\div(f_2)=2(5P)-2(\O)$ and $\div(f_5)=5(2P)-5(\O)$. Using Proposition \ref{pro:dual_coker_in_K*_mod_l-th_powers}, 
one easily checks that the embeddings 
$$\coker \eta^\vee_{\l=2,\QQ, \textnormal{tors}} \hookrightarrow H^1(\QQ,E'[\eta_{\l=2}^\vee]) \cong \QQ^*/\QQ^{*2} $$
$$\coker \eta^\vee_{\l=5,\QQ, \textnormal{tors}} \hookrightarrow H^1(\QQ,E'[\eta_{\l=5}^\vee]) \cong \QQ^*/\QQ^{*5}$$
are given by $f_2$, respectively $f_5$.  
As $f_2(P)\equiv d(d^2+d-1) \mod \QQ^{*2}$ and $f_5(P)\equiv d^4(d-1)(d+1)^3 \mod \QQ^{*5}$ and $P$ generates $\coker \eta^\vee_{\QQ, \textnormal{tors}}$ we get (i).

For (ii) note that by assumption on the torsion groups of $E$ and $E'$, we get that $\coker \eta_{\l=5,\QQ, \textnormal{tors}}$ is trivial. For $\coker \eta_{\l=2,\QQ, \textnormal{tors}}$ we have to work in the field extension $\QQ(\mub_5)$ making the action of Galois on $\ker \eta_{\l=2}^\vee$ trivial. Thus, we only need to compute the embedding
$$\coker \eta_{\l=2,\QQ, \textnormal{tors}} \hookrightarrow H^1(\QQ(\mub_5),E[\eta_{\l=2}]) \cong \QQ(\mub_5)^*/\QQ(\mub_5)^{*2}.$$
The map is given by $f_2^\vee:= X + 1/4 \cdot (d^6 + 14d^5 - 5d^4 - d^2 - 14d + 1)$, as $\div(f_2^\vee)=2(5\check P)-2(\O)$. Two of the four points of order $10$ in $\ker \eta^\vee$ have $X$-coordinate equal to 
$$(\xi^3 + \xi^2 - 1)d^6 + (-3\xi^3 - 3\xi^2)d^5 + (-7\xi^3 - 7\xi^2 - 1)d^4 + (6\xi^3 + 6\xi^2 + 3)d^3$$
$$+ (7\xi^3 + 7\xi^2 + 5)d^2 + (-3\xi^3 - 3\xi^2 - 3)d + (-\xi^3 - \xi^2 - 2),$$
where $\xi \in \mub_5$ is a primitive fifth root of unity. It follows that $f_2^\vee(5\check P)\equiv (d - 1)(d + 1)(d^2 - 4d - 1) \mod \QQ(\mub_5)^{*2}$. As $5\check P$ generates $\coker \eta_{\l=2,\QQ, \textnormal{tors}}$ we get (ii). \EP

\begin{lem}\label{lem:red-type-for-k-equal-to-10}
Let $E/\QQ$ be an elliptic curves with a rational $10$-torsion point $P=(d^3-d,(d^3-d)^2)$ corresponding to the parameter $d=u/v\in \QQ \setminus\{0,1,1/9\}$, with $u,v \in \ZZ$ coprime.
Let $p$ be a prime number.

(i) If $p|uv$ then $\ker \eta_{\l=2,p} \nsubseteq E_0(\QQ_p)$ and $\ker \eta_{\l=5,p} \nsubseteq E_0(\QQ_p)$ and $E$ has split multiplicative reduction at $p$.

(ii) If $p|(u-v)(u+v)$ then $\ker \eta_{\l=2,p} \subseteq E_0(\QQ_p)$ and $\ker \eta_{\l=5,p} \nsubseteq E_0(\QQ_p)$. Further, if $p\neq 2$ then $E$ has split multiplicative reduction at $p$.

(iii) If $p|u^2+uv-v^2$ then $\ker \eta_{\l=2,p} \nsubseteq E_0(\QQ_p)$ and $\ker \eta_{\l=5,p} \subseteq E_0(\QQ_p)$. Further for $p\neq 5$, if $uv-2v^2$ is a quadratic residue modulo $p$ then $E$ has split multiplicative reduction at $p$, and if $uv-2v^2$ is a quadratic non-residue modulo $p$ then $E$ has non-split multiplicative reduction at $p$ with
$$c(E'_{\l=2})_p/c(E)_p=\begin{cases}
                        1/2, & v_p(u^2+uv-v^2) \text{ odd}\\
                        2/2, & v_p(u^2+uv-v^2) \text{ even}.
                       \end{cases}$$
                       
(iii) If $p|u^2-4uv-v^2$ then $\ker \eta_{\l=2,p} \subseteq E_0(\QQ_p)$ and $\ker \eta_{\l=5,p} \subseteq E_0(\QQ_p)$. 
\end{lem}
\proof Computing the partial derivatives of the equation for $E$ one sees that $5P$ reduces to a singular point if and only if $-u^5v^5(u^2+uv-v^2) \equiv 0 \mod p$, and that $2P$ reduces to a singular point if and only if $-u^2v^2(u-v)(u+v) \equiv 0 \mod p$. This determines whether $\ker \eta_{\l=2,p}$ and $\ker \eta_{\l=5,p}$ lie on $E_0(\QQ_p)$. For the remaining statements use Tate's algorithm \cite{Tate_algorithm} on $E$ and $E'_{\l=2}$. For (iii) note that moving $5P$ to $(0,0)$ results in $b_2=(u^2+v^2)\left[ (u^2+uv-v^2)(u^2-3uv-v^2)-u^2v^2 \right]$, hence if $p \mid b_2$ then $p \mid u-2v$ and thus $p=5$. We leave the rest as an exercise. \EP

Now we give an unconditional example of a non-simple abelian surface $B$ over $\QQ$, such that $\#\sha(B/\QQ)=10\square$. As in all the examples of $N=6$, both elliptic curves have analytic rank equal to $0$, hence we can avoid computing the regulator quotient and get the finiteness of the Tate-Shafarevich groups. The local quotient at infinity equals $1$ for $\phi_{\l=2}$ and $1/5$ for $\phi_{\l=5}$, by Lemma \ref{lem:cokernel-at-infinity-for-non-simple-abelian-surfaces} as $E_1$ and $E_2$ have negative discriminant. Further, $E_1'(\QQ)\tor \cong E_2'(\QQ)\tor \cong \ZZ/2\ZZ$, hence we can use Lemma \ref{lem:torsion-for-k-equal-to-10} to compute the torsion quotient.

\begin{Exa} ($k=10$) 
Choose $d_1=5/2$ and $d_2=8/5$. Then $\#\sha(B/\QQ)=10\square$.

The Cremona label of $E_1$ is $123690by1$ and the conductor of $E_2$ is $338910$. By Lemma \ref{lem:torsion-for-k-equal-to-10}, we get 
$$\coker \eta^\vee_{1,\QQ, \textnormal{tors}} = \langle 2 \cdot 5 \cdot 31 \rangle \times \langle 2^2 \cdot 3 \cdot 5^4 \cdot 7^3 \rangle,$$
$$\coker \eta^\vee_{2,\QQ, \textnormal{tors}} = \langle 2 \cdot 5 \cdot 79\rangle \times \langle 2^2 \cdot 3 \cdot 5^2 \cdot 13^2\rangle ,$$
$$\coker \eta_{1,\QQ, \textnormal{tors}} = \langle - 3 \cdot 7 \cdot 19 \rangle, \ \coker \eta_{2,\QQ, \textnormal{tors}} = \langle -3  \cdot 13\rangle .$$ 
We conclude that 
$$\#\coker \phi^\vee_{\QQ, \textnormal{tors}}=2^2\cdot 5^2 \ \ \mbox{ and } \ \ \#\coker \phi_{\QQ, \textnormal{tors}}=1.$$ 
Note that by construction $\# \ker \phi_{\QQ} / \# \ker \phi_{\QQ}^\vee = 3$, thus the torsion quotient equals $2^2 \cdot 5^3$. The conductor of $E_1$ is $2 \cdot 3 \cdot 5 \cdot 7 \cdot 19 \cdot 31$ and of $E_2$ is $2 \cdot 3 \cdot 5 \cdot 11 \cdot 13 \cdot 79$. For the primes $11$ and $19$ use the equations of $E_1$, $E'_{1,\l=2}$, and $E_2$ to determine that the reduction type of $E_1$ at $p=19$ is non-split multiplicative with $c(E'_{1,\l=2})_p/c(E_1)_p=2$ and that the reduction type of $E_2$ at $p=11$ is split multiplicative. This can be done for example with Sage. The rest of the following table can be read off from Lemma \ref{lem:red-type-for-k-equal-to-10} together with Theorem \ref{thm:cardinality-of-cokernel}. Using Remark \ref{rem:can-decompose-an-isogeny-into-prime-power-degree-isogenies} it follows that the local quotient equals $2^{-3} \cdot 5^{-6}$. 
\begin{center}
    \begin{tabular}{|r||r|r|r|r|r|r|r|r|r|r|}
    \hline
    $p=$               & 2  & 3  & 5  & 7  & 11 & 13 & 19 & 31 & 79 & $\infty$\\
    \hline
red. type of $E_1$ & $\nsubseteq$, $\nsubseteq$  & $\subseteq$, $\nsubseteq$  & $\nsubseteq$, $\nsubseteq$  & $\subseteq$, $\nsubseteq$ & good  & good & $\frac{c'}{c}=2$ & $\nsubseteq$, $\subseteq$ & good &\\ 
\hline
red. type of $E_2$ & $\nsubseteq$, $\nsubseteq$   & $\subseteq$, $\nsubseteq$ & $\nsubseteq$, $\nsubseteq$ & good & $\subseteq$, $\subseteq$  &  $\subseteq$, $\nsubseteq$ & good & good & $\frac{c'}{c}=1/2$ &\\
\hline
$\frac{\# \coker \phi_{\l=2,p}}{\# \ker \phi_{\l=2,p}}=$  & 1/2    & 2  & 1/2  & 1  & 1 & 1 & 1 & 1/2 & 1/2 &1\\
\hline
$\frac{\# \coker \phi_{\l=3,p}}{\# \ker \phi_{\l=3,p}}=$  & 1/5    & 1/5  & 1/5  & 1/5  & 1 & 1/5 & 1 & 1 & 1 &1/5\\
\hline
         \end{tabular}
\end{center}
In total we have $\#\sha(B/\QQ) = 2 \cdot 5^3 \cdot \#\sha(E_1 \times E_2/\QQ) = 10 \square$.
\end{Exa} 


\subsection{Appendix. Cyclic isogenies with diagonal kernel, ($k=13$)}

We will loosen an assumption in the construction undertaken in Setting \ref{set:setting}. Instead of requiring that all points of the cyclic subgroup $G_i \subseteq E_i$ are $\QQ$-rational, we will merely demand that the $G_i$ are Galois invariant, hence instead of working with the moduli space $X_1(N)$ we will work with $X_0(N)$. Thus, we look at arbitrary cyclic isogenies $\phi:E_1 \times E_2 \rightarrow B$ with diagonal kernel. 
The next example completes the proof of Theorem \ref{thm:main-theorem}, as it shows the construction of a non-simple non-principally polarised abelian surface $B/\QQ$, such that $\#\sha(B/\QQ)=13 \cdot \square$. 

\begin{Exa}($k=13$)
Consider the following two elliptic curves over $\QQ$
$$E_1: Y^2 = X^3 - X^2 - 1829X - 32115,$$
$$E_2: Y^2 = X^3 - X^2 - 1117108895940162813412069X$$
$$- 454455515899292368353596150814715571.$$
The first curve has Cremona Label $2352j1$, where $2352=2^4 \cdot 3 \cdot 7^2$. The second curve is of conductor $135694178256= 2^4 \cdot 3 \cdot 7^2 \cdot 13 \cdot 251 \cdot 17681$. The two elliptic curves have cyclic $13$-isogenies $\eta_i:E_i \rightarrow E_i'$ with isomorphic kernels, which is due to Noam D. Elkies \cite{Elkies_13_isogeny}, as $E_1$ and $E_2$ are the quadratic twists with respect to $D=7$ of Elkies' example.

Denote by $\phi:E_1 \times E_2 \rightarrow B$ the diagonal isogeny with respect to a Galois equivariant isomorphism $\alpha: \ker \eta_1 \rightarrow \ker \eta_2$. Recall, that $\#\sha(B/\QQ)$ is independent of the choice of $\alpha$ by Proposition \ref{pro:order-of-sha-independent-of-alpha}. We claim that $\#\sha(B/\QQ)=13 \cdot \square$. 

As usual, we will compute the global and the local quotient of the Cassels-Tate equation (\ref{equ:Cassels-Tate-equation}). 
The Mordell-Weil groups of all four elliptic curves $E_1$, $E_2$, $E'_1$, and $E'_2$ are trivial and the analytic ranks are all equal to $0$. It is easy to see that this implies that the global quotient equals $1$ and we know that $\sha(B/\QQ)$ is finite. We claim that the local quotient at infinity equals $1$. As $2 \nmid \deg \phi=13$, we get that $\coker \phi_\infty$ is trivial. To prove the claim it is sufficient to show that $\ker \eta_{1,\infty}$ is trivial, too. The kernel polynomial of $\eta_1$ is 
$$(X^3 - X^2 - 1829X + 6301)(X^3 + 195X^2 + 7187X + 71569).$$
Denote by $g_1(X)$ the first factor and by $g_2(X)$ the second factor of this kernel polynomial and by $f(X):=X^3 - X^2 - 1829X - 32115$ the defining polynomial of $E_1$. All six roots of $g_1$ and $g_2$ are real numbers and both factors $g_1$ and $g_2$ generate the same totally real Galois field of degree $3$. Let $x_0$ be a zero of $g_1(X)$. As $y_0^2=f(x_0)=g_1(x_0)-38416=0-2^4 \cdot 7^4$, we get that $y_0= \pm 2^2 \cdot 7^2 \cdot \sqrt{-1} \in \CC \setminus \RR$, which shows that $\ker \eta_{1,\infty}$ is trivial. 

Among the four elliptic curves $E_i$ and $E'_i$, there are exactely two Tamagawa numbers which are divisible by $13$. These are $c(E_2)_{13}=13$ and $c(E_2')_{17681}=13$. Note that $c(E_2')_{13}=1$ and $c(E_2)_{17681}=1$. 
One easily verifies that $|\eta_i'(0)|_p=1$, for all primes $p$ and both $i$, hence by Corollary \ref{cor:criteria-for-unramified-cokernel-if-K=QQ} we conclude that $\coker \eta_{1,p}$ is maximally unramified for all primes $p$ and that $\coker \eta_{2,p}$ is maximally unramified for all $p \neq 13, 17681$.  

Using Hensel's Lemma one easily checks that $g_1(X)$ and $g_2(X)$ both factor into linear factors in $\QQ_{13}[X]$ and $\QQ_{17681}[X]$. Since $\sqrt{-1}$ also lies in $\QQ_{13}$ and in $\QQ_{17681}$ it follows that $\ker \eta_{i,13}$ and $\ker \eta_{i,17681}$ both have $13$ elements for both $i$, hence $H^1(\QQ_{13},E_i[\eta_i]) \cong H^1(\QQ_{17681},E_i[\eta_i]) \cong (\ZZ/13\ZZ)^2$ by Corollary \ref{cor:H1}. 
As $\#\coker \eta_{i,p} / \# \ker \eta_{i,p} = c(E'_i)_p / c(E_i)_p$ by Corollary \ref{cor:isogenies-on-k_v-points}, we immediately deduce that $\coker \eta_{2,13}$ is trivial and that $\coker \eta_{2,17681}$ is maximal.

Applying the Key Lemma \ref{lem:hitting-both-H-1-general-version} we deduce that the local quotient equals $1$, for all $p \neq 13, 17681$, as in this case $\coker \phi_{p}$ is maximally unramified. Further the Key Lemma implies that $\coker \phi_{13}$ is trivial and hence the local quotient for $p=13$ equals $1/13$, and that $\coker \phi_{17681}$ is maximally unramified and thus the local quotient for $p=17681$ equals $1$.

Putting everything together gives $\#\sha(B/\QQ)=13 \cdot \#\sha(E_1 \times E_2/\QQ) = 13 \cdot \square$. 
\end{Exa}

\begin{Rem}
We claim that Theorem \ref{thm:main-theorem} covers all cyclic cases, i.e. if $E_1$ and $E_2$ are elliptic curves over $\QQ$ with finite Tate-Shafarevich groups and $\phi:E_1 \times E_2 \rightarrow B$ is a cyclic isogeny, then the non-square part of $\#\sha(B/\QQ)$ equals one of the eight values $\{1,2,3,5,6,7,10,13\}$. This is ongoing work in progress.
\end{Rem}



\bibliographystyle{alpha}
\bibliography{bib_examples_of_non-square_sha}

\end{document}